\documentclass[12pt]{article}
\usepackage{full page, 
amssymb, amscd, graphicx, 
}
    \title{{\bf Logarithmic
tensor category theory, V: Convergence condition for
intertwining maps and the corresponding compatibility
condition}}
    \author{Yi-Zhi Huang, James Lepowsky and Lin Zhang}
    \date{}

\begin{document}
    \bibliographystyle{alpha}
    \maketitle

    \newtheorem{rema}{Remark}[section]
    \newtheorem{propo}[rema]{Proposition}
    \newtheorem{theo}[rema]{Theorem}
   \newtheorem{defi}[rema]{Definition}
    \newtheorem{lemma}[rema]{Lemma}
    \newtheorem{corol}[rema]{Corollary}
     \newtheorem{exam}[rema]{Example}
\newtheorem{assum}[rema]{Assumption}
     \newtheorem{nota}[rema]{Notation}
        \newcommand{\ba}{\begin{array}}
        \newcommand{\ea}{\end{array}}
        \newcommand{\be}{\begin{equation}}
        \newcommand{\ee}{\end{equation}}
        \newcommand{\bea}{\begin{eqnarray}}
        \newcommand{\eea}{\end{eqnarray}}
        \newcommand{\nno}{\nonumber}
        \newcommand{\nn}{\nonumber\\}
        \newcommand{\lbar}{\bigg\vert}
        \newcommand{\p}{\partial}
        \newcommand{\dps}{\displaystyle}
        \newcommand{\bra}{\langle}
        \newcommand{\ket}{\rangle}
 \newcommand{\res}{\mbox{\rm Res}}
\newcommand{\wt}{\mbox{\rm wt}\;}
\newcommand{\swt}{\mbox{\scriptsize\rm wt}\;}
 \newcommand{\pf}{{\it Proof}\hspace{2ex}}
 \newcommand{\epf}{\hspace{2em}$\square$}
 \newcommand{\epfv}{\hspace{1em}$\square$\vspace{1em}}
        \newcommand{\ob}{{\rm ob}\,}
        \renewcommand{\hom}{{\rm Hom}}
\newcommand{\C}{\mathbb{C}}
\newcommand{\R}{\mathbb{R}}
\newcommand{\Z}{\mathbb{Z}}
\newcommand{\N}{\mathbb{N}}
\newcommand{\A}{\mathcal{A}}
\newcommand{\Y}{\mathcal{Y}}
\newcommand{\Arg}{\mbox{\rm Arg}\;}
\newcommand{\comp}{\mathrm{COMP}}
\newcommand{\lgr}{\mathrm{LGR}}

\newcommand{\dlt}[3]{#1 ^{-1}\delta \bigg( \frac{#2 #3 }{#1 }\bigg) }

\newcommand{\dlti}[3]{#1 \delta \bigg( \frac{#2 #3 }{#1 ^{-1}}\bigg) }

 \makeatletter
\newlength{\@pxlwd} \newlength{\@rulewd} \newlength{\@pxlht}
\catcode`.=\active \catcode`B=\active \catcode`:=\active \catcode`|=\active
\def\sprite#1(#2,#3)[#4,#5]{
   \edef\@sprbox{\expandafter\@cdr\string#1\@nil @box}
   \expandafter\newsavebox\csname\@sprbox\endcsname
   \edef#1{\expandafter\usebox\csname\@sprbox\endcsname}
   \expandafter\setbox\csname\@sprbox\endcsname =\hbox\bgroup
   \vbox\bgroup
  \catcode`.=\active\catcode`B=\active\catcode`:=\active\catcode`|=\active
      \@pxlwd=#4 \divide\@pxlwd by #3 \@rulewd=\@pxlwd
      \@pxlht=#5 \divide\@pxlht by #2
      \def .{\hskip \@pxlwd \ignorespaces}
      \def B{\@ifnextchar B{\advance\@rulewd by \@pxlwd}{\vrule
         height \@pxlht width \@rulewd depth 0 pt \@rulewd=\@pxlwd}}
      \def :{\hbox\bgroup\vrule height \@pxlht width 0pt depth
0pt\ignorespaces}
      \def |{\vrule height \@pxlht width 0pt depth 0pt\egroup
         \prevdepth= -1000 pt}
   }
\def\endsprite{\egroup\egroup}
\catcode`.=12 \catcode`B=11 \catcode`:=12 \catcode`|=12\relax
\makeatother

\def\hboxtr{\FormOfHboxtr} 
\sprite{\FormOfHboxtr}(25,25)[0.5 em, 1.2 ex] 

:BBBBBBBBBBBBBBBBBBBBBBBBB |
:BB......................B |
:B.B.....................B |
:B..B....................B |
:B...B...................B |
:B....B..................B |
:B.....B.................B |
:B......B................B |
:B.......B...............B |
:B........B..............B |
:B.........B.............B |
:B..........B............B |
:B...........B...........B |
:B............B..........B |
:B.............B.........B |
:B..............B........B |
:B...............B.......B |
:B................B......B |
:B.................B.....B |
:B..................B....B |
:B...................B...B |
:B....................B..B |
:B.....................B.B |
:B......................BB |
:BBBBBBBBBBBBBBBBBBBBBBBBB |

\endsprite

\def\shboxtr{\FormOfShboxtr} 
\sprite{\FormOfShboxtr}(25,25)[0.3 em, 0.72 ex] 

:BBBBBBBBBBBBBBBBBBBBBBBBB |
:BB......................B |
:B.B.....................B |
:B..B....................B |
:B...B...................B |
:B....B..................B |
:B.....B.................B |
:B......B................B |
:B.......B...............B |
:B........B..............B |
:B.........B.............B |
:B..........B............B |
:B...........B...........B |
:B............B..........B |
:B.............B.........B |
:B..............B........B |
:B...............B.......B |
:B................B......B |
:B.................B.....B |
:B..................B....B |
:B...................B...B |
:B....................B..B |
:B.....................B.B |
:B......................BB |
:BBBBBBBBBBBBBBBBBBBBBBBBB |

\endsprite


\begin{abstract}
This is the fifth part in a series of papers in which we introduce and
develop a natural, general tensor category theory for suitable module
categories for a vertex (operator) algebra.  In this paper (Part V),
we study products and iterates of intertwining maps and of logarithmic
intertwining operators and we begin the development of our analytic
approach.
\end{abstract}


\tableofcontents
\vspace{2em}

In this paper, Part V of a series of eight papers on logarithmic
tensor category theory, we study products and iterates of intertwining
maps and of logarithmic intertwining operators and we begin the
development of our analytic approach.  The sections, equations,
theorems and so on are numbered globally in the series of papers
rather than within each paper, so that for example equation (a.b) is
the b-th labeled equation in Section a, which is contained in the
paper indicated as follows: In Part I \cite{HLZ1}, which contains
Sections 1 and 2, we give a detailed overview of our theory, state our
main results and introduce the basic objects that we shall study in
this work.  We include a brief discussion of some of the recent
applications of this theory, and also a discussion of some recent
literature.  In Part II \cite{HLZ2}, which contains Section 3, we
develop logarithmic formal calculus and study logarithmic intertwining
operators.  In Part III \cite{HLZ3}, which contains Section 4, we
introduce and study intertwining maps and tensor product bifunctors.
In Part IV \cite{HLZ4}, which contains Sections 5 and 6, we give
constructions of the $P(z)$- and $Q(z)$-tensor product bifunctors
using what we call ``compatibility conditions'' and certain other
conditions.  The present paper, Part V, contains Sections 7 and 8.  In
Part VI \cite{HLZ6}, which contains Sections 9 and 10, we construct
the appropriate natural associativity isomorphisms between triple
tensor product functors.  In Part VII \cite{HLZ7}, which contains
Section 11, we give sufficient conditions for the existence of the
associativity isomorphisms.  In Part VIII \cite{HLZ8}, which contains
Section 12, we construct braided tensor category structure.

\paragraph{Acknowledgments}
The authors gratefully
acknowledge partial support {}from NSF grants DMS-0070800 and
DMS-0401302.  Y.-Z.~H. is also grateful for partial support {}from NSF
grant PHY-0901237 and for the hospitality of Institut des Hautes 
\'{E}tudes Scientifiques in the fall of 2007.

\renewcommand{\theequation}{\thesection.\arabic{equation}}
\renewcommand{\therema}{\thesection.\arabic{rema}}
\setcounter{section}{6}
\setcounter{equation}{0}
\setcounter{rema}{0}

\section{The convergence condition for intertwining maps and
convergence and analyticity for logarithmic intertwining
operators}\label{convsec}

Now that we have constructed tensor product modules and functors, our
next goal is to construct natural associativity isomorphisms for our
category ${\cal C}$ (recall Assumption \ref{assum-c}).  More
precisely, under suitable conditions, we shall construct a natural
isomorphism between two functors {}from ${\cal C}\times{\cal
C}\times{\cal C}$ to ${\cal C}$, one given by
\[
(W_1, W_2, W_3)\mapsto (W_1\boxtimes_{P(z_1-z_{2})}
W_2)\boxtimes_{P(z_2)}W_3,
\]
and the other by
\[
(W_1, W_2, W_3)\mapsto
W_1\boxtimes_{P(z_1)} (W_2\boxtimes_{P(z_2)}W_3),
\]
where $W_1$, $W_2$ and $W_3$ are objects of ${\cal C}$ and $z_1$ and
$z_2$ are suitable complex numbers.  This will give us natural module
isomorphisms
\[
\alpha_{P(z_{1}), P(z_{2})}^{P(z_{1}-z_{2}), P(z_{2})}:
(W_1\boxtimes_{P(z_1-z_{2})} W_2)\boxtimes_{P(z_2)}W_3 \to
W_1\boxtimes_{P(z_1)} (W_2\boxtimes_{P(z_2)}W_3)
\]
and their inverses, which we will call the ``associativity
isomorphisms.''  We have seen that geometric data plays a crucial role
in the tensor product itself, and we will see that it continues to be
a crucial ingredient in the construction of the associativity
isomorphisms.

We will mainly follow, and considerably generalize, the ideas
developed in \cite{tensor4}, and it will be natural for us to work
only in the case where all tensor products involved are of type
$P(z)$, for various nonzero complex numbers $z$ (recall Remark
\ref{motivate-Mobius}).

As we have stated in Sections 4 and 5, in the remainder of this work,
in particular in this section, Assumptions \ref{assum} and
\ref{assum-c} hold.  We shall also introduce a new one, Assumption
\ref{assum-exp-set}, in this section.

In this section we study one of the prerequisites for the existence of
the associativity isomorphisms. As we have discussed in Section
1.4, in order to construct the associativity isomorphisms between
tensor products of three objects, we must have that the intertwining
maps involved are ``composable,'' which means that certain convergence
conditions have to be satisfied.  We formulate two such
conditions---one for products of suitable intertwining maps and the
other for iterates---and we prove their equivalence (Proposition
\ref{convergence}); we call the resulting single condition the
``convergence condition for intertwining maps'' (Definition
\ref{conv-conditions}).

Then we develop crucial analytic principles, including Proposition
\ref{real-exp-set} on what we call ``unique expansion sets''
(Definition \ref{uniqueexpset}), and Proposition
\ref{log-coeff-conv<=>iterate-conv} and Corollary
\ref{double-conv<=>iterate-conv} on ensuring absolute convergence of
double sums involving powers of both $z$ and $\log z$.  These
principles enable us to uniquely determine the coefficients of the
monomials in suitable variables and their logarithms obtained from
products and iterates of logarithmic intertwining operators.  We
establish the fundamental analyticity properties of suitably-evaluated
products and iterates of logarithmic intertwining operators we derive
consequences of this analyticity.

More precisely, we first need to consider the composition of a
$P(z_1)$-intertwining map and a $P(z_2)$-intertwining map for suitable
nonzero complex numbers $z_1$ and $z_2$.  Geometrically, these
compositions correspond to sewing operations (see \cite{H0} and \cite{H1}) of
Riemann surfaces with punctures and local coordinates.  Compositions
(that is, products and iterates) of maps of this type have been
defined in \cite{tensor4} for intertwining maps among ordinary
modules.  The same definitions carry over to the greater generality of
this work:

Recall {}from Definition \ref{im:imdef} the space
\[
\mathcal{M}[P(z)]_{W_{1}W_{2}}^{W_{3}}
\]
of $P(z)$-intertwining maps of type ${W_3}\choose {W_1W_{2}}$ for 
$z\in \C^{\times}$ and $W_{1}, W_{2}, W_{3}$ objects of $\mathcal{C}$.
Let $W_1$, $W_2$, $W_3$, $W_4$ and $M_1$ be objects of ${\cal C}$.
Let $z_1,z_2 \in {\mathbb C}^{\times}$, 
$I_1\in \mathcal{M}[P(z_{1})]_{W_{1}M_{1}}^{W_{4}}$ and 
$I_2\in \mathcal{M}[P(z_{2})]_{W_{2}W_{3}}^{M_{1}}$.
If for any $w_{(1)}\in W_1$,
$w_{(2)}\in W_2$, $w_{(3)}\in W_3$ and $w'_{(4)}\in W'_4$, the series
\begin{equation}\label{convp}
\sum_{n\in {\mathbb C}}\langle w'_{(4)}, I_1(w_{(1)}\otimes
\pi_n(I_2(w_{(2)}\otimes w_{(3)})))\rangle_{W_4}
\end{equation}
(recall the notation $\pi_n$ {}from (\ref{pi_n}) and Definition
\ref{Wbardef} and note that
$\pi_n(I_2(w_{(2)}\otimes w_{(3)}))\in M_1$) is absolutely 
convergent, then the sums of these series give a linear map
\[
W_1\otimes W_2\otimes W_3 \to (W'_4)^{*}.
\]
Recalling the arguments in Lemmas \ref{4.36} and
\ref{IlambdatoJlambda}, we see that the image of this map is actually
in $\overline{W_{4}}$, so that we obtain a linear map
\[
W_1\otimes W_2\otimes W_3 \to \overline{W_{4}}.
\]
Analogously, let $W_1$, $W_2$, $W_3$, $W_4$ and $M_2$ be objects of
${\cal C}$. let $z_2,z_0 \in {\mathbb C}^{\times}$, 
$I^1\in \mathcal{M}[P(z_{2})]_{M_{2}W_{3}}^{W_{4}}$ and 
$I^2\in \mathcal{M}[P(z_{0})]_{W_{1}W_{2}}^{M_{2}}$.
If for any
$w_{(1)}\in W_1$, $w_{(2)}\in W_2$, $w_{(3)}\in W_3$ and $w'_{(4)}\in
W'_4$, the series
\begin{equation}\label{convi}
\sum_{n\in {\mathbb C}}\langle w'_{(4)}, I^1(\pi_n
(I^2(w_{(1)}\otimes w_{(2)}))\otimes w_{(3)})\rangle_{W_4}
\end{equation}
is absolutely convergent, then the sums of these series also 
give a linear map
\[
W_1\otimes W_2\otimes W_3 \to \overline{W_{4}}.
\]

\begin{defi}\label{productanditerateexisting}{\rm Let $W_1$, $W_2$, $W_3$, 
$W_4$ and $M_1$ be objects of ${\cal C}$.  Let $z_1,z_2 \in {\mathbb
C}^{\times}$, $I_1\in \mathcal{M}[P(z_{1})]_{W_{1}M_{1}}^{W_{4}}$ and
$I_2\in \mathcal{M}[P(z_{2})]_{W_{2}W_{3}}^{M_{1}}$.  We say that {\it
the product of $I_1$ and $I_2$ exists} if for any $w_{(1)}\in W_1$,
$w_{(2)}\in W_2$, $w_{(3)}\in W_3$ and $w'_{(4)}\in W'_4$, the series
(\ref{convp}) is absolutely convergent. In this case, we denote the
sum (\ref{convp}) by
\begin{equation}\label{I-prod}
\langle w'_{(4)}, I_1(w_{(1)}\otimes
I_2(w_{(2)}\otimes w_{(3)}))\rangle.
\end{equation}
We call
the map
\[
W_1\otimes W_2\otimes W_3 \to \overline{W}_4,
\]
defined by (\ref{I-prod})
the {\it product} of $I_1$ and $I_2$ and denote
it by 
\[
I_1 \circ (1_{W_{1}}\otimes I_{2}).
\]
In particular, we have
\[
\langle w'_{(4)}, I_1(w_{(1)}\otimes
I_2(w_{(2)}\otimes w_{(3)}))\rangle=\langle w'_{(4)},
(I_1 \circ (1_{W_{1}}\otimes I_{2}))(w_{(1)}\otimes w_{(2)}\otimes
w_{(3)})\rangle.
\]
Analogously, let $W_1$, $W_2$, $W_3$, $W_4$ and $M_2$ be objects of
${\cal C}$, and let $z_2,z_0 \in {\mathbb C}^{\times}$, 
$I^1\in \mathcal{M}[P(z_{2})]_{M_{2}W_{3}}^{W_{4}}$ and 
$I^2\in \mathcal{M}[P(z_{0})]_{W_{1}W_{2}}^{M_{2}}$.
We say that {\it the iterate of $I^1$ and $I^2$ exists} if for any
$w_{(1)}\in W_1$, $w_{(2)}\in W_2$, $w_{(3)}\in W_3$ and $w'_{(4)}\in
W'_4$, the series (\ref{convi}) is absolutely convergent. In this
case, we denote the sum (\ref{convi}) by
\begin{equation}\label{I-iter}
\langle w'_{(4)},
I^1(I^2(w_{(1)}\otimes w_{(2)})\otimes w_{(3)})\rangle
\end{equation}
and we call the map
\[
W_1\otimes W_2\otimes W_3\to
\overline{W}_4
\]
defined by (\ref{I-iter})
the {\it iterate} of
$I^1$ and $I^2$ and denote
it by 
\[
I^1\circ
(I^2\otimes 1_{W_3}).
\]
In particular, we have
\[
\langle w'_{(4)},
I^1(I^2(w_{(1)}\otimes w_{(2)})\otimes w_{(3)})\rangle
=\langle w'_{(4)}, (I^1\circ
(I^2\otimes 1_{W_3}))(w_{(1)}\otimes w_{(2)}\otimes w_{(3)})\rangle.
\]  }
\end{defi}

\begin{rema}\label{grad-comp-prod-iter}
{\rm Note that {}from the grading compatibility condition
(\ref{grad-comp}) for $P(z)$-intertwining maps, the product and the
iterate defined above, when they exist, also satisfy the following
{\it grading compatibility conditions}: With the notation as in
Definition \ref{productanditerateexisting}, suppose that $w_{(1)}\in
W_1^{(\beta)}$, $w_{(2)}\in W_2^{(\gamma)}$ and $w_{(3)}\in
W_3^{(\delta)}$, where $\beta, \gamma, \delta \in \tilde A$.  Then
\[
(I_1 \circ (1_{W_{1}}\otimes I_{2}))(w_{(1)}\otimes 
w_{(2)}\otimes w_{(3)})\in \overline{W_{4}^{(\beta+\gamma+\delta)}}
\]
if the product of $I_1$ and $I_2$ exists, and 
\[
(I^1\circ
(I^2\otimes 1_{W_3}))(w_{(1)}\otimes
w_{(2)}\otimes w_{(3)})\in \overline{W_{4}^{(\beta+\gamma+\delta)}}
\]
if the iterate of $I^1$ and $I^2$ exists.}
\end{rema}

\begin{propo}\label{convergence}
The following two conditions are equivalent:
\begin{enumerate}
\item Let $W_1$, $W_2$, $W_3$, $W_4$ and $M_1$ be arbitrary
objects of ${\cal
C}$ and let $z_1$ and $z_2$ be arbitrary nonzero complex numbers satisfying
\[
|z_1|>|z_2|>0.
\] 
Then for any $I_1\in \mathcal{M}[P(z_{1})]_{W_{1}M_{1}}^{W_{4}}$ and
$I_2\in \mathcal{M}[P(z_{2})]_{W_{2}W_{3}}^{M_{1}}$, the product of
$I_1$ and $I_2$ exists.

\item Let $W_1$, $W_2$, $W_3$, $W_4$ and $M_2$ be arbitrary
objects of ${\cal
C}$ and let  $z_0$ and $z_2$ be arbitrary nonzero complex numbers satisfying
\[
|z_2|>|z_0|>0.
\]
Then for any $I^1\in \mathcal{M}[P(z_{2})]_{M_{2}W_{3}}^{W_{4}}$ and
$I^2\in \mathcal{M}[P(z_{0})]_{W_{1}W_{2}}^{M_{2}}$, the iterate of
$I^1$ and $I^2$ exists.
\end{enumerate}
\end{propo}
\pf We shall use the isomorphism $\Omega_0$ given by (\ref{Omega_r})
and its inverse $\Omega_{-1}$ (recall Proposition \ref{log:omega}) to
prove this result.  Suppose that Condition 1 holds.  Let $z_0$ and
$z_2$ be any nonzero complex numbers. For any intertwining maps $I^1$
and $I^2$ as in the statement of Condition 2, let 
${\cal Y}^1=\Y_{I^{1}, 0}$ and
${\cal Y}^2=\Y_{I^{2}, 0}$ be the logarithmic intertwining 
operators corresponding
to $I^1$ and $I^2$, respectively, according to Proposition
\ref{im:correspond}.  We need to prove that when $|z_2|>|z_0|>0$, the
series (\ref{convi}), which can now be written as
\begin{equation}\label{4itm}
\sum_{n\in \C}\left(\langle w'_{(4)}, 
{\cal Y}^1(\pi_{n}({\cal Y}^2(w_{(1)}, x_0)w_{(2)}),
x_2)w_{(3)} \rangle_{W_4}\lbar_{x_0=z_0, \;x_2=z_2}\right)
\end{equation}
(recall the ``substitution'' notation {}from (\ref{im:f(z)}), where we
choose $p=0$ for both substitutions), 
is absolutely convergent for any $w_{(1)}\in
W_1$, $w_{(2)}\in W_2$, $w_{(3)}\in W_3$ and $w'_{(4)}\in W'_4$.

Using the linear isomorphism $\Omega_{-1}: {\cal V}_{W_3 M_2}^{W_4}
\to {\cal V}_{M_2 W_3}^{W_4}$ (see (\ref{Omega_r})),
\[
\Omega_{-1}({\cal Y})(w, x)w_{(3)}=e^{xL(-1)}{\cal Y}(w_{(3)}, e^{-\pi
i}x)w,
\]
for ${\cal Y}\in {\cal V}_{W_3\,
M_2}^{W_4}$, $w\in M_2$ and 
$w_{(3)}\in W_3$, and its inverse $\Omega_0:{\cal V}_{M_2 W_3}^{W_4}\to
{\cal V}_{W_3 M_2}^{W_4}$, we have
\begin{eqnarray}\label{nosub}
\lefteqn{\langle w'_{(4)}, {\cal Y}^1({\cal Y}^2(w_{(1)},
x_0)w_{(2)}, x_2)w_{(3)} \rangle_{W_4}}\nno\\
&&=\langle w'_{(4)}, \Omega_{-1}(\Omega_0({\cal Y}^1))({\cal
Y}^2(w_{(1)}, x_0)w_{(2)}, x_2)w_{(3)} \rangle_{W_4} \nno\\
&&=\langle w'_{(4)}, e^{x_2L(-1)}\Omega_0({\cal Y}^1) (w_{(3)},
e^{-\pi i}x_2) {\cal Y}^2(w_{(1)}, x_0)w_{(2)}\rangle_{W_4}
\nno\\
&&=\langle e^{x_2L'(1)}w'_{(4)}, \Omega_0({\cal Y}^1) (w_{(3)},
e^{-\pi i}x_2) {\cal Y}^2(w_{(1)}, x_0)w_{(2)}\rangle_{W_4}
\end{eqnarray}
for $w_{(1)}\in W_1$, $w_{(2)}\in W_2$, $w_{(3)}\in W_3$ and
$w'_{(4)}\in W'_4$. Hence for $n\in \C$, 
\begin{eqnarray*}
\lefteqn{\langle w'_{(4)}, {\cal Y}^1(\pi_{n}({\cal Y}^2(w_{(1)},
x_0)w_{(2)}), x_2)w_{(3)} \rangle_{W_4}\lbar_{x_0= z_0,\;
x_2=z_2}}\nno\\
&&=\langle e^{x_2L'(1)}w'_{(4)}, \Omega_0({\cal Y}^1) (w_{(3)},
e^{-\pi i}x_2) \pi_{n}({\cal Y}^2(w_{(1)},
x_0)w_{(2)})\rangle_{W_4}\lbar_{x_0= z_0,\; x_2=z_2} \nno\\
&&=\langle e^{z_2L'(1)}w'_{(4)}, \Omega_0({\cal Y}^1)
(w_{(3)}, x_2) \pi_{n}({\cal Y}^2(w_{(1)},
x_0)w_{(2)})\rangle_{W_4}\lbar_{x_0= z_0,\; x_2=
-z_2},\nno\\ &&
\end{eqnarray*}
where in the last expression we take $p=0$ (respectively, $p=-1$) in
(\ref{im:f(z)}) and (\ref{log:IYp}) for the substitution $x_2=-z_2$
when $\pi\le \arg z_2 <2 \pi$, in which case $\log (-z_2)=\log z_2-\pi
i$ (respectively, when $0\le \arg z_2 < \pi$, in which case $\log
(-z_2)=\log z_2+\pi i$); cf. the corresponding considerations in
Example \ref{expl-wv}.  For brevity, let us write this last expression
as
\[
\langle e^{z_2L'(1)}w'_{(4)}, \Omega_0({\cal Y}^1)
(w_{(3)}, x_2) \pi_{n}({\cal Y}^2(w_{(1)},
x_0)w_{(2)})\rangle_{W_4}\lbar_{x_0= z_0,\; x_2=e^{-\pi i} 
z_2};
\]
that is, the substitution $x_2=e^{-\pi i} z_2$ refers to the indicated
procedure, which amounts to substituting
\[
e^{\log z_2 - \pi i}
\]
for $x_2$.  Thus
\begin{eqnarray}\label{i2p}
\lefteqn{\sum_{n\in \C}\left(\langle w'_{(4)}, {\cal Y}^1(\pi_{n}({\cal Y}^2(w_{(1)},
x_0)w_{(2)}), x_2)w_{(3)} \rangle_{W_4}\lbar_{x_0= z_0,\;
x_2=z_2}\right)}\nn
&&=\sum_{n\in \C}\left(\langle e^{z_2L'(1)}w'_{(4)}, \Omega_0({\cal Y}^1)
(w_{(3)}, x_2) \pi_{n}({\cal Y}^2(w_{(1)},
x_0)w_{(2)})\rangle_{W_4}\lbar_{x_0= z_0,\; x_2=e^{-\pi i}
z_2}\right).\nno\\ &&
\end{eqnarray}
Since the last expression is equal to the product of a
$P(-z_2)$-intertwining map and a $P(z_0)$-intertwining map evaluated
at $w_{(3)}\otimes w_{(1)}\otimes w_{(2)}\in W_3\otimes W_1\otimes
W_2$ and paired with $e^{z_2L'(1)}w'_{(4)}\in W'_4$, it converges
absolutely when $|-z_2|>|z_0|>0$, or equivalently, when
$|z_2|>|z_0|>0$.

Conversely, suppose that Condition 2 holds, and let $z_1$ and $z_2$ be
any nonzero complex numbers. For any intertwining maps $I_1$ and $I_2$
as in the statement of Condition 1, let ${\cal Y}_1=\Y_{I_{1}, 0}$ and
${\cal Y}_2=\Y_{I_{2}, 0}$ be the logarithmic intertwining operators
corresponding to $I_1$ and $I_2$, respectively.  We need to prove that
when $|z_1|>|z_2|>0$, the series (\ref{convp}), which can now be
written as
\begin{equation}\label{4prm}
\sum_{n\in \C}\left(\langle w'_{(4)}, {\cal Y}_1(w_{(1)}, x_1) 
\pi_{n}({\cal Y}_2(w_{(2)},
x_2)w_{(3)})\rangle_{W_4}\lbar_{x_1=z_1, \; x_2=z_2}\right),
\end{equation}
is absolutely convergent for any $w_{(1)}\in W_1$, $w_{(2)}\in W_2$,
$w_{(3)}\in W_3$ and $w'_{(4)}\in W'_4$.

Using the linear isomorphism $\Omega_0: {\cal V}_{M_1 W_1}^{W_4}
\to {\cal V}_{W_1 M_1}^{W_4}$,
\[
\Omega_0({\cal Y})(w_{(1)}, x)w=e^{xL(-1)}{\cal Y}(w, e^{\pi
i}x)w_{(1)},
\]
for  ${\cal Y}\in {\cal V}_{M_1\,
W_1}^{W_4}$, $w_{(1)}\in W_1$ and $w\in M_1$, 
and its inverse $\Omega_{-1}:{\cal V}_{W_1 M_1}^{W_4}\to
{\cal V}_{M_1 W_1}^{W_4}$, we have
\begin{eqnarray}\label{nosub2}
\lefteqn{\langle w'_{(4)}, {\cal Y}_1 (w_{(1)},
x_1) {\cal Y}_2(w_{(2)}, x_2)w_{(3)}\rangle_{W_4}}\nno\\
&&=\langle w'_{(4)}, \Omega_0(\Omega_{-1}({\cal Y}_1)) (w_{(1)},
x_1) {\cal Y}_2(w_{(2)}, x_2)w_{(3)}\rangle_{W_4}\nno\\
&&=\langle w'_{(4)}, e^{x_1L(-1)}\Omega_{-1}({\cal Y}_1)({\cal
Y}_2(w_{(2)}, x_2)w_{(3)}, e^{\pi i}x_1)w_{(1)} \rangle_{W_4} \nno\\
&&=\langle e^{x_1L'(1)}w'_{(4)}, \Omega_{-1}({\cal Y}_1)({\cal
Y}_2(w_{(2)},
x_2)w_{(3)}, e^{\pi i}x_1)w_{(1)} \rangle_{W_4}
\end{eqnarray}
for $w_{(1)}\in W_1$, $w_{(2)}\in W_2$, $w_{(3)}\in W_3$ and
$w'_{(4)}\in W'_4$. Hence for $n\in \C$,
\begin{eqnarray*}
\lefteqn{\langle w'_{(4)}, {\cal Y}_1 (w_{(1)},
x_1) \pi_{n}({\cal Y}_2(w_{(2)}, x_2)w_{(3)})\rangle_{W_4}\lbar_{x_1= z_1,\;
x_2=z_2}}\nno\\
&&=\langle e^{x_1L'(1)}w'_{(4)}, \Omega_{-1}({\cal Y}_1)(
\pi_{n}({\cal
Y}_2(w_{(2)},
x_2)w_{(3)}), e^{\pi i}x_1)w_{(1)} \rangle_{W_4}\lbar_{x_1= z_1,\;
x_2=z_2}\nno\\
&&=\langle e^{z_1L'(1)}w'_{(4)}, \Omega_{-1}({\cal Y}_1)(
\pi_{n}({\cal
Y}_2(w_{(2)},
x_2)w_{(3)}), x_1)w_{(1)} \rangle_{W_4}\lbar_{x_1= e^{\pi i}z_1,\;
x_2=z_2},\nno\\ &&
\end{eqnarray*}
where the substitution $x_1=e^{\pi i} z_1$ is interpreted as above,
namely, we substitute
\[
e^{\log z_1 + \pi i}
\]
for $x_1$; here $p=0$ (respectively, $p=1$) when $0\le \arg z_1 < \pi$
(respectively, when $\pi\le \arg z_1 <2 \pi$) (cf. above).  Thus
\begin{eqnarray}\label{p2i}
\lefteqn{\sum_{n\in \C}\left(\langle w'_{(4)}, {\cal Y}_1 (w_{(1)},
x_1) \pi_{n}({\cal Y}_2(w_{(2)}, x_2)w_{(3)})\rangle_{W_4}\lbar_{x_1= z_1,\;
x_2=z_2}\right)}\nn
&&=\sum_{n\in \C}\left(\langle e^{z_1L'(1)}w'_{(4)}, \Omega_{-1}({\cal Y}_1)(
\pi_{n}({\cal
Y}_2(w_{(2)},
x_2)w_{(3)}), x_1)w_{(1)} \rangle_{W_4}\lbar_{x_1= e^{\pi i}z_1,\;
x_2=z_2}\right).\nno\\ &&
\end{eqnarray}
Since the last expression is equal to the iterate of a
$P(-z_1)$-intertwining map and a $P(z_2)$-intertwining map
evaluated at $w_{(2)}\otimes w_{(3)}\otimes w_{(1)}\in W_2\otimes
W_3\otimes W_1$ and paired with $e^{z_1L'(1)}w'_{(4)}\in W'_4$,
it converges absolutely when $|-z_1|>|z_2|>0$, or equivalently,
when $|z_1|>|z_2|>0$.
\epfv

For convenience, we shall use the notations 
\begin{equation}\label{iter-abbr-pq}
\langle w'_{(4)}, 
{\cal Y}^1({\cal Y}^2(w_{(1)}, x_0)w_{(2)},
x_2)w_{(3)} \rangle_{W_4}\lbar_{x_0^{n}=e^{nl_{p}(z_0)},\;
\log x_{0}=l_{p}(z_{0}),\; \;x_2^{n}=e^{nl_{q}(z_2)},
\; \log x_{2}=l_{q}(z_{2})}
\end{equation}
and 
\begin{equation}\label{prod-abbr-pq}
\langle w'_{(4)}, {\cal Y}_1(w_{(1)}, x_1) 
{\cal Y}_2(w_{(2)},
x_2)w_{(3)}\rangle_{W_4}\lbar_{x_1^{n}=e^{nl_{p}(z_1)},\;
\log x_{1}=l_{p}(z_{1}),\; \;x_2^{n}=e^{nl_{q}(z_2)},
\; \log x_{2}=l_{q}(z_{2})}
\end{equation}
to denote
\[
\sum_{n\in \C}\left(\langle w'_{(4)}, 
{\cal Y}^1(\pi_{n}({\cal Y}^2(w_{(1)}, x_0)w_{(2)}),
x_2)w_{(3)} \rangle_{W_4}\lbar_{x_0^{n}=e^{nl_{p}(z_0)},\;
\log x_{0}=l_{p}(z_{0}),\; \;x_2^{n}=e^{nl_{q}(z_2)},
\; \log x_{2}=l_{q}(z_{2})}\right)
\]
and 
\[
\sum_{n\in \C}\left(\langle w'_{(4)}, {\cal Y}_1(w_{(1)}, x_1) 
\pi_{n}({\cal Y}_2(w_{(2)},
x_2)w_{(3)})\rangle_{W_4}\lbar_{x_1^{n}=e^{nl_{p}(z_1)},\;
\log x_{1}=l_{p}(z_{1}),\; \;x_2^{n}=e^{nl_{q}(z_2)},
\; \log x_{2}=l_{q}(z_{2})}\right),
\]
respectively. We shall further use the notations 
\begin{equation}\label{iterabbr}
\langle w'_{(4)}, 
{\cal Y}^1({\cal Y}^2(w_{(1)}, x_0)w_{(2)},
x_2)w_{(3)} \rangle_{W_4}\lbar_{x_0=z_0, \;x_2=z_2}
\end{equation}
and 
\begin{equation}\label{prodabbr}
\langle w'_{(4)}, {\cal Y}_1(w_{(1)}, x_1) 
{\cal Y}_2(w_{(2)},
x_2)w_{(3)}\rangle_{W_4}\lbar_{x_1=z_1, \; x_2=z_2},
\end{equation}
or even more simply, the notations
\begin{equation}\label{iterateabbreviation}
\langle w'_{(4)}, 
{\cal Y}^1({\cal Y}^2(w_{(1)}, z_0)w_{(2)},
z_2)w_{(3)} \rangle_{W_4}
\end{equation}
and 
\begin{equation}\label{productabbreviation}
\langle w'_{(4)}, {\cal Y}_1(w_{(1)}, z_1) 
{\cal Y}_2(w_{(2)},
z_2)w_{(3)}\rangle_{W_4}
\end{equation}
to denote (\ref{4itm}) and (\ref{4prm}), respectively, where 
we are taking $p=0$ in the  notation of 
(\ref{im:f(z)}) for both substitutions, except for occasions when we
explicitly specify different values of $p$, such as in the proof above.
We shall also use similar notations to denote 
series obtained {}from products and iterates of more than 
two intertwining operators. 

\begin{defi}\label{conv-conditions}
{\rm We call either of the two equivalent conditions in Proposition
\ref{convergence} the {\it convergence condition for intertwining maps
in the category ${\cal C}$}.}
\end{defi}

We need the following concept concerning unique expansion of an
analytic function in terms of powers of $z$ and $\log z$  
(recall our choice of the branch of $\log z$ in
(\ref{branch1}) and thus the branch of $z^{\alpha}$, $\alpha\in \C$):

\begin{defi}\label{uniqueexpset}
{\rm We call a subset ${\cal S}$ of ${\mathbb C}\times{\mathbb C}$ a
{\em unique expansion set} if the absolute convergence to $0$ on some
nonempty open subset of ${\mathbb C}^{\times}$ of any series
\[
\sum_{(\alpha,\beta)\in{\cal S}} a_{\alpha,\beta}z^\alpha(\log
z)^\beta, \;\;\; a_{\alpha,\beta} \in {\mathbb C},
\]
implies that $a_{\alpha,\beta}=0$ for all $(\alpha,\beta)\in{\cal
S}$. }
\end{defi}

Of course, a subset of a unique expansion set is again a unique
expansion set.

\begin{rema}
{\rm It is easy to show that ${\mathbb Z}\times\{ 0,\dots,N\}$ is a
unique expansion set for any $N\in{\mathbb N}$; this is also a
consequence of Proposition \ref{real-exp-set} below.  On the other
hand, it is known that ${\mathbb C}\times\{0\}$ is {\em not} a unique
expansion set\footnote{We thank A.~Eremenko for informing us of this
result.}.}
\end{rema}

For the reader's convenience, we give the following generalization of
a standard result about Laurent series:

\begin{lemma}\label{po-ser-an}
Let $D$ be a subset of $\R$ and let
\[
\sum_{\alpha\in D}a_{\alpha}z^{\alpha} \;\; (a_{\alpha} \in {\mathbb C})
\]
be absolutely convergent on a (nonempty) open
subset of $\C^{\times}$.  Then
\[
\sum_{\alpha\in D}a_{\alpha}\alpha z^{\alpha}
\]
is absolutely and uniformly convergent near any $z$ in the
open subset. In particular, the sum $\sum_{\alpha\in
D}a_{\alpha}z^{\alpha}$ as a function of $z$ is analytic in the sense that 
it is analytic at $z$ when $z$ is in the open subset of $\C^{\times}$ and 
$\arg z>0$, and that it can be analytically extended to an analytic 
function in a neighborhood of $z$
when $z$ is in the intersection of the open subset and the 
positive real line.
More generally, let $E$ be an index set and let the multisum
\[
\sum_{\alpha\in D}\sum_{\beta\in E}a_{\alpha, \beta}z^{\alpha} \;\;
(a_{\alpha, \beta}\in \C)
\]
converge absolutely on a (nonempty) open
subset of $\C^{\times}$.  Then the conclusions above hold for the
multisums $\sum_{\alpha\in D}\sum_{\beta\in E}a_{\alpha, \beta}\alpha
z^{\alpha}$ and $\sum_{\alpha\in D}\sum_{\beta\in E}a_{\alpha,
\beta}z^{\alpha}$.
\end{lemma}
\pf
We prove only the case for the series $\sum_{\alpha\in
D}a_{\alpha}z^{\alpha}$;
the general case is completely analogous.  We need only prove that 
$\sum_{\alpha\in D}a_{\alpha}\alpha z^{\alpha}$
is absolutely and uniformly convergent near any 
$z$ in the open subset. Note that since the original series is 
absolutely convergent on an open subset of $\C^{\times}$, 
$\sum_{\alpha\in D, \; \alpha\ge 0}a_{\alpha}z^{\alpha}$
and 
$\sum_{\alpha\in D, \; \alpha< 0}a_{\alpha}z^{\alpha}$
are also absolutely 
convergent on the set. For any fixed $z_{0}$ in the 
set, we can always find $z_{1}$ and $z_{2}$ in the set
such that $|z_{1}|<|z_{0}|<|z_{2}|$ and both
$\sum_{\alpha\in D, \; \alpha\ge 0}a_{\alpha}z_{2}^{\alpha}$
and 
$\sum_{\alpha\in D, \; \alpha< 0}a_{\alpha}z_{1}^{\alpha}$
are absolutely convergent. Let $r_{1}$ and $r_{2}$ be numbers 
such that $|z_{1}|<r_{1}<|z_{0}|<r_{2}<|z_{2}|$. 
Since 
\[
\lim_{\alpha\to \infty}~^{\alpha}\!\!\!\sqrt{\alpha}=1,
\]
we can find $M>0$ such that 
\[
^{\alpha}\!\!\!\sqrt{\alpha}<\min\left(\frac{|z_{2}|}{r_{2}}, 
\frac{r_{1}}{|z_{1}|}\right)
\]
when  $\alpha>M$. 
But when $r_{1}<|z|<r_{2}$, 
\[
^{\alpha}\!\!\!\sqrt{\alpha}<\min\left(\frac{|z_{2}|}{r_{2}}, 
\frac{r_{1}}{|z_{1}|}\right)< \min\left(\frac{|z_{2}|}{|z|}, 
\frac{|z|}{|z_{1}|}\right)
\]
for $\alpha>M$, so for $z$ in the open subset and satisfying
$r_{1}<|z|<r_{2}$, we have
\begin{eqnarray}\label{po-ser-an-1}
\sum_{\alpha\in D, \; \alpha>M}|a_{\alpha}\alpha z^{\alpha}|
&=&\sum_{\alpha\in D, \; \alpha>M}|a_{\alpha}|
|^{\alpha}\!\!\!\sqrt{\alpha}z|^{\alpha}\nn
&\leq&\sum_{\alpha\in D, \; \alpha>M}|a_{\alpha}z_{2}^{\alpha}|
\end{eqnarray}
and 
\begin{eqnarray}\label{po-ser-an-2}
\sum_{\alpha\in D, \; \alpha< -M}|a_{\alpha}\alpha z^{\alpha}|
&=&\sum_{\alpha\in D, \; \alpha<-M}|a_{\alpha}|
\left|\frac{^{-\alpha}\!\!\!\sqrt{-\alpha}}{z}\right|^{-\alpha}\nn
&\leq&\sum_{\alpha\in D, \; \alpha<-M}|a_{\alpha}z_{1}^{\alpha}|.
\end{eqnarray}
On the other hand, we have
\begin{eqnarray}\label{po-ser-an-3}
\sum_{\alpha\in D, \; 0\le \alpha\le M}|a_{\alpha}\alpha z^{\alpha}|
&\le& M\sum_{\alpha\in D, \; 0\le \alpha\le M}|a_{\alpha} z^{\alpha}|\nn
&\le & M\sum_{\alpha\in D, \; 0\le \alpha\le M}|a_{\alpha} z_{2}^{\alpha}|
\end{eqnarray}
and 
\begin{eqnarray}\label{po-ser-an-4}
\sum_{\alpha\in D, \; 0>\alpha\ge  -M}|a_{\alpha}\alpha z^{\alpha}|
&\le &M\sum_{\alpha\in D, \; 0>\alpha\ge  -M}|a_{\alpha} z^{\alpha}|\nn
&\le & M\sum_{\alpha\in D, \; 0>\alpha\ge  -M}|a_{\alpha} z_{1}^{\alpha}|.
\end{eqnarray}
{}From (\ref{po-ser-an-1})--(\ref{po-ser-an-4}), we see that 
$\sum_{\alpha\in D}a_{\alpha}\alpha z^{\alpha}$ is absolutely and uniformly
convergent in the neighborhood of $z_{0}$ consisting of $z$ in the 
open subset satisfying $r_{1}<|z|<r_{2}$. 
\epfv

\begin{propo}\label{real-exp-set}
For any $N \in \N$, $\mathbb{R}\times\{0,\dots,N\}$ is a unique
expansion set.  In particular, for any subset $D$ of $\mathbb{R}$,
$D\times\{0,\dots,N\}$ is a unique expansion set.
\end{propo}
\pf Let $a_{n, i}\in \C$ for $n\in \R$ and $i=0, \dots, N$, and
suppose that
\[
\sum_{n\in \R}\sum_{i=0}^{N}a_{n, i}z^{n}(\log z)^{i}=
\sum_{n\in \R}\sum_{i=0}^{N}a_{n, i}e^{n\log z}(\log z)^{i}
\]
is absolutely convergent to $0$ for $z$ in some nonempty open subset
of $\C^{\times}$.  We want to prove that each $a_{n, i}=0$.

Fix $n_{0}\in \R$.  We shall prove that $a_{n_{0}, N}=0$, and thus
the result will follow by induction on $N$; the case $N=0$ is a
special case of the proof below.

On the given open set, both
\[
\sum_{n\ge n_{0}}\sum_{i=0}^{N}a_{n, i}e^{(n-n_{0})\log z}(\log z)^{i}
\]
and 
\[
-\sum_{n< n_{0}}\sum_{i=0}^{N}a_{n, i}e^{(n-n_{0})\log z}(\log z)^{i}
\]
are absolutely convergent and we have
\begin{equation}\label{r-n-1}
\sum_{n\ge n_{0}}\sum_{i=0}^{N}a_{n, i}e^{(n-n_{0})\log z}(\log z)^{i}
=-\sum_{n< n_{0}}\sum_{i=0}^{N}a_{n, i}e^{(n-n_{0})\log z}(\log z)^{i}.
\end{equation}
Moreover, deleting $z=1$ {}from the open set if necessary, we observe
that for each $i=0, \dots, N$, the series
\[
\sum_{n\in \R}a_{n, i}e^{n\log z}
\]
is absolutely convergent on our open set.

Choose $z_{1}$ and $z_{2}$ in the open set satisfying
$|z_{1}|>|z_{2}|>0$; then for each $i=0, \dots, N$, $\sum_{n\in
\R}a_{n, i}e^{n\log z_{1}}$ and $\sum_{n\in \R}a_{n, i}e^{n\log
z_{2}}$ are absolutely convergent.

Since for $z'\in \C$ satisfying $|e^{z'}|\le |z_{1}|$ and $i=0, \dots, N$,
\begin{eqnarray*}
\sum_{n\ge n_{0}}|a_{n, i}||e^{(n-n_{0})z'}|
&= &\sum_{n\ge n_{0}}|a_{n, i}||e^{z'}|^{n-n_{0}}\nn
&\le &\sum_{n\ge n_{0}}|a_{n, i}||z_{1}|^{n-n_{0}},
\end{eqnarray*}
which is convergent, the series $\sum_{n\ge n_{0}}a_{n,
i}e^{(n-n_{0})z'}$ is absolutely convergent, and in particular, the
series
\begin{equation}\label{r-n->}
\sum_{n\ge n_{0}} a_{n, i}e^{(n-n_{0})\log z}
\end{equation}
is absolutely convergent for $z \in \C^{\times}$ satisfying $|z| \le
|z_{1}|$.  Thus by Lemma \ref{po-ser-an}, (\ref{r-n->}) defines an
analytic function on the region $0 < |z| < |z_1|$, with $0 \le \arg z
< 2\pi$.  Hence we have a single-valued analytic function
\[
f_{1}(z')=\sum_{n\ge n_{0}}\sum_{i=0}^{N}a_{n, i}e^{(n-n_{0})z'}(z')^{i}
\]
on the region $|e^{z'}|< |z_{1}|$, or equivalently, $\Re{(z')}< \log
|z_{1}|$.

Similarly, for $z'$ in the region $|e^{z'}|\ge |z_{2}|$ and $i=0,
\dots, N$,
\begin{eqnarray*}
\sum_{n< n_{0}}|a_{n, i}||e^{(n-n_{0})z'}|
&= &\sum_{n< n_{0}}|a_{n, i}||e^{z'}|^{n-n_{0}}\nn
&\le &\sum_{n< n_{0}}|a_{n, i}||z_{2}|^{n-n_{0}},
\end{eqnarray*}
so that the series $-\sum_{n< n_{0}}a_{n, i}e^{(n-n_{0})z'}$
is absolutely convergent.  Thus the series
\[
-\sum_{n < n_{0}}a_{n, i}e^{(n-n_{0})\log z}
\]
is absolutely convergent for $z \in \C^{\times}$ satisfying $|z|\ge
|z_{2}|$, defining, as above, a multivalued analytic function on the
region $|z| > |z_2|$ and hence a single-valued analytic function
\[
f_{2}(z')=-\sum_{n < n_{0}}\sum_{i=0}^{N}a_{n, i}e^{(n-n_{0})z'}(z')^{i}
\]
on the region $|e^{z'}|> |z_{2}|$, or equivalently, $\Re{(z')}> \log
|z_{2}|$.

We now define a single-valued analytic function $f(z')$ of $z'$ on the
whole plane $\C$ as follows: For $z'$ satisfying $|z_{2}|< |e^{z'}|<
|z_{1}|$, or equivalently, $\log |z_{2}| < \Re{(z')}< \log |z_{1}|$,
we have $f_{1}(z')=f_{2}(z')$, since $f_1$ and $f_2$, defined and
analytic on this region, agree on a nonempty open subset of this
region, in view of (\ref{r-n-1}). Thus we obtain a single-valued
analytic function $f(z')$ defined on the whole $z'$-plane by
\[
f(z')=\left\{\begin{array}{ll}
f_{1}(z'),& \Re{(z')} < \log |z_{1}|\\
f_{2}(z'),& \Re{(z')} > \log |z_{2}|.
\end{array}\right.
\]

When $\Re{(z')}< \log |z_{1}|$, we have
\begin{eqnarray*}
|f(z')|&\le &\sum_{n\ge n_{0}}\sum_{i=0}^{N}|a_{n, i}||e^{(n-n_{0})z'}|
|z'|^{i}\nn
&\le & \sum_{n\ge n_{0}}\sum_{i=0}^{N}|a_{n, i}||z_{1}|^{n-n_{0}}|z'|^{i}\nn
&=& \sum_{i=0}^{N}\left(\sum_{n\ge n_{0}}|a_{n, i}||z_{1}|^{n-n_{0}}\right)|z'|^{i}
\end{eqnarray*}
and when $\Re{(z')}> \log |z_{2}|$,
\begin{eqnarray*}
|f(z')|&\le &\sum_{n< n_{0}}\sum_{i=0}^{N}|a_{n, i}||e^{(n-n_{0})z'}|
|z'|^{i}\nn
&\le& \sum_{n< n_{0}}\sum_{i=0}^{N}|a_{n, i}||z_{2}|^{n-n_{0}}|z'|^{i}\nn
&=&\sum_{i=0}^{N}\left(\sum_{n< n_{0}}|a_{n, i}||z_{2}|^{n-n_{0}}\right)|z'|^{i}.
\end{eqnarray*}
Let 
\[
M_{i}=\max\left(\sum_{n\ge n_{0}}|a_{n, i}||z_{1}|^{n-n_{0}}, 
\sum_{n< n_{0}}|a_{n, i}||z_{2}|^{n-n_{0}}\right)
\]
for $i=0, \dots, N$. Then for $z'\in \C$ with $|z'| \ge 1$,
\[
|f(z')| \le \sum_{i=0}^{N}M_{i}|z'|^{i} \le
\left(\sum_{i=0}^{N}M_{i}\right)|z'|^{N},
\]
so that $f(z')$ is a polynomial of degree at most $N$ and in
particular, $\lim_{z'\to \infty}(z')^{-N}f(z')$ exists.

We now take the limit of $(z')^{-N}f(z')$ as $z'\to \infty$ along the
positive real line.  Let $M > \max(0,\log |z_{2}|)$.  When $z'\ge M$,
$f(z')=f_{2}(z')$, and for such $z'$,
\begin{eqnarray*}
\left|\sum_{n<n_{0}}\sum_{i=0}^{N}a_{n, i}e^{(n-n_{0})z'}(z')^{i-N}\right|
&\le& \sum_{n<n_{0}}\sum_{i=0}^{N}|a_{n, i}|e^{(n-n_{0})z'}(z')^{i-N}\nn
&\le&\sum_{n<n_{0}}\sum_{i=0}^{N}|a_{n, i}|e^{(n-n_{0})M} M^{i-N}.
\end{eqnarray*}
Since the right-hand side is convergent, the series
\[
-\sum_{n<n_{0}}\sum_{i=0}^{N}a_{n, i}e^{(n-n_{0})z'}(z')^{i-N}
\]
is uniformly convergent for  $z'\ge M$. Thus
\begin{eqnarray*}
\lefteqn{\lim_{z' \ge M, \; z'\to \infty}
-\sum_{n < n_{0}}\sum_{i=0}^{N}a_{n, i}e^{(n-n_{0})z'}(z')^{i-N}}\nn
&&=-\sum_{n<n_{0}}\sum_{i=0}^{N}\lim_{z' \ge M, \; z'\to \infty}
a_{n, i}e^{(n-n_{0})z'}(z')^{i-N}\nn
&&=0
\end{eqnarray*}
and so
\begin{equation}\label{r-n-5}
\lim_{z'\to \infty}(z')^{-N}f(z')
=\lim_{z' > 0, \; z'\to \infty}(z')^{-N}f(z')=0.
\end{equation}

Now let $M' < \min(0,\log |z_{1}|)$.  When $z'\le M'$,
$f(z')=f_{1}(z')$, and for such $z'$,
\begin{eqnarray*}
\left|\sum_{n\ge n_{0}}\sum_{i=0}^{N}a_{n, i}e^{(n-n_{0})z'}(z')^{i-N}\right|
&\le& \sum_{n\ge n_{0}}\sum_{i=0}^{N}|a_{n, i}|e^{(n-n_{0})z'}(-z')^{i-N}\nn
&\le&\sum_{n\ge n_{0}}\sum_{i=0}^{N}|a_{n, i}|e^{(n-n_{0})M'} (-M')^{i-N}.
\end{eqnarray*}
Since the right-hand side is convergent, the series 
\[
\sum_{n\ge n_{0}}\sum_{i=0}^{N}a_{n, i}e^{(n-n_{0})z'}(z')^{i-N}
\]
is uniformly convergent for  $z'\le M'$. Thus 
\begin{eqnarray*}
\lefteqn{\lim_{z' \le M', \; z'\to -\infty}
\sum_{n\ge n_{0}}\sum_{i=0}^{N}a_{n, i}e^{(n-n_{0})z'}(z')^{i-N}}\nn
&&=\sum_{n\ge n_{0}}\sum_{i=0}^{N}\lim_{z' \le M', \; z'\to -\infty}
a_{n, i}e^{(n-n_{0})z'}(z')^{i-N}\nn
&&=a_{n_{0}, N}
\end{eqnarray*}
and so we also have 
\begin{equation}\label{r-n-7}
\lim_{z'\to \infty}(z')^{-N}f(z')
=\lim_{z' < 0, \; z'\to -\infty}(z')^{-N}f(z')
=a_{n_{0}, N}.
\end{equation}
{}From (\ref{r-n-5}) and (\ref{r-n-7}), we obtain $a_{n_{0}, N}=0$.
\epfv

We will also need the following proposition and corollary, which
ensure that certain double sums converge when the corresponding
iterated sums and their derivatives converge:

\begin{propo}\label{log-coeff-conv<=>iterate-conv}
Let $D$ be a subset of $\R$ and $N$ a nonnegative integer. Then the series 
\begin{equation}\label{log-coeff-series}
\sum_{\alpha\in D}a_{\alpha, \beta}
z^\alpha 
\end{equation}
for $\beta=0, \dots, N$ 
are all absolutely convergent on some (nonempty) open subset of $\C^{\times}$
if and only if the series
\begin{equation}\label{iterate-series}
\sum_{\alpha\in D}\left(\sum_{\beta=0}^{N}a_{\alpha, \beta}
(\log z)^\beta\right) z^\alpha
\end{equation}
and the corresponding series of first and higher derivatives with respect to 
$z$,
viewed as series whose terms are the expressions
\[
\left(\sum_{\beta=0}^{N}a_{\alpha, \beta}
(\log z)^\beta\right) z^\alpha
\]
and their derivatives with respect to $z$, are absolutely convergent
on the same open subset.  The series of derivatives of
(\ref{iterate-series}) have the same format as (\ref{iterate-series}),
except that for the $n$-th derivative, the outer sum is over the set
$D-n$ and the inner sum has new coefficients in $\C$.
\end{propo}
\pf The last assertion is clear.

Assume the absolute convergence of (\ref{log-coeff-series})
for $\beta=0, \dots, N$.  Then the double series
\begin{equation}\label{double-series}
\sum_{\alpha\in D}\sum_{\beta=0}^{N}a_{\alpha, \beta}
z^\alpha (\log z)^\beta
\end{equation}
is absolutely convergent. We also know that the absolute convergence of
(\ref{double-series}) and
its (higher) derivatives implies the absolute convergence of
(\ref{iterate-series}) and its derivatives.  But using Lemma \ref{po-ser-an}
we see that the (higher) derivatives of (\ref{log-coeff-series}) are 
absolutely convergent. Since the (higher) derivatives of 
\begin{equation}\label{log-coeff-series-1}
\sum_{\alpha\in D}a_{\alpha, \beta}
z^\alpha (\log z)^{\beta}
\end{equation}
are (finite) linear combinations of the (higher) derivatives
of (\ref{log-coeff-series}) with coefficients containing integer
powers of $\log z$ and $z$, the (higher) derivatives of
(\ref{log-coeff-series-1}) are also absolutely convergent.
Thus the (higher) derivatives of (\ref{double-series})
are also absolutely convergent, and so (\ref{iterate-series}) and its 
derivatives are absolutely convergent.

Conversely, assume that (\ref{iterate-series}) and its derivatives are
absolutely convergent.  We need to show that (\ref{log-coeff-series})
is absolutely convergent at any $z_{0}$ in the open subset.  We
consider the series
\begin{equation}\label{iterate-series-1}
\sum_{\alpha\in D, \;\alpha\ge 0}\left(\sum_{\beta=0}^{N}a_{\alpha, \beta}
z_{2}^\beta\right) z_{1}^\alpha
\end{equation}
of functions
\[
\left(\sum_{\beta=0}^{N}a_{\alpha, \beta}
z_{2}^\beta\right) z_{1}^\alpha
\]
in two variables $z_{1}$ and $z_{2}$. Since $z_{0}$ is in the open
subset, we can find a smaller open subset inside the original one such
that for $z$ in this smaller one, $|z_{0}|<|z|$ and $|\log
z_{0}|<|\log z|$. We know that the series (\ref{iterate-series-1}) is
absolutely convergent when $z_{1}=z$, $z_{2}=\log z$ and $z$ is in the
original open subset.  For any $z_{1}$ and $z_{2}$ satisfying
$0<|z_{1}|<|z|$ and $z_{2}=\log z$ where $z$ is in the smaller open
subset,
\[
\sum_{\alpha\in D, \;\alpha\ge 0}\left|\sum_{\beta=0}^{N}a_{\alpha, \beta}
z_{2}^\beta\right| |z_{1}^\alpha|
\le \sum_{\alpha\in D, \;\alpha\ge 0}\left|\sum_{\beta=0}^{N}a_{\alpha, \beta}
(\log z)^\beta\right| |z^\alpha|
\]
is convergent. So in this case (\ref{iterate-series-1}) is absolutely
convergent. Since for any fixed $z_{2}=\log z$ 
where $z$ is in the smaller open subset, the numbers $z_{1}$ satisfying 
$0<|z_{1}|<|z|$ form an open subset, we can apply Lemma \ref{po-ser-an}
to obtain that
\begin{equation}\label{iterate-series-2}
\sum_{\alpha\in D, \;\alpha\ge 0}\frac{\partial}{\partial z_{1}}
\left(\left(\sum_{\beta=0}^{N}a_{\alpha, \beta}
z_{2}^\beta\right)  z_{1}^{\alpha}\right)
\end{equation}
is also absolutely convergent for any $z_{1}$ and $z_{2}$ satisfying
$0<|z_{1}|<|z|$ and $z_{2}=\log z$ where $z$ is in the smaller open
subset.

Also, by assumption, 
\begin{equation}\label{iterate-series-3}
\sum_{\alpha\in D, \; \alpha\ge 0}\left(\frac{\partial}{\partial z_{1}}
+\frac{1}{z_{1}}
\frac{\partial}{\partial z_{2}}\right)\left(
\left(\sum_{\beta=0}^{N}a_{\alpha, \beta}
z_{2}^{\beta}\right) z_{1}^\alpha\right)
\end{equation}
is absolutely convergent when $z_{1}=z$ and $z_{2}=\log z$ when $z$ is
in the original open subset.  For $z$ in the smaller open subset
and any $z_{1}$ and $z_{2}$ satisfying $0<|z_{1}|<|z|$ and $z_{2}=\log z$,
\begin{eqnarray*}
\lefteqn{\sum_{\alpha\in D, \; \alpha\ge 0}\left|
\left(\frac{\partial}{\partial z_{1}}
+\frac{1}{z_{1}}
\frac{\partial}{\partial z_{2}}\right)\left(
\left(\sum_{\beta=0}^{N}a_{\alpha, \beta}
z_{2}^{\beta}\right) z_{1}^\alpha\right)\right|}\nn
&&=\sum_{\alpha\in D, \; \alpha\ge 0}\left|
\left(\left(\sum_{\beta=0}^{N}a_{\alpha, \beta}
z_{2}^{\beta}\right) \alpha z_{1}^{\alpha-1}\right)
+\left(\left(\sum_{\beta=0}^{N}a_{\alpha, \beta}\beta
z_{2}^{\beta-1}\right) z_{1}^{\alpha-1}\right)\right|\nn
&&=\sum_{\alpha\in D, \; \alpha\ge 0}\left|
\left(\left(\sum_{\beta=0}^{N}a_{\alpha, \beta}
z_{2}^{\beta}\right)\alpha 
+\sum_{\beta=0}^{N}a_{\alpha, \beta}\beta
z_{2}^{\beta-1} \right)\right|\left|z_{1}^{\alpha-1}\right|\nn
&&\le |z {z_1}^{-1}|
\sum_{\alpha\in D, \; \alpha\ge 0}\left|
\left(\left(\sum_{\beta=0}^{N}a_{\alpha, \beta}
(\log z)^{\beta}\right)\alpha 
+\sum_{\beta=0}^{N}a_{\alpha, \beta}\beta
(\log z)^{\beta-1} \right)\right|\left|z^{\alpha-1}\right|\nn
&&=|z {z_1}^{-1}|
\sum_{\alpha\in D, \; \alpha\ge 0}\left|
\left(\left(\sum_{\beta=0}^{N}a_{\alpha, \beta}
(\log z)^{\beta}\right)\alpha z^{\alpha-1}
+\left(\sum_{\beta=0}^{N}a_{\alpha, \beta}\beta
(\log z)^{\beta-1}\right)z^{\alpha-1} \right)\right|\nn
&&=|z {z_1}^{-1}|
\sum_{\alpha\in D, \; \alpha\ge 0}\left|\frac{\partial}{\partial z}
\left(\left(\sum_{\beta=0}^{N}a_{\alpha, \beta}
(\log z)^{\beta}\right) 
z^{\alpha}\right)\right|
\end{eqnarray*}
(where we keep in mind that $\alpha - 1$ could be negative)
is convergent, so that (\ref{iterate-series-3}) is absolutely convergent
for such $z_{1}$ and $z_{2}$.
Thus, subtracting, we see that
\begin{equation}\label{iterate-series-4}
\sum_{\alpha\in D, \;\alpha\ge 0}
\frac{\partial}{\partial z_{2}}\left(
\left(\sum_{\beta=0}^{N}a_{\alpha, \beta}
z_{2}^{\beta}\right) z_{1}^\alpha\right)=
\sum_{\alpha\in D, \;\alpha\ge 0}
\left(\frac{\partial}{\partial z_{2}}
\left(\sum_{\beta=0}^{N}a_{\alpha, \beta}
z_{2}^{\beta}\right)\right) z_{1}^\alpha
\end{equation}
is also absolutely convergent 
for such $z_{1}$ and $z_{2}$. By Lemma \ref{po-ser-an},
\[
\sum_{\alpha\in D, \;\alpha\ge 0}
\frac{\partial}{\partial z_{1}}
\frac{\partial}{\partial z_{2}}\left(
\left(\sum_{\beta=0}^{N}a_{\alpha, \beta}
z_{2}^{\beta}\right) z_{1}^\alpha\right)
\]
is absolutely convergent 
for such $z_{1}$ and $z_{2}$.

Since $\frac{\partial}{\partial z_{1}}
+\frac{1}{z_{1}}
\frac{\partial}{\partial z_{2}}$ and $\frac{\partial}{\partial z_{2}}$
commute with each other,  we have 
\begin{eqnarray}\label{iterate-series-4.5}
\lefteqn{\sum_{\alpha\in D, \;\alpha\ge 0}
\left(\frac{\partial}{\partial z_{1}}
+\frac{1}{z_{1}}
\frac{\partial}{\partial z_{2}}\right)
\frac{\partial}{\partial z_{2}}\left(
\left(\sum_{\beta=0}^{N}a_{\alpha, \beta}
z_{2}^{\beta}\right) z_{1}^\alpha\right)}\nn
&&=\sum_{\alpha\in D, \;\alpha\ge 0}\frac{\partial}{\partial z_{2}}
\left(\frac{\partial}{\partial z_{1}}
+\frac{1}{z_{1}}
\frac{\partial}{\partial z_{2}}\right)
\left(\left(\sum_{\beta=0}^{N}a_{\alpha, \beta}
z_{2}^{\beta}\right) z_{1}^{\alpha}\right)\nn
&&=z_{1}\sum_{\alpha\in D, \;\alpha\ge 0}
\left(\frac{\partial}{\partial z_{1}}
+\frac{1}{z_{1}}
\frac{\partial}{\partial z_{2}}\right)^{2}
\left(\left(\sum_{\beta=0}^{N}a_{\alpha, \beta}
z_{2}^{\beta}\right) z_{1}^{\alpha}\right)\nn
&&\quad-z_{1}\sum_{\alpha\in D, \;\alpha\ge 0}\frac{\partial}{\partial z_{1}}
\left(\frac{\partial}{\partial z_{1}}
+\frac{1}{z_{1}}
\frac{\partial}{\partial z_{2}}\right)
\left(\left(\sum_{\beta=0}^{N}a_{\alpha, \beta}
z_{2}^{\beta}\right) z_{1}^{\alpha}\right).
\end{eqnarray}
By assumption, the first term on the right-hand side of 
(\ref{iterate-series-4.5})
is absolutely convergent when $z_{1}=z$ and $z_{2}=\log z$
and $z$ is in the original open subset, and then, by the same argument
as above,  is also absolutely 
convergent for $z_{1}$ and $z_{2}$ satisfying $0<|z_{1}|<|z|$,
$z_{2}=\log z$ and $z$ in the smaller open subset. 
By Lemma \ref{po-ser-an} and the absolute convergence of 
(\ref{iterate-series-3}) for such $z_{1}$ and $z_{2}$, the second 
term on the right-hand side of 
(\ref{iterate-series-4.5})
is also absolutely convergent for $z_{1}$ and $z_{2}$ 
satisfying $0<|z_{1}|<|z|$,
$z_{2}=\log z$ and $z$ in the smaller open subset. 
So the left-hand side of (\ref{iterate-series-4.5}) is absolutely 
convergent for such $z_{1}$ and $z_{2}$. Thus 
\begin{eqnarray*}
\lefteqn{\sum_{\alpha\in D, \;\alpha\ge 0}
\left(\left(\frac{\partial}{\partial z_{2}}\right)^{2}
\left(\sum_{\beta=0}^{N}a_{\alpha, \beta}
z_{2}^{\beta}\right)\right) z_{1}^\alpha}\nn
&&=\sum_{\alpha\in D, \;\alpha\ge 0}
\left(\frac{\partial}{\partial z_{2}}\right)^{2}\left(
\left(\sum_{\beta=0}^{N}a_{\alpha, \beta}
z_{2}^{\beta}\right) z_{1}^\alpha\right)\nn
&&=z_{1}\sum_{\alpha\in D, \;\alpha\ge 0}
\left(\frac{\partial}{\partial z_{1}}
+\frac{1}{z_{1}}
\frac{\partial}{\partial z_{2}}\right)
\frac{\partial}{\partial z_{2}}\left(
\left(\sum_{\beta=0}^{N}a_{\alpha, \beta}
z_{2}^{\beta}\right) z_{1}^\alpha\right)\nn
&&\quad -z_{1}\sum_{\alpha\in D, \;\alpha\ge 0}
\frac{\partial}{\partial z_{1}}
\frac{\partial}{\partial z_{2}}\left(
\left(\sum_{\beta=0}^{N}a_{\alpha, \beta}
z_{2}^{\beta}\right) z_{1}^\alpha\right)
\end{eqnarray*}
is absolutely convergent for $z$ in the smaller open subset
and any $z_{1}$ and $z_{2}$ satisfying $0<|z_{1}|<|z|$ and 
$z_{2}=\log z$.

Repeating these arguments, we obtain that 
\begin{equation}\label{iterate-series-5}
\sum_{\alpha\in D, \;\alpha\ge 0}
\left(\left(\frac{\partial}{\partial z_{2}}\right)^{k}
\left(\sum_{\beta=0}^{N}a_{\alpha, \beta}
z_{2}^{\beta}\right)\right) z_{1}^\alpha
\end{equation}
is absolutely convergent for such $z_{1}$ and $z_{2}$ and for $k\in
\N$.  Taking $k=N$, we see that
\[
\sum_{\alpha\in D, \; \alpha\ge 0}
a_{\alpha, N}z_{1}^\alpha
\]
is absolutely convergent for such $z_{1}$.
Continuing this process with $k=N-1, \dots, 0$ we obtain that 
\[
\sum_{\alpha\in D, \; \alpha\ge 0}
a_{\alpha, \beta}z_{1}^\alpha
\]
is absolutely convergent for such $z_{1}$ and each $\beta=0,
\dots, N$. Since $0<|z_{0}|<|z|$, we see that in the case $z_{1}=z_{0}$,
\[
\sum_{\alpha\in D, \;\alpha\ge 0}
a_{\alpha, \beta}z_{0}^\alpha
\]
is absolutely convergent for $\beta=0, \dots, N$.  

We also need to prove the absolute convergence of
\[
\sum_{\alpha\in D, \;\alpha< 0}
a_{\alpha, \beta}z_{0}^\alpha
\]
for $\beta=0, \dots, N$. The proof is completely analogous to the 
proof above except that 
we take a smaller open subset such that for $z$ in this smaller 
one, $|z_{0}|> |z|>0$ and $|\log z_{0}|>|\log z|$ instead of 
$|z_{0}|< |z|$ and $|\log z_{0}|<|\log z|$. 
Thus
\[
\sum_{\alpha\in D}
a_{\alpha, \beta}z_{0}^\alpha
\]
is absolutely convergent for $\beta=0, \dots, N$.  \epfv

\begin{corol}\label{double-conv<=>iterate-conv}
Let $D$ be a subset of $\R$ and $N$ a nonnegative integer.  Then the
double series (\ref{double-series})
is absolutely convergent on some (nonempty) open subset of $\C^{\times}$
if and only if the series
(\ref{iterate-series})
and the corresponding series of first and higher derivatives with respect to 
$z$,
viewed as series whose terms are the expressions
\[
\left(\sum_{\beta=0}^{N}a_{\alpha, \beta}(\log z)^\beta\right)
z^{\alpha} 
\]
and their derivatives with respect to $z$, are absolutely convergent
on the same open subset.  
\end{corol}
\pf By Proposition \ref{log-coeff-conv<=>iterate-conv}, we need only
prove that the absolute convergence of the double series
(\ref{double-series}) is equivalent to the absolute convergence of
each of the series (\ref{log-coeff-series}).

It is clear that the absolute convergence of each of the series
(\ref{log-coeff-series}) implies the absolute convergence of the
double series (\ref{double-series}).  Now assume the absolute
convergence of (\ref{double-series}). If $z\ne 1$, it is clear that
each of the series (\ref{log-coeff-series}) is absolutely convergent.
If $z=1$ is in the open subset, we can find $z_{1}$ and $z_{2}$ in the
open subset such that $|z_{1}|<1<|z_{2}|$. Then
\begin{eqnarray*}
\sum_{\alpha\in D}|a_{\alpha, \beta}|
&=&\sum_{\alpha\in D,\; \alpha\le 0}|a_{\alpha, \beta}|
+\sum_{\alpha\in D,\; \alpha > 0}|a_{\alpha, \beta}|\nn
&\le&\sum_{\alpha\in D,\; \alpha\le 0}|a_{\alpha, \beta}|
|z_{1}|^\alpha +\sum_{\alpha\in D,\; \alpha > 0}|a_{\alpha, \beta}|
|z_{2}|^\alpha.
\end{eqnarray*}
Since the right-hand side is convergent, the left-hand side is also
convergent. Thus each of the series (\ref{log-coeff-series}) is
absolutely convergent for all $z$ in the open subset. \epfv

\begin{assum}\label{assum-exp-set}
Throughout the remainder of this work, we shall assume that ${\cal C}$
satisfies the condition that for any object of ${\cal C}$, all the
(generalized) weights are real numbers and in addition there exists
$K\in \Z_{+}$ such that
\[
(L(0)-L(0)_{s})^{K}=0
\]
on the generalized module; when
$\mathcal{C}$ is in $\mathcal{M}_{sg}$ (recall Notation $\ref{MGM}$),
the latter assertion holds vacuously.
\end{assum}

In practice, ``virtually all the interesting examples'' satisfy this
assumption.

\begin{propo}\label{exp-set}
We have:
\begin{enumerate}

\item For any object $W$ of $\mathcal{C}$, the set $\{(n, i)\in
\C\times \N\;|\; (L(0)-n)^{i}W_{[n]}\ne 0\}$ is included in a (unique
expansion) set of the form $\R\times \{0, \dots, N\}$; when
$\mathcal{C}$ is in $\mathcal{M}_{sg}$, the set $\{(n, 0)\in \C\times
\N\;|\; W_{(n)}\ne 0\}$ is included in the (unique expansion) set $\R
\times \{0\}$.

\item For any objects $W_{1}$, $W_{2}$ and $W_{3}$ of $\mathcal{C}$,
any logarithmic intertwining operator $\Y$ of type ${W_{3}\choose
W_{1}W_{2}}$, and any $w_{(1)}\in W_{1}$, $w_{(2)}\in W_{2}$ and
$w'_{(3)}\in W'_{3}$, the powers of $x$ and $\log x$ occurring in
\begin{equation}\label{Yw1xw2}
\langle w'_{(3)}, \Y(w_{(1)}, x)w_{(2)}\rangle
\end{equation}
form a subset of a (unique expansion) set of the form $\R\times \{0,
\dots, N\}$, where $N$ depends only on $W_{1}$, $W_{2}$ and $W_{3}$
(and is independent of the three elements and independent of $\Y$);
when $\mathcal{C}$ is in $\mathcal{M}_{sg}$, the powers of $x$ in
(\ref{Yw1xw2}) form a (unique expansion) set of real numbers.

\end{enumerate}
\end{propo}
\pf This result follows immediately {}from Proposition
\ref{real-exp-set} and Proposition \ref{log:logwt} (see Remark
\ref{log:compM}, which specifies a value of $N$ for the second
assertion).  \epfv

\begin{rema}\label{assum-int}
{\rm The first assertion in Proposition \ref{exp-set} is a restatement
of Assumption \ref{assum-exp-set}, in view of Propostion
\ref{real-exp-set}.}
\end{rema}

Recall again the projections $\pi_{n}$ for $n\in \C$ {}from (\ref{pi_n}) and
Definition \ref{Wbardef} and also recall the notations 
(\ref{iter-abbr-pq}) and (\ref{prod-abbr-pq}).
We now prove the analyticity of products and iterates of intertwining
maps, in the following sense:

\begin{propo}\label{analytic}
Assume the convergence condition for intertwining maps in
$\mathcal{C}$ (recall Definition \ref{conv-conditions}), and let
$W_{1}$, $W_{2}$, $W_{3}$, $W_{4}$, $M_{1}$ and $M_{2}$ be objects of
${\cal C}$.
\begin{enumerate}

\item
Let $\Y_1\in \mathcal{V}_{W_{1}M_{1}}^{W_{4}}$ and $\Y_2\in
\mathcal{V}_{W_{2}W_{3}}^{M_{1}}$.  Then for any
$w_{(1)}\in W_1$, $w_{(2)}\in W_2$, $w_{(3)}\in
W_3$, $w'_{(4)}\in W'_4$ and $p, q\in \Z$, the sum of the absolutely 
convergent series 
\begin{eqnarray}\label{prod-p}
\lefteqn{\langle w'_{(4)}, {\cal Y}_{1}(w_{(1)}, x_1){\cal
Y}_{2}(w_{(2)},x_2)w_{(3)}\rangle\lbar_{x^{n}_{1}=e^{nl_{p}(z_{1})},
\;\log x_{1}=l_{p}(z_{1}),\;x^{n}_{2}=e^{nl_{q}(z_{2})},\;
\log x_{2}=l_{q}(z_{2})}}\nn
&&=\sum_{n\in \R}\left(\langle w'_{(4)}, {\cal Y}_{1}(w_{(1)}, x_1)
\pi_{n}({\cal
Y}_{2}(w_{(2)},x_2)w_{(3)})\rangle\lbar_{x^{n}_{1}=e^{nl_{p}(z_{1})},
\;\log x_{1}=l_{p}(z_{1}),\;x^{n}_{2}=e^{nl_{q}(z_{2})},\;
\log x_{2}=l_{q}(z_{2})}\right)\nn
\end{eqnarray}
is a single-valued analytic function on the region
given by $|z_{1}|>|z_{2}|>0$ and $0< \arg z_{1}, \arg z_{2}<2\pi$,
and for $k, l\in \N$, 
\begin{eqnarray}\label{prod-p-kl}
\lefteqn{\frac{\partial^{k+l}}{\partial z_{1}^{k}\partial z_{2}^{l}}
\langle w'_{(4)}, {\cal Y}_{1}(w_{(1)}, x_1){\cal
Y}_{2}(w_{(2)},x_2)w_{(3)}\rangle\lbar_{x^{n}_{1}=e^{nl_{p}(z_{1})},
\;\log x_{1}=l_{p}(z_{1}),\;x^{n}_{2}=e^{nl_{q}(z_{2})},\;
\log x_{2}=l_{q}(z_{2})}}\nn
&&\!\!\!\!\!\!\!\!\!=\langle w'_{(4)}, {\cal Y}_{1}(L(-1)^{k}w_{(1)}, x_1){\cal
Y}_{2}(L(-1)^{l}w_{(2)},x_2)w_{(3)}\rangle\lbar_{x^{n}_{1}=e^{nl_{p}(z_{1})},
\;\log x_{1}=l_{p}(z_{1}),\;x^{n}_{2}=e^{nl_{q}(z_{2})},\;
\log x_{2}=l_{q}(z_{2})}.\nn
\end{eqnarray}
Moreover, these analytic functions in (\ref{prod-p}) for different $p,
q\in \Z$ are different branches of a multivalued analytic function
defined on the region $|z_{1}|>|z_{2}|>0$ with the cut $\arg z_{1}=0$,
$\arg z_{2}=0$; similarly for each $k$ and $l$ for the analytic
functions in (\ref{prod-p-kl}).

\item
Analogously, let
$\Y^1\in \mathcal{V}_{M_{2}W_{3}}^{W_{4}}$ 
and $\Y^2\in
\mathcal{V}_{W_{1}W_{2}}^{M_{2}}$
Then for any
$w_{(1)}\in W_1$, $w_{(2)}\in W_2$, $w_{(3)}\in
W_3$, $w'_{(4)}\in W'_4$ and $p, q\in \Z$, the sum of the absolutely 
convergent series 
\begin{eqnarray}\label{iter-p}
\lefteqn{\langle w'_{(4)}, {\cal Y}^{1}({\cal
Y}^{2}(w_{(1)}, x_0)w_{(2)},x_2)w_{(3)}
\rangle\lbar_{x^{n}_{0}=e^{nl_{p}(z_{0})},
\;\log x_{0}=l_{p}(z_{0}),\;x^{n}_{2}=e^{nl_{q}(z_{2})},\;
\log x_{2}=l_{q}(z_{2})}}\nn
&&=\sum_{n\in \R}\left(\langle w'_{(4)}, {\cal Y}^{1}(\pi_{n}({\cal
Y}^{2}(w_{(1)}, x_0)w_{(2)}),x_2)w_{(3)}
\rangle\lbar_{x^{n}_{0}=e^{nl_{p}(z_{0})},
\;\log x_{0}=l_{p}(z_{0}),\;x^{n}_{2}=e^{nl_{q}(z_{2})},\;
\log x_{2}=l_{q}(z_{2})}\right)\nn
\end{eqnarray}
is a single-valued analytic function on the region
given by $|z_{2}|>|z_{0}|>0$ and $0< \arg z_{0}, \arg z_{2}<2\pi$,
and for $k, l\in \N$,
\begin{eqnarray}\label{iter-p-kl}
\lefteqn{\frac{\partial^{k+l}}{\partial z_{0}^{k}\partial z_{2}^{l}}
\langle w'_{(4)}, {\cal Y}^{1}({\cal
Y}^{2}(w_{(1)}, x_0)w_{(2)},x_2)w_{(3)}
\rangle\lbar_{x^{n}_{0}=e^{nl_{p}(z_{0})},
\;\log x_{0}=l_{p}(z_{0}),\;x^{n}_{2}=e^{nl_{q}(z_{2})},\;
\log x_{2}=l_{q}(z_{2})}}\nn
&&\!\!\!\!\!\!\!\!\!=\sum_{j=0}^{l}
{l\choose j}\langle w'_{(4)}, {\cal Y}^{1}({\cal
Y}^{2}(L(-1)^{k+j}w_{(1)}, x_0)\cdot\nn
&&\quad\quad\quad\quad\quad\quad\quad\quad
\cdot L(-1)^{l-j}w_{(2)},x_2)w_{(3)}\rangle
\lbar_{x^{n}_{0}=e^{nl_{p}(z_{0})},
\;\log x_{0}=l_{p}(z_{0}),\;x^{n}_{2}=e^{nl_{q}(z_{2})},\;
\log x_{2}=l_{q}(z_{2})}.\nn
\end{eqnarray}
Moreover, these analytic functions in (\ref{iter-p}) for different $p,
q\in \Z$, are different branches of a multivalued analytic function
defined on the region $|z_{2}|>|z_{0}|>0$ with the cut $\arg z_{0}=0$,
$\arg z_{2}=0$; similarly for each $k$ and $l$ for the analytic
functions in (\ref{iter-p-kl}).

\end{enumerate}
\end{propo}
\pf
We assume that $w_{(1)}$, $w_{(2)}$, $w_{(3)}$, 
$w'_{(4)}$ are homogeneous with respect to the 
generalized weight grading. The general case follows by
linearity. 

Since the series (\ref{prod-p}) is absolutely convergent in the region
$|z_{1}|>|z_{2}|>0$ and every term of this series is a single-valued
function on the region given by $|z_{1}|>|z_{2}|>0$ and $0< \arg
z_{1}, \arg z_{2}<2\pi$, its sum gives a single-valued function in the
same region.  By the $L(-1)$-derivative property for logarithmic
intertwining operators,
\begin{eqnarray}\label{prod-p-kl-1}
\lefteqn{\sum_{n\in \R}
\frac{\partial^{k+l}}{\partial z_{1}^{k}\partial z_{2}^{l}}
\langle w'_{(4)}, {\cal Y}_{1}(w_{(1)}, x_1)\pi_{n}({\cal
Y}_{2}(w_{(2)},x_2)w_{(3)})\rangle\lbar_{x^{n}_{1}=e^{nl_{p}(z_{1})},
\;\log x_{1}=l_{p}(z_{1}),\;x^{n}_{2}=e^{nl_{q}(z_{2})},\;
\log x_{2}=l_{q}(z_{2})}}\nn
&&=\sum_{n\in \R}\langle w'_{(4)}, {\cal Y}_{1}(L(-1)^{k}w_{(1)}, x_1)\cdot\nn
&&\quad\quad\quad\quad\quad\quad\quad
\cdot \pi_{n}({\cal
Y}_{2}(L(-1)^{l}w_{(2)},x_2)w_{(3)})\rangle\lbar_{x^{n}_{1}=e^{nl_{p}(z_{1})},
\;\log x_{1}=l_{p}(z_{1}),\;x^{n}_{2}=e^{nl_{q}(z_{2})},\;
\log x_{2}=l_{q}(z_{2})}\nn
\end{eqnarray}
for $k, l\in \N$, and by assumption, the series (\ref{prod-p-kl-1})
are absolutely convergent in the region $|z_{1}|>|z_{2}|>0$.

By Proposition \ref{exp-set} and the assumption that the vectors are
homogeneous, for fixed $z_2$ (so that $l_q(z_2)$ is also fixed), (\ref{prod-p}) is a
series of the form
\begin{equation}\label{prod-p-1}
\sum_{n\in \R}\left(\sum_{i=0}^{N}a_{n, i} e^{nl_{p}(z_{1})}
(l_{p}(z_{1}))^{i}\right)
\end{equation}
with $a_{n, i}\in \C$.  This series is absolutely
convergent in the region for $z_{1}$ given by $|z_{1}|>|z_{2}|>0$ and
$0< \arg z_{1}<2\pi$.  Replacing the logarithmic intertwining operator
${\cal Y}_{1}(\cdot,x)$ in (\ref{prod-p}) by the logarithmic
intertwining operator ${\cal Y}_{1}(\cdot,e^{-2\pi ip}x)$ (recall
Remark \ref{formalinvariance}; cf. Remark \ref{Ypp'}), we see that the
resulting analogue of (\ref{prod-p}) equals
\begin{equation}\label{prod-p-principalbranch}
\sum_{n\in \R}\left(\sum_{i=0}^{N}a_{n, i} z_{1}^n (\log z_1)^i\right),
\end{equation}
with the same coefficients $a_{n, i}$ as in (\ref{prod-p-1}) but using
the principal branch $\log z_1$ instead of $l_{p}(z_{1})$, and this
series is also absolutely convergent in the region for $z_1$ given by
$|z_{1}|>|z_{2}|>0$ and $0< \arg z_{1}<2\pi$.  Since for each $k \in
{\mathbb N}$ and $l=0$, (\ref{prod-p-kl-1}) is absolutely convergent
in the region $|z_{1}|>|z_{2}|>0$, as is the analogue of
(\ref{prod-p-kl-1}) with ${\cal Y}_{1}(\cdot,x)$ replaced by ${\cal
Y}_{1}(\cdot,e^{-2\pi ip}x)$ as above, the series
\[
\sum_{n\in \R}a_{n, i}z_{1}^{n}
\]
for $i=0, \dots, N$ are all absolutely convergent in the region for
$z_1$ given by $|z_{1}|>|z_{2}|>0$ and $0< \arg z_{1}<2\pi$, by
Proposition \ref{log-coeff-conv<=>iterate-conv}, and in particular,
the double series
\begin{equation}\label{prod-p-2}
\sum_{n\in \R}\sum_{i=0}^{N}a_{n, i} z_{1}^{n}(\log z_{1})^{i}
\end{equation}
is also absolutely convergent in the same region (as Corollary
\ref{double-conv<=>iterate-conv} states).  Then by Lemma
\ref{po-ser-an}, $\sum_{n\in \R}a_{n, i}z_{1}^{n}$ for $i=0, \dots, N$
as functions of $z_{1}$ are analytic in the same region.  Since
(\ref{prod-p-principalbranch}) as a function of $z_{1}$ is equal to
(\ref{prod-p-2}), it is also analytic in the same region.  Thus for
fixed $z_{2}$, the sum of the analogue of (\ref{prod-p}) as a function
of $z_{1}$ is analytic in the region given by $|z_{1}|>|z_{2}|>0$ and
$0< \arg z_{1}<2\pi$, and thus so is (\ref{prod-p}) itself; moreover,
(\ref{prod-p}) for different values of $p$ are different branches of
the same multivalued analytic function (\ref{prod-p-2}) (with $z_1$
replaced by $l_{p}(z_{1}))$ in the region for $z_1$ given by
$|z_{1}|>|z_{2}|>0$ with the cut $\arg z_{1}=0$, and the derivatives
of this function are given by the branches of the multivalued analytic
function (\ref{prod-p-kl}).  Across the cut $\arg z_{1}=0$, the
analytic function is given by (\ref{prod-p}) for adjacent values of
$p$.

The same argument shows that for fixed $z_{1}$, the sum of
(\ref{prod-p}) as a function of $z_{2}$ for different values of $q$
are different branches of the same multivalued analytic function in
the region for $z_2$ given by $|z_{1}|>|z_{2}|>0$ with the cut $\arg
z_{2}=0$, with derivatives given by (\ref{prod-p-kl}).  Thus the sum
of (\ref{prod-p}) as a function of $z_{1}$ and $z_{2}$ for different
values of $p$ and $q$ are the branches of a multivalued analytic
function in the region $|z_{1}|>|z_{2}|>0$, with derivatives given by
(\ref{prod-p-kl}).

An analogous argument proves the second half of the proposition, for
${\cal Y}^{1}$ and ${\cal Y}^{2}$.  \epfv

\begin{rema}
{\rm As usual, we shall use the same notation to denote an absolutely
convergent series and its sum. In particular, (\ref{prod-p}) and
(\ref{iter-p}) denote either the series or the sums of the series. The
proposition above says that these sums are in fact analytic functions
in $z_{1}$ and $z_{2}$ and can be analytically extended to multivalued
analytic functions on the regions $|z_{1}|>|z_{2}|>0$ and
$|z_{2}|>|z_{0}|>0$, respectively.}
\end{rema}

Again recall (\ref{pi_n}) and Definition \ref{Wbardef}, and recall the
notations (\ref{iterateabbreviation}) and (\ref{productabbreviation}).
Using Corollary \ref{double-conv<=>iterate-conv}, Propositions 
\ref{exp-set} and \ref{analytic}, we shall prove the
following result:

\begin{propo}\label{prod=0=>comp=0}
Assume the convergence condition for intertwining maps
in $\mathcal{C}$.  Let $z_1$,
$z_2$ be two nonzero complex numbers satisfying 
\[
|z_1|>|z_2|>0,
\]
and let $I_1\in \mathcal{M}[P(z_{1})]_{W_{1}M_{1}}^{W_{4}}$ and $I_2\in
\mathcal{M}[P(z_{2})]_{W_{2}W_{3}}^{M_{1}}$.  Let $w_{(1)}\in W_1$, $w_{(3)}\in
W_3$ and $w'_{(4)}\in W'_4$ be homogeneous elements with respect to 
the (generalized) weight gradings. Suppose that for all homogeneous
$w_{(2)}\in W_2$,
\[
\langle w'_{(4)}, I_1(w_{(1)}\otimes I_2(w_{(2)}\otimes w_{(3)}))\rangle=0.
\]
Then
\[
\langle w'_{(4)},  I_1({w_{(1)}}\otimes \pi_p I_2({w_{(2)}}\otimes 
w_{(3)}))\rangle=0
\]
for all $p\in \R$ and all $w_{(2)}\in W_2$. In particular, 
\[
\langle w'_{(4)},  \pi_{p} I_1({w_{(1)}}\otimes \pi_q I_2({w_{(2)}}\otimes 
w_{(3)}))\rangle=0
\]
for all $p, q\in \R$ and all $w_{(2)}\in W_2$. 
\end{propo}
\pf Recall the correspondence between $P(z)$-intertwining maps and
logarithmic intertwining operators of the same type 
(Proposition \ref{im:correspond}),
and the 
notation $\mathcal{Y}_{I, p}$, $p\in \Z$, for the logarithmic intertwining
operators corresponding to a $P(z)$-intertwining map $I$ ((\ref{YIp}) and 
(\ref{recover})). 

By Proposition \ref{analytic}, 
\begin{equation}\label{z2=>z}
\langle w'_{(4)}, \mathcal{Y}_{I_1,0}(w_{(1)},
z_1)\mathcal{Y}_{I_2,0}(w_{(2)}, z)w_{(3)}\rangle
\end{equation}
is a single-valued analytic function of $z$ on the region given by
$|z_{1}|>|z|>0$ and $0< \arg z<2\pi$, and its derivatives are given by
(\ref{prod-p-kl}) with $p=q=0$ and $k=0$, and with $z_{2}$ replaced by
$z$.  If $0< \arg z_{2}<2\pi$, then by the Taylor expansion of this
analytic function of $z$ at $z=z_{2}$ and (\ref{prod-p-kl}), we see
that for all $w_{(2)}\in W_2$ and for all $z$ in a sufficiently small
neighborhood of $z_{2}$,
\begin{eqnarray}\label{w2z}
\lefteqn{\langle w'_{(4)}, \mathcal{Y}_{I_1,0}(w_{(1)},
z_1)\mathcal{Y}_{I_2,0}(w_{(2)},z)w_{(3)}\rangle}\nn
&&=\sum_{i\in {\mathbb N}}\frac{(z-z_{2})^i}{i!}\langle w'_{(4)},
\mathcal{Y}_{I_1,0}(w_{(1)},z_1)
\mathcal{Y}_{I_2,0}(L(-1)^iw_{(2)},z_2)w_{(3)}\rangle\nn
&&=\sum_{i\in {\mathbb N}}\frac{(z-z_{2})^i}{i!}0=0.
\end{eqnarray}
If $\arg z_{2}=0$, that is, if $z_{2}$ is a positive real number, then
by Proposition \ref{analytic}, (\ref{z2=>z}) can be 
analytically extended to a single-valued analytic function of $z$
on a neighborhood of $z_{2}$ such that in the intersection of this 
neighborhood with the region $0< \arg z<\pi$, this function 
is equal to (\ref{z2=>z}). Then the same argument as above
shows that in this intersection, (\ref{w2z}) holds. 

In either case, we see that (\ref{w2z}) holds in a nonempty open
subset of the region $0< \arg z<2\pi$.  Now by Proposition 
\ref{exp-set}, Proposition \ref{log:logwt}(b) and the meaning of
the absolutely convergent series on the left-hand side of (\ref{w2z}),
for fixed $z_{1}\ne 0$, this series on the left-hand side is of the
form (\ref{iterate-series}), with $D\subset \R$, a unique expansion set.  By
Proposition \ref{analytic}, the higher-derivative series of the
left-hand side of (\ref{w2z}) are absolutely convergent.  Thus we can
apply Corollary \ref{double-conv<=>iterate-conv} to obtain that the
double series obtained {}from the left-hand side of (\ref{w2z}) by
taking the terms to be monomials in $z$ and $\log z$ is also
absolutely convergent to $0$ for $z$ in the open subset.  By the
definition of unique expansion set, we see that all of the
coefficients of the monomials in $z$ and $\log z$ of this double
series must be zero.  Hence we get
\[
\langle w'_{(4)}, \mathcal{Y}_{I_1,0}(w_{(1)}, z_1)
({w_{(2)}}^{\mathcal{Y}_{I_2,0}}_{n;\,k}
w_{(3)})\rangle=0
\]
for any homogeneous $w_{(2)}\in W_2$, $n\in \R$ and $k\in \N$. 
Since $w_{(1)}$, $w_{(3)}$ and $w'_{(4)}$ are homogeneous, we obtain
\[
\langle w'_{(4)}, I_1(w_{(1)}\otimes 
\pi_p I_2({w_{(2)}}\otimes  w_{(3)}))
\rangle=0
\]
for any homogeneous $w_{(2)}\in W_2$ and $p\in \R$, 
in view of Proposition
\ref{log:logwt}(b) and Proposition \ref{im:correspond}, and this
remains true for any $w_{(2)}\in W_2$.  The last statement is clear.
\epf

\begin{corol}\label{prospan}
Assume the convergence condition for intertwining maps
in $\mathcal{C}$ (recall Definition \ref{conv-conditions}).  Let $z_1$,
$z_2$ be two nonzero complex numbers satisfying 
\[
|z_1|>|z_2|>0.
\]
Suppose that the $P(z_2)$-tensor product of $W_2$ and $W_3$ and the
$P(z_1)$-tensor product of $W_1$ and $W_2 \boxtimes_{P(z_2)} W_3$ both
exist (recall Definition \ref{pz-tp}).
Then $W_1\boxtimes_{P(z_1)} (W_2\boxtimes_{P(z_2)}
W_3)$ is spanned (as a vector space) by all the elements of the form
\[
\pi_{n}(w_{(1)}\boxtimes_{P(z_1)}(w_{(2)}\boxtimes_{P(z_2)}w_{(3)}))
\] 
where $w_{(1)}\in W_1$, $w_{(2)}\in W_2$ and $w_{(3)}\in W_3$ are 
homogeneous with respect to the (generalized) weight gradings 
and $n\in \R$ (recall the notation
(\ref{boxtensorofelements})). 
\end{corol}
\pf Let $w'_{(4)}\in (W_1\boxtimes_{P(z_1)} (W_2\boxtimes_{P(z_2)}
W_3))'$ be homogeneous such that
\[
\langle w'_{(4)}, w_{(1)}\boxtimes_{P(z_1)} (w_{(2)}\boxtimes_{P(z_2)}
w_{(3)})\rangle=0
\]
for all homogeneous 
$w_{(1)}\in W_1$, $w_{(2)}\in W_2$ and $w_{(3)}\in W_3$. {}From
Proposition \ref{prod=0=>comp=0} we see that
\[
\langle w'_{(4)}, \pi_p({w_{(1)}}\boxtimes_{P(z_1)}
\pi_q({w_{(2)}}\boxtimes_{P(z_2)} w_{(3)}))\rangle=0
\]
for all $p, q\in \R$ and all 
$w_{(1)}\in W_1$, $w_{(2)}\in W_2$ and
$w_{(3)}\in W_3$. Since by Proposition \ref{span}, the set
\[
\{ \pi_p({w_{(1)}}\boxtimes_{P(z_1)} \pi_q({w_{(2)}}\boxtimes_{P(z_2)}
w_{(3)})) |\;p,q\in \R, w_{(1)}\in W_1, w_{(2)}\in W_2,
w_{(3)}\in W_3\}
\]
spans the space $W_1\boxtimes_{P(z_1)} (W_2\boxtimes_{P(z_2)}
W_3)$, we must have $w'_{(4)}=0$, and the result follows.
\epfv

Analogously, by similar proofs we have:

\begin{propo}\label{iter=0=>comp=0}
Assume the convergence condition for intertwining maps in ${\cal
C}$. Let $z_0$, $z_2$ be two nonzero complex numbers satisfying
\[
|z_2|>|z_0|>0,
\] 
and let $I^1\in \mathcal{M}[P(z_{2})]_{M_{2}W_{3}}^{W_{4}}$ and $I^2\in
\mathcal{M}[P(z_{0})]_{W_{1}W_{2}}^{M_{2}}$. Let $w'_{(4)}\in W'_4$, $w_{(2)}\in
W_2$ and $w_{(3)}\in W_3$ be homogeneous with respect to the (generalized)
weight gradings. Suppose that for all homogeneous $w_{(1)}\in W_1$,
\[
\langle w'_{(4)}, I^1(I^2(w_{(1)}\otimes w_{(2)}) \otimes w_{(3)})\rangle=0.
\]
Then
\[
\langle w'_{(4)}, I^1(\pi_p I^2(w_{(1)}\otimes w_{(2)})\otimes w_{(3)})
\rangle=0
\]
for all $p\in \R$ and all $w_{(1)}\in W_1$. In particular, 
\[
\langle w'_{(4)}, \pi_{p} I^1(\pi_q I^2(w_{(1)}\otimes w_{(2)})\otimes w_{(3)})
\rangle=0
\]
for all $p, q\in \R$ and all $w_{(1)}\in W_1$. \epf
\end{propo}

\begin{corol}\label{iterspan}
Assume the convergence condition for intertwining maps in ${\cal
C}$. Let $z_0$, $z_2$ be two nonzero complex numbers satisfying
\[
|z_2|>|z_0|>0.
\] 
Suppose that the $P(z_0)$-tensor product of $W_1$ and $W_2$ and the
$P(z_2)$-tensor product of $W_1 \boxtimes_{P(z_0)} W_2$ and $W_3$ both
exist. Then $(W_1\boxtimes_{P(z_0)} W_2)\boxtimes_{P(z_2)}
W_3$ is spanned by all the elements of the form
\[
\pi_{n}((w_{(1)}\boxtimes_{P(z_0)} w_{(2)})\boxtimes_{P(z_2)} w_{(3)})
\]
where $w_{(1)}\in W_1$, $w_{(2)}\in W_2$ and $w_{(3)}\in W_3$ are
homogeneous and 
$n\in \R$.\epf
\end{corol}

In addition to the definitions (\ref{prod-p}) and (\ref{iter-p}) of
the indicated products and iterates of intertwining maps, there is a
different, natural candidate for interpretations of the left-hand
sides of (\ref{prod-p}) and (\ref{iter-p}), involving multiple as
opposed to iterated sums, and we now show that these other
interpretations indeed agree with the definitions of these
expressions:

\begin{propo}\label{formal=proj}
Assume the convergence condition for intertwining maps in ${\cal
C}$. Let $z_1$, $z_2$ be two nonzero complex numbers satisfying
$|z_1|>|z_2|>0$ and let $\Y_1\in \mathcal{V}_{W_{1}M_{1}}^{W_{4}}$ 
and $\Y_2\in
\mathcal{V}_{W_{2}W_{3}}^{M_{1}}$.  Then for any $p, q\in \Z$, 
$w_{(1)}\in W_1$, $w_{(2)}\in W_2$, $w_{(3)}\in
W_3$ and $w'_{(4)}\in W'_4$, the series 
obtained by substituting $e^{nl_{p}(z_{1})}$, $e^{nl_{q}(z_{2})}$,
$l_{p}(z_{1})$ and $l_{q}(z_{2})$ for $x_{1}^{n}$, $x_{2}^{n}$, $\log x_{1}$
and $\log x_{2}$, respectively, in the formal series
\[
\langle w'_{(4)}, {\cal Y}_{1}(w_{(1)}, x_1){\cal
Y}_{2}(w_{(2)},x_2)w_{(3)}\rangle
\]
is absolutely convergent and its sum is equal to 
\begin{eqnarray*}
\lefteqn{\langle w'_{(4)}, {\cal Y}_{1}(w_{(1)}, x_1){\cal
Y}_{2}(w_{(2)},x_2)w_{(3)}\rangle\lbar_{x^{n}_{1}=e^{nl_{p}(z_{1})},\;
\log x_{1}=l_{p}(z_{1}), \;x^{n}_{2}=e^{nl_{q}(z_{2})},\;
\log x_{2}=l_{q}(z_{2})}}\nn
&&=\sum_{n\in \R}\left(\langle w'_{(4)}, {\cal Y}_{1}(w_{(1)}, x_1)\pi_{n}({\cal
Y}_{2}(w_{(2)},x_2)w_{(3)})\rangle\lbar_{x^{n}_{1}=e^{nl_{p}(z_{1})},\;
\log x_{1}=l_{p}(z_{1}), \;x^{n}_{2}=e^{nl_{q}(z_{2})},\;
\log x_{2}=l_{q}(z_{2})}\right).
\end{eqnarray*}
Analogously, let $z_0$, $z_2$ be two nonzero complex numbers satisfying
$|z_2|>|z_0|>0$ and let $\Y^1\in \mathcal{V}_{M_{2}W_{3}}^{W_{4}}$ 
and $\Y^2\in
\mathcal{V}_{W_{1}W_{2}}^{M_{2}}$. Then for any $p, q\in \Z$,
$w_{(1)}\in W_1$, $w_{(2)}\in W_2$, $w_{(3)}\in
W_3$ and $w'_{(4)}\in W'_4$, the series 
obtained by substituting $e^{nl_{p}(z_{0})}$, $e^{nl_{q}(z_{2})}$,
$l_{p}(z_{0})$ and $l_{q}(z_{2})$ for $x_{0}^{n}$, $x_{2}^{n}$, $\log x_{0}$
and $\log x_{2}$, respectively, in the formal series
\[
\langle w'_{(4)}, {\cal Y}^{1}({\cal
Y}^{2}(w_{(1)}, x_0)w_{(2)},x_2)w_{(3)}\rangle
\]
is absolutely convergent and its sum is equal to 
\begin{eqnarray*}
\lefteqn{\langle w'_{(4)}, {\cal Y}^{1}({\cal
Y}^{2}(w_{(1)}, x_0)w_{(2)},x_2)w_{(3)}
\rangle\lbar_{x^{n}_{0}=e^{nl_{p}(z_{0})},\;
\log x_{0}=l_{p}(z_{0}), \;x^{n}_{2}=e^{nl_{q}(z_{2})},\;
\log x_{2}=l_{q}(z_{2})}}\nn
&&=\sum_{n\in \R}\left(\langle w'_{(4)}, {\cal Y}^{1}(\pi_{n}({\cal
Y}^{2}(w_{(1)}, x_0)w_{(2)}),x_2)w_{(3)}
\rangle\lbar_{x^{n}_{0}=e^{nl_{p}(z_{0})},\;
\log x_{0}=l_{p}(z_{0}), \;x^{n}_{2}=e^{nl_{q}(z_{2})},\;
\log x_{2}=l_{q}(z_{2})}\right).
\end{eqnarray*}
\end{propo}
\pf
We prove only the first part, the second part being similar. 

If $w_{(1)}$, $w_{(2)}$, $w_{(3)}$ and $w'_{(4)}$ are homogeneous with
respect to the generalized weight gradings, then by Proposition
\ref{exp-set}, the first series is the triple series (recall
(\ref{log:map}) and Proposition \ref{log:logwt}(b))
\begin{equation}\label{triple-sum}
\sum_{n\in \R}\sum_{j=0}^{M}\sum_{i=0}^{N}\langle w'_{(4)}, 
(w_{(1)})^{\mathcal{Y}_{1}}_{\Delta-n-2, j}
(w_{(2)})^{{\cal
Y}_{2}}_{n, i}w_{(3)}\rangle e^{(-\Delta+n+1)l_{p}(z_{1})}
l_{p}(z_{1})^{j}e^{(-n-1)l_{q}(z_{2})}l_{q}(z_{2})^{i}
\end{equation}
and the second 
series is the corresponding iterated series 
\[
\sum_{n\in \R}\left(\sum_{j=0}^{M}\sum_{i=0}^{N}\langle w'_{(4)}, 
(w_{(1)})^{\mathcal{Y}_{1}}_{\Delta-n-2, j}
(w_{(2)})^{{\cal
Y}_{2}}_{n, i}w_{(3)}\rangle l_{p}(z_{1})^{j}
l_{q}(z_{2})^{i}\right)e^{(-\Delta+n+1)l_{p}(z_{1})}
e^{(-n-1)l_{q}(z_{2})},
\]
with
\[
\Delta=-\wt w_{(4)}'+\wt w_{(1)}+\wt w_{(2)}+\wt w_{(3)} \in \R.
\]
By assumption, the iterated series is absolutely convergent.
By Proposition \ref{analytic}, all the derivatives of 
the iterated series are also absolutely convergent.
Replacing $\Y_{2}(w_{(2)}, x)$ by $\Y_{2}(w_{(2)}, e^{-2\pi i q}x)$,
we see that the resulting iterated series 
\[
\sum_{n\in \R}\left(\sum_{j=0}^{M}\sum_{i=0}^{N}\langle w'_{(4)}, 
(w_{(1)})^{\mathcal{Y}_{1}}_{\Delta-n-2, j}
(w_{(2)})^{{\cal
Y}_{2}}_{n, i}w_{(3)}\rangle l_{p}(z_{1})^{j}
(\log z_{2})^{i}\right)e^{(-\Delta+n+1)l_{p}(z_{1})}
z_{2}^{-n-1}
\]
and its derivatives
are still absolutely convergent. 
Then by Proposition \ref{log-coeff-conv<=>iterate-conv},
with $z=z_{2}$,  the iterated 
series 
\[
\sum_{n\in \R}\left(\sum_{j=0}^{M}\langle w'_{(4)}, 
(w_{(1)})^{\mathcal{Y}_{1}}_{\Delta-n-2, j}
(w_{(2)})^{{\cal
Y}_{2}}_{n, i}w_{(3)}\rangle l_{p}(z_{1})^{j}\right)
e^{(-\Delta+n+1)l_{p}(z_{1})}
e^{(-n-1)l_{q}(z_{2})}l_{q}(z_{2})^{i}
\]
for $i=0, \dots, N$ are absolutely convergent. 
By the same argument
but with  $\Y_{1}(w_{(1)}, x)$ replaced by $\Y_{1}(w_{(1)}, e^{-2\pi i p}x)$
and with $z=z_{1}$ in Proposition \ref{log-coeff-conv<=>iterate-conv}, 
we see that the double series 
\[
\sum_{n\in \R}\sum_{j=0}^{M}\langle w'_{(4)}, 
(w_{(1)})^{\mathcal{Y}_{1}}_{\Delta-n-2, j}
(w_{(2)})^{{\cal
Y}_{2}}_{n, i}w_{(3)}\rangle l_{p}(z_{1})^{j}
e^{(-\Delta+n+1)l_{p}(z_{1})}
e^{(-n-1)l_{q}(z_{2})}l_{q}(z_{2})^{i}
\]
for $i=0, \dots, N$ are absolutely convergent. Thus
the triple series (\ref{triple-sum}), as a finite sum of 
these series, is also absolutely convergent.

In the general case, $w_{(1)}$, $w_{(2)}$, $w_{(3)}$ 
and $w'_{(4)}$ are finite sums of homogeneous vectors. 
Thus we have 
\begin{equation}\label{formal=proj-1}
\langle w'_{(4)}, {\cal Y}_{1}(w_{(1)}, x_1){\cal
Y}_{2}(w_{(2)},x_2)w_{(3)}\rangle
=\sum_{i=1}^{k}\langle w^{i \prime}_{(4)}, {\cal Y}_{1}(w^{i}_{(1)}, x_1)
{\cal
Y}_{2}(w^{i}_{(2)},x_2)w^{i}_{(3)}\rangle
\end{equation}
and 
\begin{eqnarray}\label{formal=proj-2}
\lefteqn{\langle w'_{(4)}, {\cal Y}_{1}(w_{(1)}, x_1){\cal
Y}_{2}(w_{(2)},x_2)w_{(3)}\rangle\lbar_{x^{n}_{1}=e^{nl_{p}(z_{1})},\;
\log x_{1}=l_{p}(z_{1}), \;x^{n}_{2}=e^{nl_{q}(z_{2})},\;
\log x_{2}=l_{q}(z_{2})}}\nn
&&=\sum_{i=1}^{k}\langle w^{i \prime}_{(4)}, 
{\cal Y}_{1}(w^{i}_{(1)}, x_1)
{\cal
Y}_{2}(w^{i}_{(2)},x_2)w^{i}_{(3)}\rangle
\lbar_{x^{n}_{1}=e^{nl_{p}(z_{1})},\;
\log x_{1}=l_{p}(z_{1}), \;x^{n}_{2}=e^{nl_{q}(z_{2})},\;
\log x_{2}=l_{q}(z_{2})},\nn
\end{eqnarray}
where $w^{i}_{(1)}\in W_{1}$, $w^{i}_{(2)}\in W_{2}$, 
$w^{i}_{(3)}\in W_{3}$ and $w^{i \prime}_{(4)}\in W'_{4}$ 
are homogeneous
with respect to the generalized weight gradings.  Substituting 
$e^{nl_{p}(z_{1})}$, $e^{nl_{q}(z_{2})}$,
$l_{p}(z_{1})$ and $l_{q}(z_{2})$
for $x_{1}^{n}$, $x_{2}^{n}$, $\log x_{1}$
and $\log x_{2}$, respectively, in
each term in the right-hand side of (\ref{formal=proj-1})
gives an absolutely convergent series whose sum is equal
to the corresponding term in the 
right-hand side of (\ref{formal=proj-2}), and so making these
substitutions in
the left-hand side of (\ref{formal=proj-1}) gives 
an absolutely convergent series whose sum is equal
to the left-hand side of (\ref{formal=proj-2}).  
\epfv

\begin{rema}\label{4notations}
{\rm Proposition \ref{formal=proj} in fact justifies the notations
that we have introduced in (\ref{iter-abbr-pq}) and
(\ref{prod-abbr-pq}) (and in particular, in
(\ref{iterabbr})--(\ref{productabbreviation})).  That is, with $z_1$
and $z_2$ satisfying the appropriate inequality, for any $p, q\in \Z$,
$w_{(1)}\in W_1$, $w_{(2)}\in W_2$, $w_{(3)}\in W_3$ and $w'_{(4)}\in
W'_4$, (\ref{prod-abbr-pq}) also means the absolutely convergent sum
of the multiple series obtained by substituting $e^{nl_{p}(z_{1})}$,
$e^{nl_{q}(z_{2})}$, $l_{p}(z_{1})$ and $l_{q}(z_{2})$ for
$x_{1}^{n}$, $x_{2}^{n}$, $\log x_{1}$ and $\log x_{2}$, respectively,
in the formal series
\[
\langle w'_{(4)}, {\cal Y}_{1}(w_{(1)}, x_1){\cal
Y}_{2}(w_{(2)},x_2)w_{(3)}\rangle;
\]
similarly for $z_0$ and $z_2$, (\ref{iter-abbr-pq}) and the formal
series
\[
\langle w'_{(4)}, {\cal Y}^{1}({\cal
Y}^{2}(w_{(1)}, x_0)w_{(2)},x_2)w_{(3)}\rangle.
\]
When $p=q=0$, (\ref{prodabbr}) (or (\ref{productabbreviation})) also
means the absolutely convergent sum of the multiple series obtained by
substituting $e^{n\log z_{1}}$, $e^{n\log z_{2}}$, $\log z_{1}$ and
$\log z_{2}$ for $x_{1}^{n}$, $x_{2}^{n}$, $\log x_{1}$ and $\log
x_{2}$, respectively, in the first formal series above; similarly for
(\ref{iterabbr}) (or (\ref{iterateabbreviation})) and the second
formal series above. Since for an absolutely convergent series, we use the same 
notation to denote the series and its sum, 
(\ref{iter-abbr-pq}) and
(\ref{prod-abbr-pq}) (and in particular,
(\ref{iterabbr})--(\ref{productabbreviation})) also denote 
the analytic functions given by the sums of the corresponding 
series. Moreover, if $w_{(1)}$, $w_{(2)}$, $w_{(3)}$ and $w'_{(4)}$
are finite sums of elements of the same generalized modules, then 
the same notations also mean the finite sum of the series or sums 
obtained from the summands of $w_{(1)}$, $w_{(2)}$, $w_{(3)}$ and $w'_{(4)}$.
In the rest of this work, we shall use these notations to mean any one of
these things, depending on what we need.}
\end{rema}

{}From Proposition \ref{formal=proj}, we immediately obtain:

\begin{corol}\label{formal=proj-cor}
Assume the convergence condition for intertwining maps in ${\cal
C}$. Let $\Y_1\in \mathcal{V}_{W_{1}M_{1}}^{W_{4}}$, $\Y_2\in
\mathcal{V}_{W_{2}W_{3}}^{M_{1}}$, 
$\Y_3\in \mathcal{V}_{W_{1}\widetilde{M}_{1}}^{W_{4}}$, $\Y_4\in
\mathcal{V}_{W_{2}W_{3}}^{\widetilde{M}_{1}}$ and let
$w_{(1)}, \widetilde{w}_{(1)}\in W_1$, $w_{(2)}, 
\widetilde{w}_{(2)}\in W_2$, $w_{(3)}, \widetilde{w}_{(3)}\in
W_3$ and $w'_{(4)}, \widetilde{w}'_{(4)}\in W'_4$. If 
\[
\langle w'_{(4)}, {\cal Y}_{1}(w_{(1)}, x_1){\cal
Y}_{2}(w_{(2)},x_2)w_{(3)}\rangle
=\langle \widetilde{w}'_{(4)}, {\cal Y}_{3}(\widetilde{w}_{(1)}, x_1){\cal
Y}_{4}(\widetilde{w}_{(2)},x_2)\widetilde{w}_{(3)}\rangle,
\] 
then for any $p, q\in \Z$ and $z_{1}, z_{2}\in \C$ 
satisfying $|z_{1}|>|z_{2}|>0$,
\begin{eqnarray}
\lefteqn{\langle w'_{(4)}, {\cal Y}_{1}(w_{(1)}, x_1){\cal
Y}_{2}(w_{(2)},x_2)w_{(3)}\rangle\lbar_{x^{n}_{1}=e^{nl_{p}(z_{1})},\;
\log x_{1}=l_{p}(z_{1}), \;x^{n}_{2}=e^{nl_{q}(z_{2})},\;
\log x_{2}=l_{q}(z_{2})}}\nn
&&=\langle \widetilde{w}'_{(4)}, {\cal Y}_{3}(\widetilde{w}_{(1)}, x_1){\cal
Y}_{4}(\widetilde{w}_{(2)},x_2)\widetilde{w}_{(3)}
\rangle\lbar_{x^{n}_{1}=e^{nl_{p}(z_{1})},\;
\log x_{1}=l_{p}(z_{1}), \;x^{n}_{2}=e^{nl_{q}(z_{2})},\;
\log x_{2}=l_{q}(z_{2})}.\nn
\end{eqnarray}
Analogously, let $\Y^1\in \mathcal{V}_{M_{2}W_{3}}^{W_{4}}$,
$\Y^2\in
\mathcal{V}_{W_{1}W_{2}}^{M_{2}}$, $\Y^3\in 
\mathcal{V}_{\widetilde{M}_{2}W_{3}}^{W_{4}}$ and
$\Y^4\in
\mathcal{V}_{W_{1}W_{2}}^{\widetilde{M}_{2}}$ and let 
$w_{(1)}, \widetilde{w}_{(1)}\in W_1$, $w_{(2)}, 
\widetilde{w}_{(2)}\in W_2$, $w_{(3)}, \widetilde{w}_{(3)}\in
W_3$ and $w'_{(4)}, \widetilde{w}'_{(4)}\in W'_4$. If 
\[
\langle w'_{(4)}, {\cal Y}^{1}({\cal
Y}^{2}(w_{(1)}, x_0)w_{(2)},x_2)w_{(3)}\rangle
=\langle \widetilde{w}'_{(4)}, {\cal Y}^{3}({\cal
Y}^{4}(\widetilde{w}_{(1)}, x_0)\widetilde{w}_{(2)},x_2)
\widetilde{w}_{(3)}\rangle,
\]
then for $p, q\in \Z$ and
$z_{0}, z_{2}\in \C$
satisfying $|z_{0}|>|z_{2}|>0$, 
\begin{eqnarray}
\lefteqn{\langle w'_{(4)}, {\cal Y}^{1}({\cal
Y}^{2}(w_{(1)}, x_0)w_{(2)},x_2)w_{(3)}
\rangle\lbar_{x^{n}_{0}=e^{nl_{p}(z_{0})},\;
\log x_{0}=l_{p}(z_{0}), \;x^{n}_{2}=e^{nl_{q}(z_{2})},\;
\log x_{2}=l_{q}(z_{2})}}\nn
&&=\langle \widetilde{w}'_{(4)}, {\cal Y}^{3}({\cal
Y}^{4}(\widetilde{w}_{(1)}, x_0)\widetilde{w}_{(2)},x_2)\widetilde{w}_{(3)}
\rangle\lbar_{x^{n}_{0}=e^{nl_{p}(z_{0})},\;
\log x_{0}=l_{p}(z_{0}), \;x^{n}_{2}=e^{nl_{q}(z_{2})},\;
\log x_{2}=l_{q}(z_{2})}.\hspace{2em}\square\nn
\end{eqnarray}
\end{corol}

\begin{rema}
{\rm One can generalize the convergence condition for two intertwining
maps and the results above to products and iterates of any number of
intertwining maps. The convergence conditions for three
intertwining maps and the spanning properties in the case of four
generalized modules will be needed in Section 12 in the proof of the
commutativity of the pentagon diagram, and we will discuss these
conditions and properties in Section 12. }
\end{rema}

\begin{rema}\label{weakly-abs-conv}
{\rm The convergence studied in this section can easily be formulated
as special cases of the following general notion: Let $W$ be a
(complex) vector space and let $\langle \cdot,\cdot \rangle:
W^{*}\times W\to \C$ be the pairing between the dual space $W^{*}$ and
$W$.  Consider the weak topology on $W^{*}$ defined by this pairing,
so that $W^{*}$ becomes a Hausdorff locally convex topological vector
space.  Let $\sum_{n\in I}w^{*}_{n}$ be a formal series in $W^{*}$, 
where $I$ is an index set. We say that $\sum_{n\in I}w^{*}_{n}$ is {\it weakly
absolutely convergent} if for all $w\in W$, the formal series
\begin{equation}\label{sum-w}
\sum_{n\in I}\langle w^{*}_{n}, w\rangle
\end{equation}
of complex numbers is absolutely
convergent. Note that if $\sum_{n\in I}w^{*}_{n}$ is weakly
absolutely convergent, then (\ref{sum-w}) as $w$ ranges through $W$
defines a (unique) element of $W^{*}$, and the formal series is in
fact convergent to this element in the weak topology. This element is
the {\it sum} of the series and is denoted using the same notation
$\sum_{n\in I}w^{*}_{n}$. In this section, the convergence that we
have been discussing amounts to the weak absolute convergence of
formal series in $(W')^{*}$ for an object $W$ of $\mathcal{C}$, and
this kind of convergence will again be used in Section 8.  In Section
9, we will use this notion for $(W_{1}\otimes W_{2})^{*}$ where
$W_{1}$ and $W_{2}$ are generalized $V$-modules, and in Section 12, we
will be using more general cases.}
\end{rema}

\newpage

\setcounter{equation}{0}
\setcounter{rema}{0}

\section{$P(z_1,z_2)$-intertwining maps and the corresponding
compatibility condition}

In this section we first prove some natural identities satisfied by
products and iterates of logarithmic intertwining operators and of
intertwining maps.  These identities were first proved in
\cite{tensor4} for intertwining operators and intertwining maps among
ordinary modules.  We also prove a list of identities relating
products of formal delta functions, as was done in \cite{tensor4}.
Using all these identities as motivation, we define
``$P(z_1,z_2)$-intertwining maps'' and study their basic properties,
by analogy with the relevant parts of the study of $P(z)$-intertwining
maps in Sections 4 and 5.  The notion of $P(z_1,z_2)$-intertwining map
is new; the treatment in this section is different {}from that in
\cite{tensor4}, even for the case of ordinary intertwining operators.

At the end of this section, we show that products and iterates of
intertwining maps or of logarithmic intertwining operators ``factor
through'' suitable tensor product modules in a unique way.

It is possible to define ``tensor products of three modules,'' as
opposed to iterated tensor products, and $P(z_1,z_2)$-intertwining
maps would play the same role for such tensor products of three
modules that $P(z)$-intertwining maps play for tensor products of two
modules.  However, one would of course in addition need appropriate
natural isomorphisms between triple tensor products and the
corresponding iterated tensor products, and much more than
$P(z_1,z_2)$-intertwining maps (as defined here) would be necessary
for this; see Section 9 below, in particular.  Since we do not need
``tensor products of three modules'' in this work, we will not
formally introduce and study them.

We recall our continuing Assumptions \ref{assum}, \ref{assum-c} and
\ref{assum-exp-set} concerning our category ${\cal C}$.

Recall the Jacobi identity (\ref{log:jacobi}) in the definition of the
notion of logarithmic intertwining operator associated with
generalized modules $(W_1,Y_1)$, $(W_2,Y_2)$ and $(W_3,Y_3)$ for a
M\"obius (or conformal) vertex algebra $V$.  Suppose that we also have
generalized modules $(W_4,Y_4)$, $(M_1,Y_{M_1})$ and $(M_2,Y_{M_2})$.
Then {}from (\ref{log:jacobi}) we see that a product of logarithmic
intertwining operators of types ${W_4\choose W_1 M_1}$ and
${M_1\choose W_2 W_3}$ satisfies an identity analogous to
(\ref{log:jacobi}), as does an iterate of logarithmic intertwining
operators of types ${W_4\choose M_2 W_3}$ and ${M_2\choose W_1 W_2}$:

Let ${\cal Y}_1$ and ${\cal Y}_2$ be logarithmic intertwining
operators of types ${W_4\choose W_1 M_1}$ and ${M_1\choose W_2 W_3}$,
respectively.  Then for $v\in V$, $w_{(1)}\in W_1$, $w_{(2)}\in W_2$
and $w_{(3)}\in W_3$, the product of ${\cal Y}_1$ and ${\cal Y}_2$
satisfies the identity
\begin{eqnarray}
\lefteqn{\dlt{x_1}{x_0}{-y_1}\dlt{x_2}{x_0}{-y_2}Y_4(v,x_0){\cal
Y}_1(w_{(1)},y_1){\cal Y}_2(w_{(2)},y_2)w_{(3)}}\nno \\
&&= \dlt{y_1}{x_0}{-x_1}\dlt{x_2}{x_0}{-y_2}{\cal
Y}_1(Y_1(v,x_1)w_{(1)},y_1){\cal Y}_2(w_{(2)},y_2)w_{(3)}\nno \\
&& \quad+\dlt{x_1}{-y_1}{+x_0}\dlt{x_2}{x_0}{-y_2}{\cal
Y}_1(w_{(1)},y_1)Y_{M_1}(v,x_0){\cal Y}_2(w_{(2)},y_2)w_{(3)}\nno \\
&&= \dlt{y_1}{x_0}{-x_1}\dlt{x_2}{x_0}{-y_2}{\cal
Y}_1(Y_1(v,x_1)w_{(1)},y_1){\cal Y}_2(w_{(2)},y_2)w_{(3)}\nno \\
&& \quad+\dlt{x_1}{-y_1}{+x_0}\dlt{y_2}{x_0}{-x_2}{\cal
Y}_1(w_{(1)},y_1){\cal Y}_2(Y_2(v,x_2)w_{(2)},y_2)w_{(3)}\nno \\
&& \quad+\dlt{x_1}{-y_1}{+x_0}\dlt{x_2}{-y_2}{+x_0}{\cal
Y}_1(w_{(1)},y_1){\cal Y}_2(w_{(2)},y_2)Y_3(v,x_0)w_{(3)}.\nno \\
&& \label{Y12}
\end{eqnarray}
In addition, let ${\cal Y}^1$ and ${\cal Y}^2$ be logarithmic
intertwining operators of types ${W_4\choose M_2 W_3}$ and
${M_2\choose W_1 W_2}$, respectively.  Then for $v\in V$, $w_{(1)}\in
W_1$, $w_{(2)}\in W_2$ and $w_{(3)}\in W_3$, the iterate of ${\cal
Y}^1$ and ${\cal Y}^2$ satisfies the identity
\begin{eqnarray}
\lefteqn{\dlt{x_2}{x_0}{-y_2}\dlt{x_1}{x_2}{-y_0}Y_4(v,x_0){\cal
Y}^1({\cal Y}^2(w_{(1)},y_0)w_{(2)},y_2)w_{(3)}}\nno \\
&&= \dlt{y_2}{x_0}{-x_2}\dlt{x_1}{x_2}{-y_0}{\cal
Y}^1(Y_{M_2}(v,x_2){\cal Y}^2(w_{(1)},y_0)w_{(2)},y_2)w_{(3)}\nno
\\
&& \quad+\dlt{x_2}{-y_2}{+x_0}\dlt{x_1}{x_2}{-y_0}{\cal Y}^1({\cal
Y}^2(w_{(1)},y_0)w_{(2)},y_2)Y_3(v,x_0)w_{(3)}\nno \\
&&= \dlt{y_2}{x_0}{-x_2}\dlt{y_0}{x_2}{-x_1}{\cal Y}^1({\cal
Y}^2(Y_1(v,x_1)w_{(1)}),y_0)w_{(2)},y_2)w_{(3)}\nno \\
&& \quad+\dlt{y_2}{x_0}{-x_2}\dlt{x_1}{-y_0}{+x_2}{\cal Y}^1({\cal
Y}^2(w_{(1)},y_1)Y_2(v,x_2)w_{(2)},y_2)w_{(3)}\nno \\
&& \quad+\dlt{x_2}{-y_2}{+x_0}\dlt{x_1}{x_2}{-y_0}{\cal Y}^1({\cal
Y}^2(w_{(1)},y_0)w_{(2)},y_2)Y_3(v,x_0)w_{(3)}.\nno \\ &&
\label{Y34}
\end{eqnarray}

Under natural hypotheses motivated by Section 7, we will need to
specialize the formal variables $y_1$, $y_2$ and $y_0$ to complex
numbers $z_1$, $z_2$ and $z_0$, respectively, in (\ref{Y12}) and
(\ref{Y34}), when $|z_1|>|z_2|>0$ and $|z_2|>|z_0|>0$.  For this,
following \cite{tensor4}, we will need the next lemma, on products of
formal delta functions, with certain of the variables being complex
variables in suitable domains.  Our formulation and proof here are
different {}from those in \cite{tensor4}.  In addition to justifying the
specializations just indicated, this lemma will give us the natural
relation between the specialized expressions (\ref{F12}) and
(\ref{F34}) below.

\begin{lemma}\label{deltalemma}
Let $z_1$ and $z_2$ be complex numbers and set $z_0=z_1-z_2$.  Then
the left-hand sides of the following expressions converge absolutely
in the indicated domains, in the sense that the coefficient of each
monomial in the formal variables $x_0$, $x_1$ and $x_2$ is an
absolutely convergent series in the two variables related by the
inequalities, and the following identities hold:
\begin{eqnarray}
\dlt{x_1}{x_0}{-z_1}\dlt{x_2}{x_0}{-z_2}&=&
\dlt{x_2}{x_0}{-z_2}\dlt{x_1}{x_2}{-z_0}\nno\\
&&\mbox{\rm for arbitrary } z_1, z_2;\label{l1}\\
\dlt{z_1}{x_0}{-x_1}\dlt{x_2}{x_0}{-z_2} & = &
\dlt{x_0}{z_1}{+x_1}\dlt{x_2}{z_0}{+x_1}\nno \\
&& \mbox{\rm if }|z_1|>|z_2|;\label{l2-1}\\
\dlt{z_2}{x_0}{-x_2}\dlt{z_0}{x_2}{-x_1}&=&
\dlt{x_0}{z_1}{+x_1}\dlt{x_2}{z_0}{+x_1}\nno \\
&& \mbox{\rm if }|z_2|>|z_0|>0;\label{l2-2}\\
\dlt{x_1}{-z_1}{+x_0}\dlt{z_2}{x_0}{-x_2}&=&
\dlt{x_0}{z_2}{+x_2}\dlt{x_1}{-z_0}{+x_2}\nno \\
&& \mbox{\rm if }|z_1|>|z_2|>0;\label{l3}\\
\dlt{x_2}{-z_2}{+x_0}\dlt{x_1}{x_2}{-z_0}&=&
\dlt{x_1}{-z_1}{+x_0}\dlt{x_2}{-z_2}{+x_0}\nno \\
&& \mbox{\rm if }|z_2|>|z_0|.\label{l4}
\end{eqnarray}
(Note that the first identity does not require a restricted domain for
$z_1$, $z_2$ and $z_0$, while the others need certain conditions among
the complex numbers $z_i$ in order for the expressions on the
left-hand sides to be well defined, that is, absolutely convergent.
None of the five expressions on the right-hand sides require
restricted domains for absolute convergence.)
\end{lemma}
\pf In this proof we will use additional formal variables $y_0$,
$y_1$, $y_2$, and repeatedly use Remark 2.3.25 in \cite{LL} about
delta function substitution.

First, we have
\begin{eqnarray}\label{proofof8.3}
\dlt{x_1}{x_0}{-y_1}\dlt{x_2}{x_0}{-y_2}&=&
\dlt{x_1}{x_0}{-y_1}\dlt{x_0}{x_2}{+y_2}\nno \\
&=& \dlt{x_1}{x_2+y_2}{-y_1}\dlt{x_0}{x_2}{+y_2}\nno \\
&=& \dlt{x_2}{x_0}{-y_2}\dlt{x_1}{x_2}{-(y_1-y_2)}.
\end{eqnarray}
(Note that the notation $(x_0+y_2-y_1)^{n}$ is unambiguous: it is the
power series expansion in nonnegative powers of $y_1$ and $y_2$.)
Since it is clear that the left-hand side of this identity lies in
\[
{\mathbb C}[y_1, y_2]((x_0^{-1}))[[x_1, x_1^{-1}, x_2, x_2^{-1}]],
\]
one
can substitute any complex numbers $z_1$, $z_2$ for $y_1$, $y_2$,
respectively, and get the identity (\ref{l1}).

For (\ref{l2-1}), we have
\begin{eqnarray*}
\dlt{y_1}{x_0}{-x_1}\dlt{x_2}{x_0}{-y_2}&=&
\dlt{x_0}{y_1}{+x_1}\dlt{x_2}{x_0}{-y_2}\\
& = & \dlt{x_0}{y_1}{+x_1}\dlt{x_2}{y_1+x_1}{-y_2}\\
& = & \dlt{x_0}{y_1}{+x_1}\dlt{x_2}{(y_1-y_2)}{+x_1},
\end{eqnarray*}
and the right-hand side and hence the left-hand side lies in 
\[{\mathbb C}[y_1, y_1^{-1}, (y_1-y_2),
(y_1-y_2)^{-1}][[x_0, x_0^{-1}, x_1,x_2, x_2^{-1}]].\]
Thus if
$|z_1|>|z_2|>0$, so that the binomial expansion of $(z_1-z_2)^{n}$
converges for all $n$, we can substitute $z_1$, $z_2$ for $y_1$,
$y_2$ and obtain (\ref{l2-1}). On the other hand,
\begin{eqnarray*}
\dlt{y_2}{x_0}{-x_2}\dlt{y_0}{x_2}{-x_1}&=&
\dlt{x_0}{y_2}{+x_2}\dlt{x_2}{y_0}{+x_1}\\
&=&\dlt{x_0}{y_2}{+y_0+x_1}\dlt{x_2}{y_0}{+x_1}.
\end{eqnarray*}
It is clear {}from the right-hand side that both sides lie in 
\[{\mathbb
C}[y_0, y_0^{-1}, (y_2+y_0), (y_2+y_0)^{-1}][[x_0, x_0^{-1}, x_1, x_2,
x_2^{-1}]],\]
so if $|z_2|>|z_0|>0$ we can substitute $z_2$, $z_0$
for $y_2$, $y_0$ and obtain (\ref{l2-2}).

To prove (\ref{l3}), we see that
\begin{eqnarray*}
\dlt{x_1}{-y_1}{+x_0}\dlt{y_2}{x_0}{-x_2}&=&
\dlt{x_1}{-y_1}{+x_0}\dlt{x_0}{y_2}{+x_2}\\ &=&
\dlt{x_1}{-y_1}{+y_2+x_2}\dlt{x_0}{y_2}{+x_2},
\end{eqnarray*}
and the right-hand side and hence both sides lie in 
\[{\mathbb C}[y_2, y_2^{-1}, (y_1-y_2), (y_1-y_2)^{-1}][[x_0,
x_0^{-1}, x_1, x_1^{-1}, x_2]],\]
so that when $|z_1|>|z_2|>0$ we can
substitute $z_1$, $z_2$ for $y_1$, $y_2$ and obtain (\ref{l3}).
Finally, we have
\begin{eqnarray*}
\dlt{x_2}{-y_2}{+x_0}\dlt{x_1}{x_2}{-y_0}&=&
\dlt{x_2}{-y_2}{+x_0}\dlt{x_1}{-y_2+x_0}{-y_0}\\
&=&\dlt{x_2}{-y_2}{+x_0}\dlt{x_1}{-(y_2+y_0)}{+x_0},
\end{eqnarray*}
and {}from the right-hand side we see that both sides lie in 
\[{\mathbb
C}[(y_2+y_0), (y_2+y_0)^{-1}, y_2, y_2^{-1}][[x_0, x_1, x_1^{-1}, x_2,
x_2^{-1}]],\]
so that when $|z_2|>|z_0|$ we can substitute $z_2$,
$z_0$ for $y_2$, $y_0$ and obtain the identity (\ref{l4}).  \epfv

If we assume that the convergence condition for intertwining maps in
${\cal C}$ holds and that our generalized modules are objects of
${\cal C}$, in the setting of Section 7, then after pairing with an
element $w'_{(4)} \in W'_4$, we can specialize the formal variables
$y_1$, $y_2$ to complex numbers $z_1$, $z_2$ in (\ref{Y12}) whenever
$|z_1|>|z_2|>0$, and we can specialize $y_2$, $y_0$ to complex numbers
$z_2$, $z_0$ in (\ref{Y34}) whenever $|z_2|>|z_0|>0$, using Lemma
\ref{deltalemma}:

\begin{propo}\label{compositeJacobiforproductsanditerates}
Assume that the convergence condition for intertwining maps in ${\cal
C}$ holds and that the generalized modules entering into (\ref{Y12})
and (\ref{Y34}) are objects of ${\cal C}$.  Continuing to use the
notation of (\ref{Y12}) and (\ref{Y34}), also let $w'_{(4)} \in W'_4$.
Let $z_1$, $z_2$ be complex numbers satisfying $|z_1|>|z_2|>0$.  Then
for a $P(z_1)$-intertwining map $I_1$ of type ${W_4\choose W_1 M_1}$
and a $P(z_2)$-intertwining map $I_2$ of type ${M_1\choose W_2 W_3}$,
the following expressions are absolutely convergent, and the following
formula for the product
\[
I_{1}\circ (1_{W_{1}}\otimes I_{2})
\]
of $I_1$
and $I_2$ holds:
\begin{eqnarray}\label{F12}
\lefteqn{\Bigg\langle w'_{(4)},
\dlt{x_1}{x_0}{-z_1}\dlt{x_2}{x_0}{-z_2}Y_4(v,x_0)
(I_{1}\circ (1_{W_{1}}\otimes I_{2}))
(w_{(1)} \otimes w_{(2)} \otimes w_{(3)})\Bigg\rangle}\nno \\
&&=\Bigg\langle w'_{(4)},\dlt{z_1}{x_0}{-x_1}\dlt{x_2}{x_0}{-z_2}\cdot\nno\\
&&\hspace{5cm}\cdot(I_{1}\circ (1_{W_{1}}\otimes I_{2}))
(Y_1(v,x_1)w_{(1)} \otimes w_{(2)} \otimes
w_{(3)})\Bigg\rangle\nno\\
&&\quad+\Bigg\langle w'_{(4)},\dlt{x_1}{-z_1}{+x_0}\dlt{z_2}{x_0}{-x_2}\cdot\nno\\
&&\hspace{5cm}\cdot(I_{1}\circ (1_{W_{1}}\otimes I_{2}))
(w_{(1)} \otimes Y_2(v,x_2)w_{(2)} \otimes
w_{(3)})\Bigg\rangle\nno\\
&&\quad +\Bigg\langle w'_{(4)},\dlt{x_1}{-z_1}{+x_0}\dlt{x_2}{-z_2}{+x_0}\cdot\nno\\
&&\hspace{5cm}\cdot(I_{1}\circ (1_{W_{1}}\otimes I_{2}))
(w_{(1)} \otimes w_{(2)} \otimes Y_3(v,x_0)
w_{(3)})\Bigg\rangle . \nno\\
\end{eqnarray}
Moreover, let $z_2$, $z_0$ be complex numbers satisfying
$|z_2|>|z_0|>0$. Then for a $P(z_2)$-intertwining map $I^1$ of type
${W_4\choose M_2 W_3}$ and a $P(z_0)$-intertwining map $I^2$ of type
${M_2\choose W_1 W_2}$, the following expressions are absolutely
convergent, and the following formula for the iterate
\[
I^{1}\circ (I^{2}\otimes 1_{W_{3}})
\]
of $I^1$ and $I^2$ holds:
\begin{eqnarray}\label{F34}
\lefteqn{\Bigg\langle w'_{(4)},
\dlt{x_2}{x_0}{-z_2}\dlt{x_1}{x_2}{-z_0}Y_4(v,x_0)
(I^{1}\circ (I^{2}\otimes 1_{W_{3}}))
(w_{(1)} \otimes w_{(2)} \otimes w_{(3)})\Bigg\rangle}\nno\\
&&=\Bigg\langle w'_{(4)},\dlt{z_2}{x_0}{-x_2}\dlt{z_0}{x_2}{-x_1}\cdot\nno\\
&&\hspace{5cm}\cdot(I^{1}\circ (I^{2}\otimes 1_{W_{3}}))
(Y_1(v,x_1)w_{(1)} \otimes w_{(2)} \otimes
w_{(3)})\Bigg\rangle\nno\\
&&\quad +\Bigg\langle w'_{(4)},\dlt{z_2}{x_0}{-x_2}\dlt{x_1}{-z_0}{+x_2}\cdot\nno\\
&&\hspace{5cm}\cdot(I^{1}\circ (I^{2}\otimes 1_{W_{3}}))
(w_{(1)} \otimes Y_2(v,x_2)w_{(2)} \otimes
w_{(3)})\Bigg\rangle\nno \\
&&\quad+\Bigg\langle w'_{(4)},\dlt{x_2}{-z_2}{+x_0}\dlt{x_1}{x_2}{-z_0}\cdot\nno\\
&&\hspace{5cm}\cdot(I^{1}\circ (I^{2}\otimes 1_{W_{3}}))
(w_{(1)} \otimes w_{(2)} \otimes Y_3(v,x_0)
w_{(3)})\Bigg\rangle . \nno\\
\end{eqnarray}
\end{propo}
\pf When $y_1$ and $y_2$ are specialized to $z_1$ and $z_2$,
respectively, the product of the two delta-function expressions on the
left-hand side of (\ref{Y12}) and the three products of pairs of
delta-function expressions on the right-hand side of (\ref{Y12}) all
converge absolutely in the domain $|z_1|>|z_2|>0$, by Lemma
\ref{deltalemma}; note that for the last of the three products of
pairs of delta-function expressions on the right-hand side of
(\ref{Y12}), the convergence is immediate.  Analogously, {}from Lemma
\ref{deltalemma} we see that the corresponding statements also hold
for (\ref{Y34}), when $y_0$ and $y_2$ are specialized to $z_0$ and
$z_2$, respectively, in the domain $|z_1|>|z_2|>0$.  Recalling the
notations (\ref{4itm}), (\ref{4prm}), (\ref{iterateabbreviation}) and
(\ref{productabbreviation}), we see that the result follows {}from the
convergence condition.  \epfv

Considering the ${\mathfrak s}{\mathfrak l}(2)$-action instead of the
$V$-action, by (\ref{log:L(j)b}) we have
\begin{eqnarray}
\lefteqn{L(j)\mathcal{Y}_1(w_{(1)},y_1)\mathcal{Y}_2(w_{(2)},y_2)w_{(3)}}
\nno\\
&&=\sum_{i=0}^{j+1}{j+1\choose i}y_1^i\mathcal{Y}_1(L(j-i)w_{(1)},y_1)
\mathcal{Y}_2(w_{(2)},y_2)w_{(3)}\nno\\ &&\quad+\mathcal{Y}_1(w_{(1)},y_1)L(j)
\mathcal{Y}_2(w_{(2)},y_2)w_{(3)}\nno\\
&&=\sum_{i=0}^{j+1}{j+1\choose i}y_1^i\mathcal{Y}_1(L(j-i)w_{(1)},y_1)
\mathcal{Y}_2(w_{(2)},y_2)w_{(3)}\nno\\ &&\quad+\mathcal{Y}_1(w_{(1)},y_1)
\sum_{k=0}^{j+1}{j+1\choose k}y_2^k{\cal
Y}_2(L(j-k)w_{(2)},y_2)w_{(3)}\nno\\ &&\quad+{\cal
Y}_1(w_{(1)},y_1)\mathcal{Y}_2(w_{(2)},y_2)L(j)w_{(3)}
\end{eqnarray}
for $j=-1, 0$ and $1$.  In the setting of Proposition
\ref{compositeJacobiforproductsanditerates}, if $|z_1|>|z_2|>0$ we can
substitute $z_1$, $z_2$ for $y_1$, $y_2$, respectively, and we obtain,
setting $z_0 = z_1 - z_2$,
\begin{eqnarray}\label{zz:sl2p}
\lefteqn{\langle w'_{(4)},L(j)(I_{1}\circ (1_{W_{1}}\otimes I_{2}))
(w_{(1)} \otimes w_{(2)} \otimes w_{(3)})\rangle}\nno\\
&&=\Bigg\langle w'_{(4)},\sum_{i=0}^{j+1}{j+1\choose i}(z_2+z_0)^i
(I_{1}\circ (1_{W_{1}}\otimes I_{2}))
(L(j-i)w_{(1)} \otimes w_{(2)} \otimes w_{(3)})\Bigg\rangle\nno\\
&&\quad+\Bigg\langle w'_{(4)},\sum_{k=0}^{j+1}{j+1\choose k}z_2^k
(I_{1}\circ (1_{W_{1}}\otimes I_{2}))
(w_{(1)} \otimes L(j-k)w_{(2)} \otimes w_{(3)})\Bigg\rangle\nno\\
&&\quad+\langle w'_{(4)},
(I_{1}\circ (1_{W_{1}}\otimes I_{2}))
(w_{(1)} \otimes w_{(2)} \otimes L(j)w_{(3)})\rangle
\end{eqnarray}
for $j=-1, 0$ and $1$.

On the other hand, by (\ref{log:L(j)b}) we also have
\begin{eqnarray}\label{zz:sl2i0}
\lefteqn{L(j){\cal Y}^1({\cal Y}^2(w_{(1)},y_0)w_{(2)},y_2)w_{(3)}}
\nno\\
&&=\sum_{i=0}^{j+1}{j+1\choose i}y_2^i{\cal Y}^1(L(j-i){\cal Y}^2
(w_{(1)}, y_0)w_{(2)},y_2)w_{(3)}\nno\\ &&\quad+{\cal Y}^1({\cal
Y}^2(w_{(1)},y_0)w_{(2)},y_2)L(j)w_{(3)} \nno\\
&&=\sum_{i=0}^{j+1}{j+1\choose i}y_2^i{\cal Y}^1
\Bigg(\sum_{k=0}^{j-i+1}{j-i+1\choose k}y_0^k{\cal
Y}^2(L(j-i-k)w_{(1)},y_0)w_{(2)} ,y_2\Bigg)w_{(3)}\nno\\
&&\quad+\sum_{i=0}^{j+1}{j+1\choose i}y_2^i{\cal Y}^1( {\cal
Y}^2(w_{(1)}, y_0)L(j-i)w_{(2)} ,y_2)w_{(3)}\nno\\ &&\quad+{\cal
Y}^1({\cal Y}^2(w_{(1)},y_0)w_{(2)},y_2)L(j)w_{(3)}
\end{eqnarray}
for $j=-1, 0$ and $1$, The first term of the right-hand side is
\begin{eqnarray*}
\lefteqn{\sum_{i=0}^{j+1}{j+1\choose i}y_2^i\sum_{k=0}^{j-i+1}
{j-i+1\choose k}y_0^k{\cal Y}^1({\cal Y}^2(L(j-i-k)w_{(1)},y_0)w_{(2)}
,y_2)w_{(3)}}\nno\\
&&=\sum_{t=0}^{j+1}\sum_{k=0}^t{j+1\choose t-k}{j+1-t+k\choose k}
y_2^{t-k}y_0^k{\cal Y}^1({\cal Y}^2(L(j-t)w_{(1)},y_0)w_{(2)}
,y_2)w_{(3)}\nno\\
&&=\sum_{t=0}^{j+1}{j+1\choose t}(y_2+y_0)^t{\cal Y}^1({\cal
Y}^2(L(j-t)w_{(1)},y_0)w_{(2)} ,y_2)w_{(3)},
\end{eqnarray*}
where we have used the identity $\displaystyle{j+1\choose t-k}{j+1-t+k
\choose k}={j+1\choose t}{t\choose k}$ in the last step.  Thus in the
setting of Proposition \ref{compositeJacobiforproductsanditerates}, if
$|z_2|>|z_0|>0$ we can substitute $z_2$, $z_0$ for $y_2$, $y_0$,
respectively, in (\ref{zz:sl2i0}), and we obtain
\begin{eqnarray}\label{zz:sl2i}
\lefteqn{\langle w'_{(4)},L(j)(I^{1}\circ (I^{2}\otimes 1_{W_3}))
(w_{(1)} \otimes w_{(2)} \otimes w_{(3)})\rangle}\nno\\
&&=\Bigg\langle w'_{(4)},\sum_{t=0}^{j+1}{j+1\choose t}(z_2+z_0)^t
(I^{1}\circ (I^{2}\otimes 1_{W_3}))
(L(j-t)w_{(1)} \otimes w_{(2)} \otimes w_{(3)})\Bigg\rangle\nno\\
&&\quad+\Bigg\langle w'_{(4)},\sum_{i=0}^{j+1}{j+1\choose i}z_2^i
(I^{1}\circ (I^{2}\otimes 1_{W_3}))
(w_{(1)} \otimes L(j-i)w_{(2)} \otimes w_{(3)})\Bigg\rangle\nno\\
&&\quad+\langle w'_{(4)},(I^{1}\circ (I^{2}\otimes 1_{W_3}))
(w_{(1)} \otimes w_{(2)} \otimes L(j)w_{(3)})\rangle
\end{eqnarray}
for $j=-1, 0$ and $1$.

Of course, in case $V$ is a conformal vertex algebra, these formulas
follow {}from the earlier computation for the $V$-action (Proposition
\ref{compositeJacobiforproductsanditerates}), by setting $v=\omega$
and taking $\res_{x_1}\res_{x_2}\res_{x_0}x_0^{j+1}$, $j=-1, 0, 1$.

Lemma \ref{deltalemma}, Proposition
\ref{compositeJacobiforproductsanditerates}, (\ref{zz:sl2p}),
(\ref{zz:sl2i}) and Remark \ref{grad-comp-prod-iter} motivate the
following definition, which is analogous to the definition of the
notion of $P(z)$-intertwining map (Definition \ref{im:imdef}):

\begin{defi}\label{Pz1z2intwmap}
{\rm Let $z_0, z_1, z_2\in {\mathbb C}^{\times }$ with $z_0=z_1-z_2$
(so that in particular $z_1\neq z_2$, $z_0\neq z_1$ and $z_0\neq
-z_2$). Let $(W_1,Y_1)$, $(W_2,Y_2)$, $(W_3,Y_3)$ and $(W_4,Y_4)$ be
generalized modules for a M\"obius (or conformal) vertex algebra $V$.
A {\it $P(z_1,z_2)$-intertwining map} is a linear map
\[
F:\, W_1\otimes W_2\otimes W_3\to \overline{W}_4
\]
such that the following conditions
are satisfied: the {\it grading compatibility condition}: For $\beta,
\gamma, \delta \in \tilde A$ and $w_{(1)}\in W_1^{(\beta)}$,
$w_{(2)}\in W_2^{(\gamma)}$, $w_{(3)}\in W_3^{(\delta)}$,
\begin{equation} \label{grad-comp-F}
F(w_{(1)}\otimes 
w_{(2)}\otimes w_{(3)})\in \overline{W_{4}^{(\beta+\gamma+\delta)}};
\end{equation}
the {\em lower truncation condition}: for any elements $w_{(1)}\in
W_1$, $w_{(2)}\in W_2$ and $w_{(3)}\in W_3$, and any $n\in {\mathbb
C}$,
\begin{equation}\label{zz:ltc}
\pi_{n-m}F(w_{(1)}\otimes w_{(2)}\otimes w_{(3)})=0\;\;\mbox{ for }\;m\in
{\mathbb N}\;\mbox{ sufficiently large}
\end{equation}
(which follows {}from (\ref{grad-comp-F}), in view of the grading
restriction condition (\ref{set:dmltc})); the {\em composite Jacobi
identity}:
\begin{eqnarray}\label{zz:Y}
\lefteqn{\dlt{x_1}{x_0}{-z_1}\dlt{x_2}{x_0}{-z_2}Y_4(v,x_0)F(w_{(1)}\otimes
w_{(2)}\otimes w_{(3)})}\nno \\
&& =\dlt{x_0}{z_1}{+x_1}\dlt{x_2}{z_0}{+x_1}F(Y_1(v,x_1)w_{(1)}\otimes
w_{(2)}\otimes w_{(3)})\nno \\
&& \quad+\dlt{x_0}{z_2}{+x_2}\dlt{x_1}{-z_0}{+x_2}F(w_{(1)}\otimes
Y_2(v,x_2)w_{(2)}\otimes w_{(3)})\nno \\
&& \quad+\dlt{x_1}{-z_1}{+x_0}\dlt{x_2}{-z_2}{+x_0}F(w_{(1)}\otimes
w_{(2)}\otimes Y_3(v,x_0)w_{(3)})\nno \\
\end{eqnarray}
for $v\in V$, $w_{(1)}\in W_1$, $w_{(2)}\in W_2$ and $w_{(3)}\in W_3$
(note that all the expressions in the right-hand side of (\ref{zz:Y})
are well defined, that none of the products of delta-function
expressions require restricted domains, and that the left-hand side of
(\ref{zz:Y}) is meaningful because any infinite linear combination of
$v_n$ ($n \in {\mathbb Z}$) of the form $\sum_{n<N}a_nv_n$ ($a_n\in
{\mathbb C}$) acts in a well-defined way on any $F(w_{(1)}\otimes
w_{(2)}\otimes w_{(3)})$, in view of (\ref{zz:ltc})); and the {\em
${\mathfrak s}{\mathfrak l}(2)$-bracket relations}: for any
$w_{(1)}\in W_1$, $w_{(2)}\in W_2$ and $w_{(3)}\in W_3$,
\begin{eqnarray}\label{zz:L}
\lefteqn{L(j)F(w_{(1)}\otimes w_{(2)}\otimes w_{(3)})}\nno\\
&&=\sum_{i=0}^{j+1}{j+1\choose i}z_{1}^{i}F
(L(j-i)w_{(1)}\otimes w_{(2)}\otimes w_{(3)})\nno\\
&&\quad+\sum_{k=0}^{j+1}{j+1\choose k}z_2^kF
(w_{(1)}\otimes L(j-k)w_{(2)}\otimes w_{(3)})\nno\\
&&\quad+F(w_{(1)}\otimes w_{(2)}\otimes L(j)w_{(3)})
\end{eqnarray}
for $j=-1, 0$ and $1$ (again, in case $V$ is a conformal
vertex algebra, this follows {}from (\ref{zz:Y}) by setting $v=\omega$
and taking $\res_{x_1}\res_{x_2}\res_{x_0}x_0^{j+1}$).  }
\end{defi}

\begin{rema}{\rm (cf. Remark \ref{Pintwmaplowerbdd})  If $W_4$ in Definition
\ref{Pz1z2intwmap} is lower bounded, then (\ref{zz:ltc}) can be
strengthened to:
\begin{equation}\label{pinF=0}
\pi_{n}F(w_{(1)}\otimes w_{(2)}\otimes w_{(3)})=0\;\;\mbox{ for }
\;\Re{(n)}\;\mbox{ sufficiently negative.}
\end{equation}}
\end{rema}

We emphasize that every term in (\ref{zz:Y}) and (\ref{zz:L}) in this
definition is purely algebraic; that is, no convergence is involved.

{}From Lemma \ref{deltalemma}, Proposition
\ref{compositeJacobiforproductsanditerates}, (\ref{zz:sl2p}),
(\ref{zz:sl2i}) and Remark \ref{grad-comp-prod-iter}, we have the
following:

\begin{propo}\label{productanditerateareintwmaps}
In the setting of Proposition
\ref{compositeJacobiforproductsanditerates}, for intertwining maps
$I_1$, $I_2$, $I^1$ and $I^2$ as indicated, when $|z_1|>|z_2|>0$,
$I_1\circ (1_{W_1}\otimes I_2)$ is a $P(z_1,z_2)$-intertwining map and
when $|z_2|>|z_0|>0$, $I^1\circ (I^2\otimes 1_{W_3})$ is a
$P(z_2+z_0,z_2)$-intertwining map. \epf
\end{propo}

Now we consider $P(z_1,z_2)$-intertwining maps {}from a ``dual''
viewpoint, and we use this to motivate an analogue $\tau_{P(z_{1},
z_{2})}$ of the action $\tau_{P(z)}$ introduced in Section 5.2.
Fix any $w'_{(4)}\in W'_4$. Then (\ref{zz:Y}) implies:
\begin{eqnarray}
\lefteqn{\Bigg\langle
w'_{(4)},\dlt{x_1}{x_0}{-z_1}\dlt{x_2}{x_0}{-z_2}Y_4(v,x_0)F(w_{(1)}
\otimes w_{(2)}\otimes w_{(3)})\Bigg\rangle} \nno \\
&&=\Bigg\langle
w'_{(4)},\dlt{x_0}{z_1}{+x_1}\dlt{x_2}{z_0}{+x_1}F(Y_1(v,x_1)w_{(1)}
\otimes w_{(2)}\otimes w_{(3)})\Bigg\rangle \nno \\
&& \quad +\Bigg\langle
w'_{(4)},\dlt{x_0}{z_2}{+x_2}\dlt{x_1}{-z_0}{+x_2}F(w_{(1)}\otimes
Y_2(v,x_2)w_{(2)}\otimes w_{(3)})\Bigg\rangle \nno \\
&& \quad +\Bigg\langle
w'_{(4)},\dlt{x_1}{-z_1}{+x_0}\dlt{x_2}{-z_2}{+x_0}F(w_{(1)}\otimes
w_{(2)}\otimes Y_3(v,x_0)w_{(3)})\Bigg\rangle .\nno\\
&& \label{cmpF}
\end{eqnarray}
The left-hand side can be written as
\begin{eqnarray*}
\Bigg\langle
\dlt{x_1}{x_0}{-z_1}\dlt{x_2}{x_0}{-z_2}Y'_4(e^{x_0L(1)}(-x_0^2)^{-L(0)}v,
x_0^{-1})w'_{(4)},
F(w_{(1)}\otimes w_{(2)}\otimes w_{(3)})\Bigg\rangle,
\end{eqnarray*}
and so by replacing $v$ by $(-x_0^2)^{L(0)}e^{-x_0L(1)}v$ and then
replacing $x_0$ by $x_0^{-1}$ in both sides of (\ref{cmpF}) we
see that
\begin{eqnarray}\label{taumot}
\lefteqn{\Bigg\langle \dlt{x_1}{x^{-1}_0}{-z_1}\dlt{x_2}{x^{-1}_0}{-z_2}
Y'_4(v,x_0)w'_{(4)},F(w_{(1)}\otimes w_{(2)}\otimes
w_{(3)})\Bigg\rangle} \nno\\ 
&&=\Bigg\langle
w'_{(4)},\dlti{x_0}{z_1}{+x_1}\dlt{x_2}{z_0}{+x_1}\cdot\nno\\
&&\qquad\qquad
F(Y_1((-x_0^{-2})^{L(0)}e^{-x_0^{-1}L(1)}v,x_1)w_{(1)}\otimes
w_{(2)}\otimes w_{(3)})\Bigg\rangle \nno\\ 
&&+\Bigg\langle
w'_{(4)},\dlti{x_0}{z_2}{+x_2}\dlt{x_1}{-z_0}{+x_2}\cdot\nno\\
&&\qquad\qquad F(w_{(1)}\otimes
Y_2((-x_0^{-2})^{L(0)}e^{-x_0^{-1}L(1)}v,x_2)w_{(2)}\otimes
w_{(3)})\Bigg\rangle \nno\\ 
&&+\Bigg\langle
w'_{(4)},\dlt{x_1}{-z_1}{+x^{-1}_0}\dlt{x_2}{-z_2}{+x^{-1}_0}\cdot\nno\\
&&\qquad\qquad F(w_{(1)}\otimes w_{(2)}\otimes
Y_3((-x_0^{-2})^{L(0)}e^{-x_0^{-1}L(1)}v,x_0^{-1})w_{(3)})\Bigg\rangle.
\end{eqnarray}

Arguing just as in (\ref{y-t-delta})--(\ref{3.19-1}),
we note that in the left-hand side of (\ref{taumot}), the coefficients of
\[
\dlt{x_1}{x^{-1}_0}{-z_1}\dlt{x_2}{x^{-1}_0}{-z_2}Y'_4(v,x_0)
\]
in powers of $x_0$, $x_1$ and $x_2$, for all $v\in V$, span
\[
\tau_{W'_4}(V\otimes \iota_{+}{\mathbb C}[t,t^{-1},(z_1^{-1}-t)^{-1},
(z_2^{-1}-t)^{-1}])
\]
(recall the notation $\tau_W$ {}from (\ref{tauW}), (\ref{tauw}),
(\ref{3.7}) and the notation $\iota_{\pm}$ {}from (\ref{iota+-})).
By analogy with the case of $P(z)$-intertwining maps, we shall
define an action of
\[
V\otimes \iota_{+}{\mathbb C}[t,t^{-1},(z_1^{-1}
-t)^{-1},(z_2^{-1}-t)^{-1}]
\]
on $(W_1\otimes W_2\otimes W_3)^{*}$. We
shall need the following analogue of Lemma \ref{tauP}, where we use the 
notations $Y_{t}$, $T_{z}$ and $o$ introduced in Section 5.1,
and where we recall that $z_{0}=z_{1}-z_{2}$:

\begin{lemma}\label{tauzzlm}
We have
\begin{eqnarray}
\lefteqn{o\bigg(\dlt{x_1}{x_0^{-1}}{-z_1}\dlt{x_2}
{x_0^{-1}}{-z_2}Y_{t}(v,x_0)\bigg)}\nno\\
&&=\dlt{x_1}{x_0^{-1}}{-z_1}\dlt{x_2}
{x_0^{-1}}{-z_2}Y^o_{t}(v,x_0),\label{zztr1}\\
\lefteqn{(\iota_+\circ\iota_-^{-1}\circ o)\bigg(\dlt{x_1}{x_0^{-1}}{-z_1}
\dlt{x_2}{x_0^{-1}}{-z_2}Y_{t}(v,x_0)\bigg)}\nno\\
&&=\dlt{x_1}{-z_1}{+x^{-1}_0}\dlt{x_2}{-z_2}{+x^{-1}_0}Y^o_t(v,x_0),
\label{zztr2}\\
\lefteqn{(\iota_+\circ T_{z_1}\circ\iota_-^{-1}\circ o)\bigg(\dlt{x_1}
{x_0^{-1}}{-z_1}\dlt{x_2}{x_0^{-1}}{-z_2}Y_{t}(v,x_0)\bigg)}\nno\\
&&=\dlti{x_0}{z_1}{+x_1}\dlt{x_2}{z_0}{+x_1}Y_t((-x_0^{-2})^{L(0)}
e^{-x_0^{-1}L(1)}v,x_1),\label{zztr3}\\
\lefteqn{(\iota_+\circ T_{z_2}\circ\iota_-^{-1}\circ o)\bigg(\dlt{x_1}
{x_0^{-1}}{-z_1}\dlt{x_2}{x_0^{-1}}{-z_2}Y_{t}(v,x_0)\bigg)}\nno\\
&&=\dlti{x_0}{z_2}{+x_2}\dlt{x_1}{-z_0}{+x_2}Y_t((-x_0^{-2})^{L(0)}
e^{-x_0^{-1}L(1)}v,x_2).\label{zztr4}
\end{eqnarray}
\end{lemma}
\pf The identity (\ref{zztr1})  immediately follows {}from (\ref{3.38}), and
(\ref{zztr2}) follows {}from (\ref{zztr1}), as in the proof of
(\ref{ztr2}). For (\ref{zztr3}), note that by (\ref{op-y-t-2}),
the coefficient of $x_1^{-m-1}x_2^{-n-1}$ in the right-hand side of
(\ref{zztr1}) is
\begin{eqnarray*}
\lefteqn{(x_0^{-1}-z_1)^m(x_0^{-1}-z_2)^n\bigg(e^{x_0L(1)}(-x_0^{-2})^{L(0)}v\otimes
x_0\delta\bigg(\frac t{x_0^{-1}}\bigg)\bigg)}\nn
&&=(t-z_1)^m(t-z_2)^n\bigg(e^{x_0L(1)}(-x_0^{-2})^{L(0)}v\otimes
x_0\delta\bigg(\frac t{x_0^{-1}}\bigg)\bigg).
\end{eqnarray*}
Acted on by $\iota_+\circ T_{z_1}\circ\iota_-^{-1}$, this becomes
\begin{eqnarray*}
\lefteqn{t^m(z_0+t)^n\bigg(e^{x_0L(1)}(-x_0^{-2})^{L(0)}v\otimes
x_0\delta\bigg(\frac {z_1+t}{x_0^{-1}}\bigg)\bigg)}\nn
&&=x_0\delta\bigg(\frac {z_1+t}{x_0^{-1}}\bigg)(z_0+t)^n\bigg(
e^{x_0L(1)}(-x_0^{-2})^{L(0)}v\otimes t^m\bigg)\nn
&&=x_0\delta\bigg(\frac {z_1+t}{x_0^{-1}}\bigg)(z_0+t)^n\bigg(
(-x_0^{-2})^{L(0)}e^{-x_0^{-1}L(1)}v\otimes t^m\bigg),
\end{eqnarray*}
by formula (5.3.1) in \cite{FHL}, and using (\ref{3.5}), we see that 
this is the coefficient of $x_1^{-m-1}x_2^{-n-1}$ in the right-hand side
of (\ref{zztr3}). The analogous
identity (\ref{zztr4}) is proved similarly. \epfv

Our analogue of Definition \ref{deftau} is:

\begin{defi}\label{tauzzdef}{\rm
Let $z_1, z_2\in {\mathbb C}^{\times}$, $z_1\neq z_2$. We define
a linear action
$\tau_{P(z_1,z_2)}$ of the space
\begin{equation}\label{thez1z2space}
V\otimes \iota _{+}{\mathbb C}[t,t^{-1},(z_1^{-1}-t)^{-1},(z_2^{-1}
-t)^{-1}]
\end{equation}
on $(W_1\otimes W_2\otimes W_3)^{*}$  by
\begin{eqnarray}\label{tauzzdef0}
\lefteqn{(\tau_{P(z_1,z_2)}(\xi)\lambda)(w_{(1)}\otimes w_{(2)}\otimes
w_{(3)})}\nno\\
&&=\lambda(\tau_{W_1}((\iota_+\circ T_{z_1}\circ\iota_-^{-1}\circ
o)\xi)w_{(1)}\otimes w_{(2)}\otimes w_{(3)})\nno\\
&&\quad +\lambda(w_{(1)}\otimes\tau_{W_2}((\iota_+\circ T_{z_2}\circ
\iota_-^{-1}\circ o)\xi)w_{(2)}\otimes w_{(3)})\nno\\
&&\quad +\lambda(w_{(1)}\otimes w_{(2)}\otimes\tau_{W_3}((\iota_+\circ
\iota_-^{-1}\circ o)\xi)w_{(3)})
\end{eqnarray}
for
\[
\xi\in V\otimes \iota _{+}{\mathbb C}[t,t^{-1},(z_1^{-1}-t)^{-1},
(z_2^{-1}-t)^{-1}],
\]
\[
\lambda \in (W_1\otimes W_2\otimes W_3)^{*},
\]
$w_{(1)}\in W_1$, $w_{(2)}\in W_2$ and $w_{(3)}\in W_3$. (The fact that
the right-hand side is in fact well defined follows immediately 
{}from the generating function reformulation of (\ref{tauzzdef0}) given
in (\ref{tauzzgf}) below.)
Denote by $Y'_{P(z_1,z_2)}$ the action of $V\otimes{\mathbb
C}[t,t^{-1}]$ on $(W_1\otimes W_2\otimes W_3)^*$ thus defined, that
is,
\begin{equation}\label{y'-zz}
Y'_{P(z_1,z_2)}(v,x)=\tau_{P(z_1,z_2)}(Y_t(v,x)).
\end{equation}
}
\end{defi}

By Lemma \ref{tauzzlm}, (\ref{3.7}) and (\ref{tauw-yto}), 
we see that (\ref{tauzzdef0}) can be written in terms of
generating functions as
\begin{eqnarray}\label{tauzzgf}
\lefteqn{\left(\tau_{P(z_1,z_2)}\left(\dlt{x_1}{x_0^{-1}}{-z_1}\dlt{x_2}
{x_0^{-1}}{-z_2}Y_{t}(v,x_0)\right)
\lambda \right)(w_{(1)}\otimes w_{(2)}\otimes
w_{(3)})}\nno\\
&=& \dlti{x_0}{z_1}{+x_1}\dlt{x_2}{z_0}{+x_1}\cdot\nno\\
&& \qquad\qquad \lambda (Y_1((-x_0^{-2})^{L(0)}e^{-x_0^{-1}L(1)}v,x_1)
w_{(1)}\otimes w_{(2)}\otimes w_{(3)})\nno\\
&& +\dlti{x_0}{z_2}{+x_2}\dlt{x_1}{-z_0}{+x_2}\cdot\nno\\
&& \qquad\qquad \lambda (w_{(1)}\otimes
Y_2((-x_0^{-2})^{L(0)}e^{-x_0^{-1}L(1)}v,x_2)w_{(2)}\otimes w_{(3)})\nno\\
&& +\dlt{x_1}{-z_1}{+x^{-1}_0}\dlt{x_2}{-z_2}{+x^{-1}_0}\lambda
(w_{(1)}\otimes w_{(2)}\otimes Y_3^o(v,x_0)w_{(3)})\nno\\
\end{eqnarray}
for $v\in V$, $\lambda \in (W_1\otimes W_2\otimes W_3)^{*}$,
$w_{(1)}\in W_1$, $w_{(2)}\in W_2$ and $w_{(3)}\in W_3$; the expansion 
coefficients in $x_{0}$, $x_{1}$ and $x_{2}$ of the left-hand side span 
the space of elements in the left-hand side of (\ref{tauzzdef0}). 
Compare this with the motivating formula (\ref{taumot}).  The
generating function form (\ref{y'-zz}) 
of the action $Y'_{P(z_1,z_2)}$ (\ref{y'-zz}) can be
obtained by taking $\res_{x_1}\res_{x_2}$ of both sides of
(\ref{tauzzgf}).

\begin{rema}{\rm
The action $\tau _{P(z_1,z_2)}$ of
\[
V\otimes \iota _{+}{\mathbb
C}[t,t^{-1},(z_1^{-1}-t)^{-1},(z_2^{-1}-t)^{-1}]
\]
on $(W_1 \otimes W_2
\otimes W_3)^*$,  defined for all $z_1,z_2\in {\mathbb C}^{\times
}$ with $z_1\neq z_2$,  coincides with the action $\tau
^{(1)}_{P(z_1,z_2)}$ when $|z_1|>|z_2|>0$, and coincides with the
action $\tau ^{(2)}_{P(z_1,z_2)}$ when $|z_2|>|z_1-z_2|>0$, where 
$\tau^{(1)}_{P(z_{1}, z_{2})}$ and $\tau^{(2)}_{P(z_{1}, z_{2})}$
are the two actions defined in Section 14 of \cite{tensor4}. The action 
$\tau_{P(z_{1}, z_{2})}$ and the related notion of 
$P(z_{1}, z_{2})$-intertwining map extend the 
corresponding considerations
in \cite{tensor4} in a natural way. }
\end{rema}

\begin{rema}\label{F-intw}
{\rm (cf. Remark \ref{I-intw})
Using the action $\tau_{P(z_1,z_2)}$, we can write
the equality (\ref{taumot}) as
\begin{eqnarray}\label{intw}
\lefteqn{\left(\dlt{x_1}{x^{-1}_0}{-z_1}\dlt{x_2}{x^{-1}_0}{-z_2}
Y'_4(v,x_0)w'_{(4)}\right)\circ F }\nno\\
&&=\tau_{P(z_1,z_2)}\left(\dlt{x_1}{x_0^{-1}}{-z_1}\dlt{x_2}
{x_0^{-1}}{-z_2}Y_{t}(v,x_0)\right)(w'_{(4)}\circ F).
\end{eqnarray}
Furthermore, using the action of $V \otimes
\iota_{+}{\mathbb C}[t,t^{-1}, (z_1^{-1}-t)^{-1},
(z_2^{-1}-t)^{-1}]$ on $W_{4}'$ (recall (\ref{tauW}), 
(\ref{tauw}) and (\ref{3.7})), we can also write
(\ref{intw}) as
\begin{eqnarray}\label{zz:Psi}
\lefteqn{\left(\tau_{W_{4}'}\left(\dlt{x_1}{x^{-1}_0}{-z_1}\dlt{x_2}{x^{-1}_0}
{-z_2} Y_{t}(v,x_0)\right) w'_{(4)}\right)\circ F}\nn
&&=\tau_{P(z_1,z_2)}\left(\dlt{x_1}{x_0^{-1}}{-z_1}\dlt{x_2}
{x_0^{-1}}{-z_2}Y_{t}(v,x_0)\right)(w'_{(4)}\circ F).
\end{eqnarray}}
\end{rema}

As in Section 5, we need to consider gradings by $A$ and $\tilde{A}$.

The space $W_1 \otimes W_2 \otimes W_3$ is naturally
$\tilde{A}$-graded, and this gives us naturally-defined subspaces
$((W_{1}\otimes W_{2}\otimes W_{3})^{*})^{(\beta)}$ for $\beta \in
\tilde{A}$, as in the discussion after Remark \ref{I-intw}.

The space (\ref{thez1z2space}) is naturally $A$-graded, {}from the
$A$-grading on $V$: For $\alpha\in A$,
\begin{equation}
(V\otimes \iota _{+}{\mathbb C}[t,t^{-1},(z_1^{-1}-t)^{-1},(z_2^{-1}
-t)^{-1}])^{(\alpha)}
=V^{(\alpha)}\otimes\iota _{+}{\mathbb C}
[t,t^{-1},(z_1^{-1}-t)^{-1},(z_2^{-1}
-t)^{-1}].
\end{equation}

\begin{defi}\label{3-mod-actioncompatible}
{\rm We call a linear action $\tau$ of 
\[
V\otimes \iota _{+}{\mathbb C}[t,t^{-1},(z_1^{-1}-t)^{-1},(z_2^{-1}
-t)^{-1}]
\]
on $(W_1 \otimes W_2 \otimes W_{3})^*$
{\it $\tilde{A}$-compatible} if 
for $\alpha\in A$, $\beta\in \tilde{A}$,
\[
\xi\in 
(V\otimes \iota _{+}{\mathbb C}[t,t^{-1},(z_1^{-1}-t)^{-1},(z_2^{-1}
-t)^{-1}])^{(\alpha)}
\]
and $\lambda\in ((W_1 \otimes W_2 \otimes W_{3})^{*})^{(\beta)}$,
\[
\tau(\xi)\lambda\in ((W_{1}\otimes W_{2} \otimes W_{3})^{*})^{(\alpha+\beta)}.
\]}
\end{defi}

{}From (\ref{tauzzdef0}) or (\ref{tauzzgf}), we have:

\begin{propo}\label{tauzz-a-comp}
The action  $\tau_{P(z_1,z_2)}$ is $\tilde{A}$-compatible. \epf
\end{propo}

Again as in Section 5, when $V$ is a conformal vertex algebra, we write
\[
Y'_{P(z_1,z_2)}(\omega ,x)=\sum _{n\in {\mathbb Z}}L'_{P(z_1,z_2)}(n)x^{-n-2}.
\]
In this case, by setting $v=\omega$ in (\ref{tauzzgf}) and taking
$\res_{x_0}{x_0}^{j+1}\res_{x_1}\res_{x_2}$ for $j=-1, 0, 1$, we see that
\begin{eqnarray}\label{LP'(j)F}
\lefteqn{(L'_{P(z_1,z_2)}(j)\lambda)(w_{(1)}\otimes w_{(2)}\otimes
w_{(3)})}\nno\\
&&=\lambda\Bigg(\Bigg(\sum_{i=0}^{1-j}{1-j\choose i}z_1^iL(-j-i)\Bigg)
w_{(1)}\otimes
w_{(2)}\otimes w_{(3)} \nno\\
&&\quad\quad\quad\quad +\sum_{i=0}^{1-j}{1-j\choose i}z_2^iw_{(1)}\otimes
L(-j-i)w_{(2)}\otimes w_{(3)}\nno\\
&&\quad\quad\quad\quad +w_{(1)}\otimes w_{(2)}\otimes L(-j)w_{(3)}\Bigg).
\end{eqnarray}
If $V$ is  a M\"obius vertex algebra, we
define the actions $L'_{P(z_1,z_2)}(j)$ on $(W_1\otimes W_2\otimes
W_3)^*$ by (\ref{LP'(j)F}) for $j=-1, 0$ and $1$.  Using these
notations, the ${\mathfrak s}{\mathfrak l}(2)$-bracket relations (\ref{zz:L})
for a $P(z_1, z_2)$-intertwining map $F$ can be written as
\begin{equation}\label{LwF=LwF}
(L'(j)w'_{(4)})\circ F= L'_{P(z_1,z_2)}(j)(w'_{(4)}\circ F)
\end{equation}
for $w'_{(4)}\in W'_4$, $j=-1, 0, 1$ (cf. Remarks \ref{I-intw2} and
\ref{F-intw}).  We have
\[
L'_{P(z_1,z_2)}(j)((W_{1}\otimes W_{2} \otimes W_{3})^{*})^{(\beta)}\subset 
((W_{1}\otimes W_{2}\otimes W_{3})^{*})^{(\beta)}
\]
for $j=-1, 0, 1$ and $\beta\in \tilde{A}$ (cf. Remark
\ref{L'jpreservesbetaspace} and Proposition \ref{tauzz-a-comp}).

For the natural analogue of Proposition \ref{pz} (see Proposition
\ref{zzcor} below), we shall use the following analogues of the
relevant notions in Sections 4 and 5:  A map
\[
F\in \hom(W_1\otimes W_2\otimes W_3, (W_{4}')^{*})
\]
is {\it $\tilde{A}$-compatible} if
\[
F\in \hom(W_1\otimes W_2\otimes W_3, \overline{W}_{4})
\]
and if $F$ satisfies the natural analogue of the condition in
(\ref{IAtildecompat}), as in (\ref{grad-comp-F}).  A map
\[
G \in \hom(W_4',(W_1\otimes W_2\otimes W_3)^{*})
\]
is {\it $\tilde{A}$-compatible} if $G$ satisfies the analogue of 
(\ref{JAtildecompat}).  Then just as in Lemma \ref{IlambdatoJlambda}
and Remark \ref{alternateformoflemma}:

\begin{rema}\label{Atildecompatcorrespondence}{\rm
We have a canonical isomorphism {}from the space of
$\tilde{A}$-compatible linear maps
\[
F:W_{1}\otimes W_{2}\otimes W_{3} \rightarrow \overline{W}_{4}
\]
to the space of $\tilde{A}$-compatible linear maps
\[
G:W'_{4} \rightarrow (W_{1}\otimes W_{2} \otimes W_{3})^{*},
\]
determined by:
\begin{equation}
\langle w'_{(4)}, F(w_{(1)}\otimes w_{(2)}\otimes w_{(3)})\rangle
=G(w'_{(4)})(w_{(1)}\otimes w_{(2)}\otimes w_{(3)})
\end{equation}
for $w_{(1)}\in W_{1}$, $w_{(2)}\in W_{2}$, $w_{(3)}\in W_{3}$  
and $w'_{(4)}\in W'_{4}$,
or equivalently,
\begin{equation}\label{wF=Gw}
w'_{(4)}\circ F = G(w'_{(4)})
\end{equation}
for $w'_{(4)} \in W'_{4}$.
}
\end{rema}

We also have the natural analogues of Definition
\ref{gradingrestrictedmapJ} and Remarks
\ref{Jcompatimpliesgradingrestr} and \ref{Jlowerbounded}:

\begin{defi}{\rm
A map $G\in \hom(W'_4, (W_1\otimes W_2\otimes W_3)^{*})$ is {\em
grading restricted} if for $n\in {\mathbb C}$, $w_{(1)}\in W_1$,
$w_{(2)}\in W_2$ and $w_{(3)}\in W_3$,
\begin{equation}\label{Ggradrestr}
G((W'_4)_{[n-m]})(w_{(1)}\otimes w_{(2)}\otimes w_{(3)})=0\;\;\mbox{ for
}\;m\in {\mathbb N}\;\mbox{ sufficiently large.}
\end{equation}
}
\end{defi}

\begin{rema}
{\rm
If $G \in {\rm Hom}(W'_4,(W_1\otimes W_2 \otimes W_3)^{*})$ is
$\tilde{A}$-compatible, then $G$ is also grading restricted.
}
\end{rema}

\begin{rema}\label{Glowerbounded}
{\rm If in addition $W_4$ (and $W'_4$) are lower bounded, then the
stronger condition
\begin{equation}\label{Glowerbdd}
G((W'_4)_{[n]})(w_{(1)}\otimes w_{(2)}\otimes w_{(3)})=0\;\;\mbox{ for }
\;\Re{(n)}\;\mbox{ sufficiently negative}
\end{equation}
holds.}
\end{rema}

As in Proposition \ref{pz} we now have:

\begin{propo}\label{zzcor}
Let $z_1,z_2\in {\mathbb C}^{\times }$, $z_1\neq z_2$. Let $W_1$,
$W_2$, $W_3$ and $W_4$ be generalized $V$-modules.  Then under the
canonical isomorphism described in Remark
\ref{Atildecompatcorrespondence}, the $P(z_1,z_2)$-intertwining maps $F$
correspond exactly to the (grading restricted) $\tilde{A}$-compatible
maps $G$ that intertwine the actions of
\[
V\otimes \iota _{+}{\mathbb C}[t,t^{-1},(z_1^{-1}-t)^{-1},(z_2^{-1}-t)^{-1}]
\]
and of $L'(j)$ and $L'_{P(z_1,z_2)}(j)$, $j=-1,0,1$, on $W_4'$ and on
$(W_1\otimes W_2\otimes W_3)^*$.  If $W_4$ is lower bounded, we may
replace the grading restrictions by (\ref{pinF=0}) and
(\ref{Glowerbdd}).
\end{propo}
\pf By (\ref{wF=Gw}), Remark \ref{F-intw} asserts that (\ref{taumot}),
or equivalently, (\ref{zz:Y}), is equivalent to the condition
\begin{eqnarray}\label{Gtau=tauG}
\lefteqn{G\left(\tau_{W_{4}'}\left(\dlt{x_1}{x^{-1}_0}{-z_1}\dlt{x_2}{x^{-1}_0}
{-z_2} Y_{t}(v,x_0)\right) w'_{(4)}\right)}\nn
&&=\tau_{P(z_1,z_2)}\left(\dlt{x_1}{x_0^{-1}}{-z_1}\dlt{x_2}
{x_0^{-1}}{-z_2}Y_{t}(v,x_0)\right)G(w'_{(4)}),
\end{eqnarray}
that is, the condition that $G$ intertwines the actions of 
\[
V\otimes \iota _{+}{\mathbb C}[t,t^{-1},(z_1^{-1}-t)^{-1},(z_2^{-1}-t)^{-1}]
\]
on $W_4'$ and on $(W_1\otimes W_2\otimes W_3)^*$.  Analogously, {}from 
(\ref{LwF=LwF}) we see that (\ref{zz:L}) is equivalent to the
condition
\begin{equation}
G(L'(j)w'_{(4)})= L'_{P(z_1,z_2)}(j)G(w'_{(4)})
\end{equation}
for $j=-1, 0, 1$, that is, the condition that $G$ intertwines the
actions of $L'(j)$ and $L'_{P(z_1,z_2)}(j)$.
\epfv

Let $W_{1}$, $W_{2}$ and $W_{3}$ be generalized $V$-modules.  By
analogy with (\ref{W1W2_[C]^Atilde}) and (\ref{W1W2_(C)^Atilde}), we
have the spaces
\begin{equation}\label{3-mod-2-gradings}
((W_1\otimes W_2\otimes W_3)^{*})_{[\C]}^{(\tilde{A})}
=\coprod_{n\in \C}\coprod_{\beta\in \tilde{A}}
((W_1\otimes W_2\otimes W_3)^{*})_{[n]}^{(\beta)}\subset
(W_1\otimes W_2\otimes W_3)^{*}
\end{equation}
and 
\begin{equation}\label{3-mod-2-s-gradings}
((W_1\otimes W_2\otimes W_3)^{*})_{(\C)}^{(\tilde{A})}
=\coprod_{n\in \C}\coprod_{\beta\in \tilde{A}}
((W_1\otimes W_2\otimes W_3)^{*})_{(n)}^{(\beta)}\subset
(W_1\otimes W_2\otimes W_3)^{*},
\end{equation}
defined by means of the operator $L'_{P(z_1,z_2)}(0)$.  Each space
\begin{equation}\label{W1W2W3beta}
((W_{1}\otimes W_{2}\otimes W_{3})^{*})^{(\beta)}
\end{equation}
is defined by analogy with (\ref{W1W2beta}).

Again by analogy with the situation in Section 5, consider the
following conditions for elements
\[
\lambda \in (W_1\otimes W_2\otimes W_3)^*:
\]

\begin{description}
\item{\bf The $P(z_1, z_2)$-compatibility condition}

(a) The \emph{$P(z_1, z_2)$-lower truncation condition}: For all $v\in
V$, the formal Laurent series $Y'_{P(z_1,z_2)}(v,x)\lambda$ involves
only finitely many negative powers of $x$.

(b) The following formula holds for all $v\in V$:
\begin{eqnarray}\label{zz:cpb}
\lefteqn{\tau_{P(z_1,z_2)}\bigg(\dlt{x_1}{x_0^{-1}}{-z_1}\dlt{x_2}{x_0^{-1}}
{-z_2}Y_{t}(v,x_0)\bigg)\lambda}\nno\\
&&=\dlt{x_1}{x_0^{-1}}{-z_1}\dlt{x_2}{x_0^{-1}}{-z_2}Y'_{P(z_1,z_2)}
(v,x_0)\lambda.
\end{eqnarray}
(Note that the two sides of (\ref{zz:cpb}) are not {\it a priori}
equal for general $\lambda\in (W_1\otimes W_2\otimes W_3)^{*}$.
Condition (a) implies that the right-hand side in Condition (b) is
well defined.)

\item{\bf The $P(z_1,z_2)$-local grading restriction condition}

(a) The {\em $P(z_1,z_2)$-grading condition}: There exists a doubly
graded subspace of the space (\ref{3-mod-2-gradings}) containing
$\lambda$ and stable under the component operators
$\tau_{P(z_1,z_2)}(v\otimes t^{m})$ of the operators
$Y'_{P(z_1,z_2)}(v,x)$ for $v\in V$, $m\in {\mathbb Z}$, and under the
operators $L'_{P(z_1,z_2)}(-1)$, $L'_{P(z_1,z_2)}(0)$ and
$L'_{P(z_1,z_2)}(1)$. In particular, $\lambda$ is a (finite) sum of
generalized eigenvectors for $L'_{P(z_1,z_2)}(0)$ that are also
homogeneous with respect to $\tilde A$.

(b) Let $W_{\lambda; P(z_{1}, z_{2})}$ be the smallest doubly graded
(or equivalently, $\tilde A$-graded) subspace of the space
(\ref{3-mod-2-gradings}) containing $\lambda$ and stable under the
component operators $\tau_{P(z_1,z_2)}(v\otimes t^{m})$ of the
operators $Y'_{P(z_1,z_2)}(v,x)$ for $v\in V$, $m\in {\mathbb Z}$, and
under the operators $L'_{P(z_1,z_2)}(-1)$, $L'_{P(z_1,z_2)}(0)$ and
$L'_{P(z_1,z_2)}(1)$ (the existence being guaranteed by Condition
(a)).  Then $W_{\lambda; P(z_{1}, z_{2})}$ has the properties
\begin{eqnarray}
&\dim(W_{\lambda; P(z_{1}, z_{2})})^{(\beta)}_{[n]}<\infty,&\\
&(W_{\lambda; P(z_{1}, z_{2})})^{(\beta)}_{[n+k]}=0\;\;
\mbox{ for }\;k\in {\mathbb Z}
\;\mbox{ sufficiently negative,}&
\end{eqnarray}
for any $n\in {\mathbb C}$ and $\beta\in \tilde A$, where the
subscripts denote the ${\mathbb C}$-grading by (generalized)
$L'_{P(z_1,z_2)}(0)$-eigenvalues and the superscripts denote the
$\tilde A$-grading.

\item{\bf The $L(0)$-semisimple $P(z_{1}, z_{2})$-local 
grading restriction condition}

(a) The {\em $L(0)$-semisimple $P(z_1,z_2)$-grading condition}: There
exists a doubly graded subspace of the space
(\ref{3-mod-2-s-gradings}) containing $\lambda$ and stable under the
component operators $\tau_{P(z_1,z_2)}(v\otimes t^{m})$ of the
operators $Y'_{P(z_1,z_2)}(v,x)$ for $v\in V$, $m\in {\mathbb Z}$, and
under the operators $L'_{P(z_1,z_2)}(-1)$, $L'_{P(z_1,z_2)}(0)$ and
$L'_{P(z_1,z_2)}(1)$. In particular, $\lambda$ is a (finite) sum of
eigenvectors for $L'_{P(z_1,z_2)}(0)$ that are also homogeneous with
respect to $\tilde A$.

(b) Consider $W_{\lambda; P(z_{1}, z_{2})}$ as above, which in this
case is in fact the smallest doubly graded subspace of the space
(\ref{3-mod-2-s-gradings}) containing $\lambda$ and stable under the
component operators $\tau_{P(z_1,z_2)}(v\otimes t^{m})$ of the
operators $Y'_{P(z_1,z_2)}(v,x)$ for $v\in V$, $m\in {\mathbb Z}$, and
under the operators $L'_{P(z_1,z_2)}(-1)$, $L'_{P(z_1,z_2)}(0)$ and
$L'_{P(z_1,z_2)}(1)$.  Then $W_{\lambda; P(z_{1}, z_{2})}$ has the
properties
\begin{eqnarray}
&\dim(W_{\lambda; P(z_{1}, z_{2})})^{(\beta)}_{(n)}
<\infty,&\label{zz-semi-lgrc1}\\
&(W_{\lambda; P(z_{1}, z_{2})})^{(\beta)}_{(n+k)}=0\;\;
\mbox{ for }\;k\in {\mathbb Z}
\;\mbox{ sufficiently negative},&\label{zz-semi-lgrc2}
\end{eqnarray}
for any $n\in {\mathbb C}$ and $\beta\in \tilde A$, where the
subscripts denote the ${\mathbb C}$-grading by
$L'_{P(z_1,z_2)}(0)$-eigenvalues and the superscripts denote the
$\tilde A$-grading.

\end{description}

Then we have the following, by analogy with the comments preceding the
statement of the $P(z)$-compatibility condition (recall
(\ref{5.18-p})) and the $P(z)$-local grading restriction conditions:

\begin{propo}\label{8.12}
Suppose that $G\in \hom (W'_4, (W_1\otimes W_2\otimes W_3)^*)$
corresponds to a $P(z_1, z_2)$-intertwining map as in Proposition
\ref{zzcor}. Then for any $w'_{(4)}\in W'_4$, $G(w'_{(4)})$ satisfies
the $P(z_1, z_2)$-compatibility condition
and the $P(z_1,z_2)$-local grading restriction condition.
If $W_{4}$ is an ordinary $V$-module, then $G(w'_{(4)})$
satisfies the $L(0)$-semisimple $P(z_1,z_2)$-local grading 
restriction condition.
\end{propo}
\pf For any $w'_{(4)}\in W'_4$, the fact that $G(w'_{(4)})$ satisfies
the $P(z_1, z_2)$-compatibility condition follows {}from
(\ref{Gtau=tauG}), just as in (\ref{5.18-p}).  Since $G$ in particular
intertwines the actions of $V\otimes{\mathbb C}[t, t^{-1}]$ and of the
$L(j)$-operators and is $\tilde A$-compatible, $G(W'_4)$ is a
generalized $V$-module and thus $G(w'_{(4)})$ satisfies the
$P(z_1,z_2)$-local grading restriction condition, and if $W_{4}$ is an
ordinary $V$-module, then $G(W'_4)$ must also be an ordinary
$V$-module and thus $G(w'_{(4)})$ satisfies the $L(0)$-semisimple
$P(z_1,z_2)$-local grading restriction condition, just as in the
comments preceding the statement of the $P(z)$-local grading
restriction conditions. \epf

\begin{rema}\label{consequenceofPz1z2compat}
{\rm In the next section we will use the following: Assume the
$P(z_1,z_2)$-compatibility condition.  By (\ref{l1}) (a ``purely
algebraic'' identity, involving no convergence issues), (\ref{zz:cpb})
can be written as
\begin{eqnarray}\label{alternatecompat}
\lefteqn{\tau_{P(z_1,z_2)}\bigg(\dlt{x_1}{x_2}{-z_0}\dlt{x_2}{x_0^{-1}}
{-z_2}Y_{t}(v,x_0)\bigg)\lambda}\nn
&&=\dlt{x_1}{x_2}{-z_0}\dlt{x_2}{x_0^{-1}}{-z_2}Y'_{P(z_1,z_2)}
(v,x_0)\lambda\nn
&&=\dlt{x_1}{x_2}{-z_0}\biggl(\dlt{x_2}{x_0^{-1}}{-z_2}Y'_{P(z_1,z_2)}
(v,x_0)\lambda\biggr),
\end{eqnarray}
and all of the indicated products exist; the definition (\ref{3.5}) of
$Y_{t}(v,x_0)$ makes it clear that the triple product in parentheses
on the left-hand side exists, and the simplest way to see that the
triple product in the middle expression exists is to repeat the proof
(\ref{proofof8.3}) of (\ref{l1}) (with the formal variables $y_1$ and
$y_2$), multiplying each step by $Y'_{P(z_1,z_2)}(v,x_0)\lambda$,
whose powers of $x_0$ are truncated from below.  We can take ${\rm
Res}_{x_{1}}$ of (\ref{alternatecompat}) to obtain
\begin{eqnarray}\label{resofconsequence}
\tau_{P(z_1,z_2)}\bigg(\dlt{x_2}{x_0^{-1}}{-z_2}Y_{t}(v,x_0)\bigg)\lambda
=\dlt{x_2}{x_0^{-1}}{-z_2}Y'_{P(z_1,z_2)} (v,x_0)\lambda,
\end{eqnarray}
which is reminiscent of the $P(z)$-compatibility condition
(\ref{cpb}) for $z=z_2$.  Now we can multiply both sides by
$\displaystyle\dlt{x_1}{x_2}{-z_0}$, giving
\begin{eqnarray*}
\lefteqn{\dlt{x_1}{x_2}{-z_0}\tau_{P(z_1,z_2)}\bigg(\dlt{x_2}{x_0^{-1}}
{-z_2}Y_{t}(v,x_0)\bigg)\lambda}\nn
&&=\dlt{x_1}{x_2}{-z_0}\dlt{x_2}{x_0^{-1}}{-z_2}Y'_{P(z_1,z_2)}
(v,x_0)\lambda,
\end{eqnarray*}
and as we have seen, these products exist.  Thus by (\ref{zz:cpb}) and
(\ref{alternatecompat}),
\begin{eqnarray}\label{consequenceofPz1z2compatformula}
\lefteqn{\tau_{P(z_1,z_2)}\bigg(\dlt{x_1}{x_0^{-1}}{-z_1}\dlt{x_2}{x_0^{-1}}
{-z_2}Y_{t}(v,x_0)\bigg)\lambda}\nn
&&=\dlt{x_1}{x_2}{-z_0}\tau_{P(z_1,z_2)}\bigg(\dlt{x_2}{x_0^{-1}}
{-z_2}Y_{t}(v,x_0)\bigg)\lambda.
\end{eqnarray}}
\end{rema}

Under the assumption that tensor products exist, we can replace
products and iterates of intertwining maps by corresponding products
and iterates for which the intermediate module is a tensor product, and in a 
unique way:

\begin{propo}\label{intermediate}
Assume that the convergence condition for intertwining maps in ${\cal
C}$ holds.  Let $W_1$, $W_2$, $W_3$, $W_4$ and $M_1$ be objects of
${\cal C}$ and let $z_1, z_2 \in {\mathbb C}$ such that
$|z_{1}|>|z_{2}|>0$.  Let $I_1\in
\mathcal{M}[P(z_{1})]_{W_{1}M_{1}}^{W_{4}}$ and $I_2\in
\mathcal{M}[P(z_{2})]_{W_{2}W_{3}}^{M_{1}}$, and assume that
$W_{2}\boxtimes_{P(z_{2})}W_{3}$ exists (in ${\cal C}$).  Then there
exists a unique
\[
\widetilde{I}_{1}\in
\mathcal{M}[P(z_{1})]_{W_{1}\;(W_{2}\boxtimes_{P(z_{2})}W_{3})}^{W_{4}}
\]
such that
\[
I_{1}\circ (1_{W_{1}}\otimes I_{2})
=\widetilde{I}_{1}\circ (1_{W_{1}}\otimes
\boxtimes_{P(z_{2})}).
\]
Analogously, let $W_1$, $W_2$, $W_3$, $W_4$ and $M_2$ be objects of
${\cal C}$, and let $z_2,z_0 \in {\mathbb C}$ such that
$|z_{2}|>|z_{0}|>0$.  Let $I^1\in
\mathcal{M}[P(z_{2})]_{M_{2}W_{3}}^{W_{4}}$ and $I^2\in
\mathcal{M}[P(z_{0})]_{W_{1}W_{2}}^{M_{2}}$, and assume that
$W_{1}\boxtimes_{P(z_{0})}W_{2}$ exists.  Then there exists a unique
\[
\widetilde{I}^{1}\in
\mathcal{M}[P(z_{2})]_{(W_{1}\boxtimes_{P(z_{0})}W_{2})\;W_{3}}^{W_{4}}
\]
such that
\[
I^{1}\circ (I^{2}\otimes 1_{W_{3}})
=\widetilde{I}^{1}\circ (\boxtimes_{P(z_{0})} \otimes
1_{W_{3}}).
\]
\end{propo}
\pf
We prove only the first part; the second part is proved 
analogously.

By Proposition \ref{pz-iso}, $I_{2}$ corresponds naturally to
an element $\eta$ of $\hom(W_{2}\boxtimes_{P(z_{2})}W_{3}, 
M_{1})$ such that $I_{2}=\overline{\eta}\circ \boxtimes_{P(z_{2})}$. 
Let
\[
\widetilde{I}_{1}=I_{1}\circ (1_{W_{1}}\otimes \eta).
\]
Then $\widetilde{I}_{1}$ is a $P(z_{1})$-intertwining map 
of type ${W_{4}\choose W_{1}\; (W_{2}\boxtimes_{P(z_{2})}W_{3})}$
and we have
\begin{eqnarray*}
I_{1}\circ  (1_{W_{1}}\otimes I_{2})
&=&I_{1}\circ  (1_{W_{1}}\otimes
(\overline{\eta}\circ \boxtimes_{P(z_{2})}))\nn
&=&(I_{1}\circ (1_{W_{1}}\otimes \eta))
\circ (1_{W_{1}}\otimes \boxtimes_{P(z_{2})})\nn
&=&\widetilde{I}_{1}\circ (1_{W_{1}}\otimes \boxtimes_{P(z_{2})}),
\end{eqnarray*}
where these expressions are understood in the sense of Definition
\ref{productanditerateexisting}. 

The equality
\[
\langle w'_{(4)}, 
I_1(w_{(1)}\otimes I_2(w_{(2)}\otimes w_{(3)}))\rangle
=\langle w'_{(4)}, \widetilde{I}_{1}(w_{(1)}\otimes 
(w_{(2)}\boxtimes_{P(z_{2})} w_{(3)}))\rangle
\]
for all $w_{(j)}\in W_{j}$ and $w'_{(4)}\in W'_{4}$ determines the
$P(z_1)$-intertwining map $\widetilde{I}_{1}$ uniquely.  Indeed, 
By Proposition
\ref{prod=0=>comp=0}, this assertion uniquely determines 
\[
\langle w'_{(4)}, \widetilde{I}_{1}(w_{(1)}\otimes 
\pi_{n}(w_{(2)}\boxtimes_{P(z_{2})} w_{(3)}))\rangle
\]
for all $n\in \R$ and for all homogeneous vectors and
hence for all vectors, and since the components
$\pi_{n}(w_{(2)}\boxtimes_{P(z_2)} w_{(3)})$ span $W_2
\boxtimes_{P(z_2)} W_3$ by Proposition \ref{span}, $\widetilde{I}_{1}$ is uniquely
determined.
\epfv

{}From Proposition \ref{im:correspond}, in which we take $p=0$, we
obtain the corresponding result for logarithmic intertwining
operators:

\begin{corol}\label{intermediate2}
Under the assumptions of Proposition \ref{intermediate}, let
$\Y_1\in \mathcal{V}_{W_{1}M_{1}}^{W_{4}}$ and
$\Y_2\in \mathcal{V}_{W_{2}W_{3}}^{M_{1}}$.  Then 
there exists a unique
\[
\widetilde{\Y}_{1}\in 
\mathcal{V}_{W_{1}\;(W_{2}\boxtimes_{P(z_{2})}W_{3})}^{W_{4}}
\]
such that for $w_{(1)}\in W_{1}$, $w_{(2)}\in W_{2}$, $w_{(3)}\in W_{3}$
and $w_{(4)}'\in W'_{4}$,
\[
\langle w_{(4)}', \Y_{1}(w_{(1)}, z_{1})\Y_{2}(w_{(2)}, z_{2})w_{(3)}
\rangle
=\langle w_{(4)}', 
\widetilde{\Y}_{1}(w_{(1)}, z_{1})\Y_{\boxtimes_{P(z_{2})}, 0}
(w_{(2)}, z_{2})w_{(3)}\rangle
\]
(recall (\ref{log:IYp}), (\ref{recover}) and
(\ref{productabbreviation})).  Analogously, let $\Y^1\in
\mathcal{V}_{M_{2}W_{3}}^{W_{4}}$ and $\Y^2\in
\mathcal{V}_{W_{1}W_{2}}^{M_{2}}$.  Then there exists a unique
\[
\widetilde{\Y}^{1}\in
\mathcal{V}_{(W_{1}\boxtimes_{P(z_{0})}W_{2})\;W_{3}}^{W_{4}}
\]
such that for $w_{(1)}\in W_{1}$, $w_{(2)}\in W_{2}$, $w_{(3)}\in W_{3}$
and $w_{(4)}'\in W'_{4}$,
\[
\langle w_{(4)}', \Y^{1}(\Y^{2}(w_{(1)}, z_{0})w_{(2)}, z_{2})w_{(3)}\rangle
=\langle w_{(4)}', \widetilde{\Y}^{1}(\Y_{\boxtimes_{P(z_{0})}, 0}
(w_{(1)}, z_{0})w_{(2)}, z_{2})w_{(3)}\rangle.
\]
(recall (\ref{iterateabbreviation})).
\epf
\end{corol}

\begin{rema}\label{factor-thr}
{\rm The first half of Proposition \ref{intermediate} in fact states
that the product of $I_{1}$ and $I_{2}$ can be rewritten as a new
product of intertwining maps such that the intermediate object of the
new product is the tensor product generalized module
$W_{2}\boxtimes_{P(z_{2})}W_{3}$ and the $P(z_{2})$-intertwining map
is $\boxtimes_{P(z_{2})}$.  The second half of the proposition can be
stated analogously, for iterates of intertwining maps.  Corollary
\ref{intermediate2} states that a product or an iterate of logarithmic
intertwining operators, evaluated at suitable points, can be expressed
as a new product or iterate for which the intermediate object is the
relevant tensor product and the second intertwining operator
corresponds to the intertwining map defining the tensor product.  Thus
these results can be viewed as saying that the product of $I_{1}$ and
$I_{2}$, or of $\Y_{1}$ and $\Y_{2}$, uniquely ``factors through'' 
$W_{2}\boxtimes_{P(z_{2})}W_{3}$
and that the iterate of $I^{1}$ and
$I^{2}$, or of $\Y^{1}$ and $\Y^{2}$, uniquely ``factors through''
$W_{1}\boxtimes_{P(z_{0})}W_{2}$. }
\end{rema}


\bigskip

\noindent {\small \sc Department of Mathematics, Rutgers University,
Piscataway, NJ 08854 (permanent address)}

\noindent {\it and}

\noindent {\small \sc Beijing International Center for Mathematical Research,
Peking University, Beijing, China}

\noindent {\em E-mail address}: yzhuang@math.rutgers.edu

\vspace{1em}

\noindent {\small \sc Department of Mathematics, Rutgers University,
Piscataway, NJ 08854}

\noindent {\em E-mail address}: lepowsky@math.rutgers.edu

\vspace{1em}

\noindent {\small \sc Department of Mathematics, Rutgers University,
Piscataway, NJ 08854}

\noindent {\em E-mail address}: linzhang@math.rutgers.edu


\begin{thebibliography}{FGST2}

\bibitem[FHL]{FHL}
I.~B. Frenkel, Y.-Z. Huang and J.~Lepowsky,
On axiomatic approaches to vertex operator algebras and modules,
preprint, 1989;
{\em Memoirs Amer. Math. Soc.} {\bf 104}, 1993.

\bibitem[H1]{H0}
Y.-Z. Huang, On the geometric interpretation of vertex operator algebras, 
Ph.D. thesis, Rutgers University, 1990.

\bibitem[H2]{tensor4}
Y.-Z. Huang, A theory of tensor products for module categories for a
vertex operator algebra, IV, {\em J. Pure Appl. Alg.} 100 (1995)
173--216.

\bibitem[H3]{H1}
Y.-Z. Huang, {\em Two-dimensional Conformal Geometry and Vertex
Operator Algebras}, Progress in Math., Vol. 148, Birkh\"{a}user, Boston, 
1997.

\bibitem[HLZ1]{HLZ1} Y.-Z.~Huang, J.~Lepowsky and L.~Zhang, Logarithmic
tensor category theory for generalized modules for a conformal
vertex algebra, I: Introduction and strongly graded
algebras and their generalized modules, to appear.

\bibitem[HLZ2]{HLZ2} Y.-Z.~Huang, J.~Lepowsky and L.~Zhang, Logarithmic
tensor category theory, II: Logarithmic formal calculus
and properties of logarithmic intertwining operators, to appear.

\bibitem[HLZ3]{HLZ3} Y.-Z.~Huang, J.~Lepowsky and L.~Zhang, Logarithmic
tensor category theory, III: Intertwining maps and tensor
product bifunctors, to appear.

\bibitem[HLZ4]{HLZ4} Y.-Z.~Huang, J.~Lepowsky and L.~Zhang, Logarithmic
tensor category theory, IV: Constructions of tensor
product bifunctors and the compatibility conditions, to appear.

\bibitem[HLZ5]{HLZ6} Y.-Z.~Huang, J.~Lepowsky and L.~Zhang, Logarithmic
tensor category theory, VI: Expansion condition, associativity of logarithmic
intertwining operators, and the associativity isomorphisms, to appear.

\bibitem[HLZ6]{HLZ7} Y.-Z.~Huang, J.~Lepowsky and L.~Zhang, Logarithmic
tensor category theory, VII: Convergence and extension
properties and applications to expansion for intertwining
maps, to appear.

\bibitem[HLZ7]{HLZ8} Y.-Z.~Huang, J.~Lepowsky and L.~Zhang, Logarithmic
tensor category theory, VIII: Braided tensor category
structure on categories of generalized modules for a
conformal vertex algebra, to appear.

\bibitem[LL]{LL}
J. Lepowsky and H. Li, {\em Introduction to Vertex Operator Algebras
and Their Representations}, Progress in Math., Vol. 227, Birkh\"auser,
Boston, 2003.

\end{thebibliography}
\end{document}